\documentclass[11pt]{article}
\usepackage{amssymb,amsmath,amsthm}
\usepackage{amscd}
\usepackage{mathrsfs}
\allowdisplaybreaks[1]
\newcommand{\C}{\mathbb{C}}
\newcommand{\N}{\mathbb{N}}
\newcommand{\R}{\mathbb{R}}

\newcommand{\dd}{\,{\rm d}}

\numberwithin{equation}{section}
\newtheorem{thm}{Theorem}[section]
\newtheorem{df}[thm]{Definition}
\newtheorem{prop}[thm]{Proposition}
\newtheorem{lem}[thm]{Lemma}
\newtheorem{rem}[thm]{Remark}
\newtheorem{cor}[thm]{Corollary}

\begin{document}

\title{On domain of Poisson operators and factorization\\
for divergence form elliptic operators}

\date{}

\author{
\null\\
Yasunori Maekawa\\    
Mathematical Institute, Tohoku University\\          
6-3 Aoba, Aramaki, Aoba, Sendai 980-8578, Japan\\
{\tt maekawa@math.tohoku.ac.jp }
\and
\\
Hideyuki Miura \\
Department of Mathematics,  Graduate School of Science,  Osaka University\\
1-1 Machikaneyama, Toyonaka, Osaka 560-0043, Japan\\
{\tt miura@math.sci.osaka-u.ac.jp}}

\maketitle

\begin{center}
{\bf Abstract}
\end{center}
\vspace{-2mm}
We consider second order uniformly elliptic operators of divergence form in $\R^{d+1}$ whose coefficients are independent of one variable. For such a class of operators we establish a factorization into a product of first order operators related with Poisson operators and Dirichlet-Neumann maps. Consequently, we obtain a solution formula for the inhomogeneous elliptic boundary value problem in the half space, which is useful to show the existence of solutions in a wider class of inhomogeneous data. We also establish $L^2$ solvability of boundary value problems for a new class of the elliptic operators.

\vspace{0.5cm}



\tableofcontents

\section{Introduction}\label{sec.intro}


We consider the second order elliptic operator of divergence form in $\R^{d+1}= \{ (x,t) \in \R^d \times \R \}$, 
\begin{equation}
\mathcal{A} = -\nabla \cdot A \nabla, ~~~~~~~~~~~~~~A = A (x) = \big (a_{i,j} (x) \big ) _{1\leq i,j\leq d+1}\,.\label{def.calA}
\end{equation}
Here $d\in \N$, $\nabla = (\nabla_x,\partial_t)^\top$ with $\nabla_x=(\partial_1,\cdots,\partial_d)^\top$,  and each $a_{i,j}$ is always assumed to be $t$-independent in this paper. The differential operator with $t$-independent coefficients like \eqref{def.calA} has been studied for a long time. 
One of the typical assumptions of $A$ is that $a_{i,j}$ are complex-valued measurable functions satisfying  the ellipticity condition
\begin{align}
{\rm Re} \langle A (x) \eta, \eta \rangle \geq \nu_1 |\eta|^2, ~~~~~~~~~~ | \langle A (x) \eta,\zeta\rangle | \leq \nu_2 |\eta| |\zeta|\label{ellipticity}
\end{align}
for all $\eta,\zeta\in \C^{d+1}$ and for some constants $\nu_1,\nu_2$ with $0<\nu_1\leq \nu_2 <\infty$. Here $\langle \cdot,\cdot \rangle$ denote the inner product of $\C^{d+1}$, i.e., $\langle \eta, \zeta\rangle = \sum_{j=1}^{d+1}\eta_j \bar{\zeta}_j$ for $\eta,\zeta\in \C^{d+1}$. The condition 
\eqref{ellipticity} is always assumed throughout this paper. 
For later use we set 
\begin{align*}
A' =(a_{i,j})_{1\leq i,j\leq d},~~~b=a_{d+1, d+1}, ~~~{\bf r_1} = ( a_{1, d+1},\cdots , a_{d, d+1} )^\top, ~~~{\bf r_2}  = ( a_{d+1,1} ,\cdots , a_{d+1,d})^\top.
\end{align*}
We call ${\bf r_1}$ and ${\bf r_2}$ the off-block vectors of $A$. Let us denote by $D_H (T)$ the domain of a linear operator $T$ in a Banach space $H$. Under the condition \eqref{ellipticity} the standard theory of sesquilinear forms  gives a realization of $\mathcal{A}$ in $L^2 (\R^{d+1})$, denoted again by $\mathcal{A}$, such as
\begin{align}
D_{L^2}(\mathcal{A}) & = \big \{ u\in H^1 (\R^{d+1})~|~{\rm there ~is}~ F \in L^2 (\R^{d+1})~{\rm such ~that}~\nonumber \\
& ~~~~~~~~~~~~~~~~~~ \langle A\nabla u, \nabla v\rangle _{L^2(\R^{d+1})} = \langle F, v\rangle _{L^2 (\R^{d+1})}~{\rm for~all}~v\in H^1 (\R^{d+1})\big \},\label{realization.A}\\
 \mathcal{A} u & = F ~~~~~ {\rm for} ~~ u\in D_{L^2}(\mathcal{A}).\nonumber 
\end{align}
Here $H^1(\R^{d+1})$ is the usual Sobolev space and $\langle f,g\rangle_{L^2(\R^{d+1})}= \int_{\R^{d+1}} f(x,t) \overline{g(x,t)} \dd x \dd t$. The realization of $\mathcal{A}'=-\nabla_x\cdot A' \nabla_x$ in $L^2 (\R^d)$ is given in the similar manner, for $A'$ satisfies the uniform ellipticity condition thanks to \eqref{ellipticity}. 
In this paper we are concerned with a factorization of the elliptic operators under weak assumptions for the matrix $A$, whose motivation lies in the application to the boundary value problems.  
The simplest example of $\mathcal{A}$ is the $(d+1)$-dimensional Laplacian $-\Delta= -\Delta_x - \partial_t^2 =  -\sum_{j=1}^d \partial_j^2 - \partial_t^2$. In this case it is straightforward to obtain the factorization
\begin{align}
-\Delta = - (\partial_t - (-\Delta_x)^\frac12 ) (\partial_t + (-\Delta_x)^\frac12 ). \label{factorization.laplace}
\end{align} 
The factorization \eqref{factorization.laplace} is  not only formal 
but also topological. 
Indeed, we see $D_{L^2}((-\Delta_x)^{1/2})=H^1 (\R^d)$, $D_{L^2}((\partial_t \pm  (-\Delta_x)^{1/2} ))=H^1 (\R^{d+1})$, and $D_{L^2}((\partial_t -  (-\Delta_x)^{1/2} ) (\partial_t + (-\Delta_x)^{1/2} ))=H^2 (\R^{d+1})$. Another key feature of \eqref{factorization.laplace} is that it is a factorization of the operator in the $t$ variable and the  $x$ variables. Hence, by the $t$-independent assumption for the coefficients,  the factorization into the first order differential operators as in 
\eqref{factorization.laplace} is easily extended to the case  when $A$ is a typical block matrix, i.e., ${\bf r_1} = {\bf r_2} ={\bf 0}$ and $b=1$, at least  in the formal level. Indeed, it suffices to replace $(-\Delta_x)^{1/2}$ by ${\mathcal{A}'}^{1/2}$, the square root of $\mathcal{A}'$ in $L^2 (\R^d)$. However, in contrast to the Laplacian case, the topological factorization is far from trivial  in this case, for one has to verify the relation $D_{L^2}({\mathcal{A}'}^{1/2})=H^1 (\R^d)$. This characterization is nothing but the Kato square root problem for divergence form elliptic operators, which was finally settled in \cite{AHLMT}. The first goal of the present paper  is to give sufficient
conditions on $A$, which may be a full entry matrix,  so that the exact topological factorization of $\mathcal{A}$ like \eqref{factorization.laplace} is verified. 
To this end we first introduce some terminologies.
 
\begin{df}\label{def.intro} {\rm (i)} For a given $h\in \mathcal{S}'(\R^d)$ we denote by $M_h: \mathcal{S}(\R^d)\rightarrow \mathcal{S}'(\R^d)$ the multiplication $M_h u = h u$.

\noindent {\rm (ii)} We denote by $E_{\mathcal{A}}: \dot{H}^{1/2} (\R^d)\rightarrow \dot{H}^1 (\R^{d+1}_+)$ the $\mathcal{A}$-extension operator, i.e., $u=E_{\mathcal{A}} g$ is the solution to the Dirichlet problem
\begin{equation}\label{eq.dirichlet0}
\begin{cases}
& \mathcal{A} u  = 0~~~~~~~{\rm in}~~~\R^{d+1}_+,\\
& \hspace{0.3cm} u  = g~~~~~~~{\rm on}~~\partial\R^{d+1}_+
=\R^d\times \{ t=0\}.
\end{cases}
\end{equation}
Here $\dot{H}^s (\Omega)$ denote the homogeneous Sobolev space on the domain $\Omega$ of order $s$.  The one parameter family of linear operators $\{E_{\mathcal{A}} (t)\}_{t\geq 0}$, defined by $E_{\mathcal{A}} (t) g = (E_{\mathcal{A}} g )(\cdot,t)$ for $g\in \dot{H}^{1/2}(\R^d)$, is called the Poisson semigroup associated with $\mathcal{A}$.
 
\noindent {\rm (iii)} We denote by $\Lambda_{\mathcal{A}}: \dot{H}^{1/2} (\R^d)\rightarrow \dot{H}^{-1/2} (\R^d)=\big (\dot{H}^{1/2} (\R^d) \big )^*$ the Dirichlet-Neumann map associated with $\mathcal{A}$, which is defined through the sesquilinear form 
\begin{equation}
\langle  \Lambda_{\mathcal{A}} g, \varphi \rangle_{\dot{H}^{-\frac12},\dot{H}^{\frac12}}  =  \langle  A\nabla E_{\mathcal{A}}g, \nabla E_{\mathcal{A}} \varphi \rangle_{L^2(\R^{d+1}_+)},~~~~~~~~~~g,\varphi \in \dot{H}^{\frac12}(\R^d).\label{def.Lambda}
\end{equation}
Here $\langle \cdot,\cdot\rangle _{\dot{H}^{-1/2},\dot{H}^{1/2}}$ denotes the duality coupling of  $\dot{H}^{-1/2}(\R^d)$ and $\dot{H}^{1/2}(\R^d)$.

\end{df}

\begin{rem}{\rm The notion of {\it semigroup} for $\{E_{\mathcal{A}}(t) \}_{t\geq 0}$ will be justified in Section \ref{subsec.poisson}. When $A$ is Hermite, the Dirichlet-Neumann map $\Lambda_{\mathcal{A}}$ is extended as a self-adjoint  operator in $L^2(\R^d)$, denoted again by $\Lambda_{\mathcal{A}}$. 
As is well known, 
even in the general case the ellipticity condition \eqref{ellipticity} ensures that $\Lambda_{\mathcal{A}}$ is extended as an injective {\it m}-sectorial operator in $L^2(\R^d)$ satisfying $D_{L^2} (\Lambda_{\mathcal{A}})\subset H^{1/2}(\R^d)$. 
}
\end{rem}
The adjoint matrix of $A$ will be denoted by $A^*$, and  $\mathcal{A}^*$ is the adjoint of $\mathcal{A}$ in $L^2 (\R^{d+1})$, which is a realization of $-\nabla \cdot A^* \nabla$ in $L^2 (\R^{d+1})$. Also 
$\mathcal{M} (\R^d)$ denotes the space of finite Radon measures, 
and $L^{p,\infty}(\R^d)$ denotes the Lorentz space $L^{p,q}(\R^d)$ with 
the exponent $q=\infty$. 
Then our first result is stated as follows.

%
%
%
%

\begin{thm}[Factorization for specific cases]
\label{thm.factorization.strong.1}  
Suppose that either 

\vspace{0.3cm}

\noindent {\rm (i)} $A$ is Hermite or both

\vspace{0.3cm}

\noindent {\rm (iia)} for $j=1,2$, $\nabla_x \cdot {\bf r}_j \in 
L^{d,\infty} (\R^d) + L^\infty (\R^d)$ if $d\ge 2$ 
(or $\nabla_x \cdot {\bf r}_j \in 
\mathcal{M}(\R) + L^\infty(\R)$ if $d=1$) 
with small $L^{d,\infty}(\R^d)$ parts
(or small $\mathcal{M}(\R)$ parts resp.) and 

\vspace{0.1cm}

\noindent {\rm (iib)}
${\rm Im} ~({\bf r_1 + r_2}) = {\bf 0}$ and ${\rm Im}~b =0$. 

\vspace{0.3cm}

\noindent Then  $H^1 (\R^d)$ is continuously embedded in $D_{L^2}(\Lambda_{\mathcal{A}})\cap D_{L^2}(\Lambda_{\mathcal{A}^*})$, and the operators $-{\bf P}_{\mathcal{A}}, ~-{\bf P}_{\mathcal{A}^*}$ defined by  
\begin{align}
& D_{L^2} ({\bf P}_{\mathcal{A}}) = H^1 (\R^d), ~~~~~~~~~ -{\bf P}_{\mathcal{A}} f =  -M_{1/b} \Lambda_{\mathcal{A}} f - M_{{\bf r_2}/b} \cdot \nabla _x f\,, \label{def.P.intro.1}\\
& D_{L^2} ({\bf P}_{\mathcal{A}^*}) = H^1 (\R^d), ~~~~~~~~-{\bf P}_{\mathcal{A}^*} f =  - M_{1/\bar{b}} \Lambda_{\mathcal{A}^*} f - M_{{\bf \bar{r}_1}/\bar{b}} \cdot \nabla_x f\,,\label{def.P.intro.2}
\end{align}
generate strongly continuous and analytic semigroups in $L^2 (\R^d)$. Moreover, the realization of $\mathcal{A}'$ in $L^2 (\R^d)$ and the realization  $\mathcal{A}$ in $L^2 (\R^{d+1})$  are respectively  factorized  as 
\begin{align}
\mathcal{A}' & = M_b \mathcal{Q}_{\mathcal{A}} {\bf P}_{\mathcal{A}},~~~~~~~~~~~~\mathcal{Q}_{\mathcal{A}} = M_{1/b} (M_{\bar{b}} {\bf P}_{\mathcal{A}^*})^*, \label{eq.factorization.strong.0}\\
\mathcal{A} & = - M_b (\partial_t -  \mathcal{Q}_{\mathcal{A}} ) (\partial_t + {\bf P}_{\mathcal{A}} ).\label{eq.factorization.strong.1}
\end{align}
Here $(M_{\bar{b}} {\bf P}_{\mathcal{A}^*})^*$ is the adjoint of $M_{\bar{b}} {\bf P}_{\mathcal{A}^*}$ in $L^2 (\R^d)$.

\end{thm}

\begin{rem}{\rm  In Theorem \ref{thm.factorization.strong.1} 
the operators $\partial_t + {\bf P}_{\mathcal{A}}$  and  
$\partial_t - \mathcal{Q}_{\mathcal{A}}$ 
are respectively defined as sum operators in $L^2(\R^{d+1})$, 
that is,
\begin{align*}
D_{L^2} (\partial_t + {\bf P}_{\mathcal{A}} ) & = \{ u \in L^2 (\R^{d+1})~|~\partial_t u, ~{\bf P}_{\mathcal{A}} u \in L^2 (\R^{d+1}) \} = H^1 (\R^{d+1}),\\
D_{L^2} (\partial_t - \mathcal{Q}_{\mathcal{A}} ) & = \{ u\in L^2 (\R^{d+1})~|~\partial_t u, ~\mathcal{Q}_{\mathcal{A}} u \in L^2 (\R^{d+1}) \}.
\end{align*}
Similarly, $M_b(\partial_t - \mathcal{Q}_{\mathcal{A}} ) (\partial_t + {\bf P}_{\mathcal{A}} )$ is defined as a product operator in $L^2 (\R^{d+1})$:
\begin{align*}
D_{L^2} (M_b(\partial_t - \mathcal{Q}_{\mathcal{A}} ) (\partial_t + {\bf P}_{\mathcal{A}} ) )= \{ u \in H^1 (\R^{d+1})~|~ (\partial_t + {\bf P}_{\mathcal{A}} ) u\in D_{L^2} (\partial_t - \mathcal{Q}_{\mathcal{A}})\}. 
\end{align*}


}
\end{rem}

\begin{rem}
\label{rem.ii}
{\rm  
When $A$ is Hermite the relation $D_{L^2}(\Lambda_{\mathcal{A}})=H^1 (\R^d)$ also holds; see Theorem \ref{thm.poisson.hermite}.  
}
\end{rem}

\begin{rem}{\rm
The regularity condition (iia) is natural in view of the scaling 
of the operators. In general one cannot remove this regularity condition on 
the divergence of the off-block vectors; see Section \ref{subsec.poisson.regular} for details. 
It should be noted that (iia) is satisfied if
 $\nabla_x \cdot {\bf r}_j$ belongs to $L^{d}(\R^d)+L^{\infty} (\R^d)$, since $C_0^\infty (\R^d)$ is dense in $L^d(\R^d)$.
The condition (iib) in Theorem \ref{thm.factorization.strong.1} 
is regarded as a kind of symmetry condition on $A$. 
Indeed, it is clear that if $A$ is Hermite or real 
valued then these conditions are satisfied. 
}
\end{rem}

\begin{rem}{\rm If the coefficients of $A$ are 
Lipschitz continuous, 
then the conclusion of Theorem \ref{thm.factorization.strong.1} 
holds even when the condition (iib) is absent; see \cite{MM2}. 
In this case we also have $D_{L^2}(\Lambda_{\mathcal{A}}) = D_{L^2}(\mathcal{Q}_{\mathcal{A}}) = H^1 (\R^d)$ with equivalent norms.
}

\end{rem}

\begin{rem}\label{rem.iii}{\rm
When $A$ possesses enough regularity it is classical in the theory of pseudo-differential operators that one looks for the factorization of $\mathcal{A}$ of the form $-M_{b}(\partial_t - \mathcal{A}_1 )  (\partial_t + \mathcal{A}_2)$ for some first order operators $\mathcal{A}_1$ and $\mathcal{A}_2$ but with modulo lower order operators; e.g. \cite{Treves}. On the other hand, \eqref{eq.factorization.strong.1} is just exact, i.e., any modifications by lower order operators are not required.
}
\end{rem}

\begin{rem}\label{rem.iv}{\rm 
The operator $-{\bf P}_{\mathcal{A}}$ is nothing but 
the generator $-\mathcal{P}_{\mathcal{A}}$ of the Poisson semigroup associated with $\mathcal{A}$. 
That is, $-{\bf P}_{\mathcal{A}} f = -\mathcal{P}_{\mathcal{A}} f := \lim_{t\downarrow 0} t^{-1} \big ( E_{\mathcal{A}} (t) f   - f \big )$ in $L^2 (\R^d)$, which is well-defined if $\{E_{\mathcal{A}}(t)\}_{t\geq 0}$ can be extended as a strongly continuous semigroup in $L^2 (\R^d)$. We will call 
$\mathbf{P}_{\mathcal{A}}$ and $\mathcal{P}_{\mathcal{A}}$
the \textit{Poisson operators}.
}
\end{rem}


When the coefficients of $A$ are merely bounded and measurable we cannot expect the representation \eqref{eq.factorization.strong.1}. 
In fact, it is revealed by \cite{KKPT} that \eqref{eq.dirichlet0} is not always solvable for boundary data in $L^2(\R^d)$ in general, which implies that the construction of  the Poisson semigroup in $L^2(\R^d)$ itself 
is impossible without additional assumptions. 
On the other hand, as is well-known, the ellipticity \eqref{ellipticity} is enough to realize the Poisson semigroup 
$\{E_{\mathcal{A}} (t)\}_{t\geq 0}$ in $H^{1/2}(\R^d)$ due to the Lax-Milgram theorem. The key feature here is that this semigroup is strongly continuous and analytic (see Proposition \ref{prop.poisson.analytic}), and  hence its generator $-\mathcal{P}_{\mathcal{A}}$ is  always well-defined in the functional setting of $H^{1/2}(\R^d)$.  Although the analyticity of  $\{E_{\mathcal{A}} (t)\}_{t\geq 0}$ in $H^{1/2}(\R^d)$ is not difficult to prove, it is an important observation since $D_{H^{1/2}}(\mathcal{P}_{\mathcal{A}})$ is shown to possess a strong regularity such as $D_{H^{1/2}}(\mathcal{P}_{\mathcal{A}})\hookrightarrow D_{L^2}(\Lambda_{\mathcal{A}})\cap H^1 (\R^d)$. 
The space $D_{H^{1/2}}(\mathcal{P}_{\mathcal{A}})$ is useful and plays an essential role in  rigorous derivation of several identities related to the factorization \eqref{eq.factorization.strong.1} without additional conditions on $A$ other than the ellipticity  \eqref{ellipticity}. These identities are described in Theorem \ref{thm.factorization} below, which can be regarded as  a counterpart of Theorem \ref{thm.factorization.strong.1} for general case. 
We denote  by  $\dot{H}^1 (\R; L^2 (\R^d))$ the function space $\{u\in L^2_{loc} (\R; L^2 (\R^d))~|~\partial_t u \in L^2 (\R^{d+1})\}$.

\begin{thm}[Identity in weak form for general coefficients  case]\label{thm.factorization} {\rm (i)} Assume that 
$$u\in L^2 (\R; D_{H^{1/2}}(\mathcal{P}_{\mathcal{A}}))\cap \dot{H}^1 (\R; L^2 (\R^d)),~~~~~v\in L^2 (\R; D_{H^{1/2}}(\mathcal{P}_{\mathcal{A}^*}))\cap \dot{H}^1 (\R; L^2 (\R^d)).$$
Then $u, v\in H^1(\R^{d+1})$ and  we have 
\begin{align}
\langle A\nabla u, \nabla v\rangle_{L^2 (\R^{d+1})} = \langle ~ (\partial_t + \mathcal{P}_{\mathcal{A}}) u, ~M_{\bar{b}} (\partial_t + \mathcal{P}_{\mathcal{A}^*} )v ~ \rangle _{L^2 (\R^{d+1})}. \label{eq.factorization}
\end{align}

\noindent {\rm (i')} Assume that $u\in L^2 (\R; D_{H^{1/2}} (\mathcal{P}_{\mathcal{A}} ))\cap \dot{H}^1(\R; H^{1/2}(\R^d))$ and $v\in H^1 (\R^{d+1})$.  Then $u\in H^1 (\R^{d+1})$ and \eqref{eq.factorization} is replaced by
\begin{align}
 \langle A\nabla u, \nabla v\rangle_{L^2 (\R^{d+1})} & = \langle ~ (\partial_t + \mathcal{P}_{\mathcal{A}}) u, ~ ( M_{\bar{b}} \partial_t +  M_{{\bf \bar{r}_1}}\cdot \nabla_x )v ~ \rangle _{L^2 (\R^{d+1})} \nonumber \\
&~~~~~ + \int_\R  \langle ~(\partial_t + \mathcal{P}_{\mathcal{A}}) u (t) , ~ \Lambda_{\mathcal{A}^*} v (t)  ~ \rangle _{\dot{H}^\frac12, \dot{H}^{-\frac12}} \dd t. \label{eq.factorization'}
\end{align}

\noindent {\rm (ii)} Assume that 
$$u\in L^2 (\R_+; D_{H^{1/2}}(\mathcal{P}_{\mathcal{A}}))\cap \dot{H}^1 (\R_+; L^2 (\R^d)),~~~~~v\in L^2 (\R_+; D_{H^{1/2}}(\mathcal{P}_{\mathcal{A}^*}))\cap \dot{H}^1 (\R_+; L^2 (\R^d)).$$
Then $u,v\in H^1(\R^{d+1}_+)$ and  we have 
\begin{align}
\langle A\nabla  u, \nabla v \rangle _{L^2 (\R^{d+1}_+)} = \langle  \gamma u, ~ \Lambda_{\mathcal{A}^*} \gamma  v \rangle _{\dot{H}^{\frac12}, \dot{H}^{-\frac12}} + \langle ~ (\partial_t + \mathcal{P}_{\mathcal{A}}) u, ~M_{\bar{b}} (\partial_t + \mathcal{P}_{\mathcal{A}^*} )v ~ \rangle _{L^2 (\R^{d+1}_+)}. \label{eq.factorization.half}
\end{align}
Here $\gamma: H^1 (\R^{d+1}_+)\rightarrow H^{1/2}(\partial\R^{d+1}_+)=H^{1/2}(\R^d)$ is  the trace operator.

\noindent {\rm (ii')} Assume that $u\in L^2 (\R_+; D_{H^{1/2}}(\mathcal{P}_{\mathcal{A}}))\cap \dot{H}^1 (\R_+; H^{1/2} (\R^d))$ and $v\in H^1 (\R^{d+1}_+)$.  Then $u\in H^1(\R^{d+1}_+)$ and \eqref{eq.factorization.half} is replaced by 
\begin{align}
\langle A \nabla  u,  \nabla v  \rangle _{L^2 (\R^{d+1}_+)} & =   \langle  \gamma u, ~ \Lambda_{\mathcal{A}^*} \gamma  v \rangle _{\dot{H}^{\frac12}, \dot{H}^{-\frac12}} + \int_0^\infty \langle~ (\partial_t + \mathcal{P}_{\mathcal{A}}) u (t) , ~ \Lambda_{\mathcal{A}^*} v (t)   ~ \rangle _{\dot{H}^{\frac12}, \dot{H}^{-\frac12}} \dd t \nonumber \\
& ~~~ +  \langle  ~(\partial_t + \mathcal{P}_{\mathcal{A}}) u (t) , ~ ( M_{\bar{b}}\partial_t + M_{\bf \bar{r}_1} \cdot \nabla_x  )  v (t)   ~ \rangle _{L^2 (\R^{d+1}_+)}. \label{eq.thm.factorization.half}
\end{align}

\end{thm} 

\begin{rem}\label{rem.thm.factorization}{\rm The identities \eqref{eq.factorization}-\eqref{eq.thm.factorization.half} are closely related to the Rellich type identity \cite{Rellich}. However, it is important to specify the class of $u$ and $v$ for which  \eqref{eq.factorization}-\eqref{eq.thm.factorization.half} are verified. Indeed,  it is not clear whether or not  \eqref{eq.factorization} holds  even for functions in $C_0^\infty (\R^{d+1})$ if we impose only \eqref{ellipticity}, since the inclusion  $C_0^\infty (\R^d)\subset D_{H^{1/2}}(\mathcal{P}_{\mathcal{A}})\cap D_{H^{1/2}}(\mathcal{P}_{\mathcal{A}^*})$ is not known  when $A$ is not Hermite and  nonsmooth. The identities \eqref{eq.factorization'} and \eqref{eq.thm.factorization.half} are sometimes more useful than \eqref{eq.factorization} and \eqref{eq.factorization.half}, for in general the detailed informations on  $D_{H^{1/2}} (\mathcal{P}_{\mathcal{A}^*} )$ are lacking, while only $H^1$ regularity is required for $v$ in  \eqref{eq.factorization'} and \eqref{eq.thm.factorization.half}. When $\{e^{-t\mathcal{P}_{\mathcal{A}}}\}_{t\geq 0}$ and $\{e^{-t\mathcal{P}_{\mathcal{A}^*}}\}_{t\geq 0}$ are extended as strongly continuous semigroups in $L^2 (\R^d)$ and $D_{L^2}(\mathcal{P}_{\mathcal{A}})=D_{L^2}(\mathcal{P}_{\mathcal{A}^*})=H^1 (\R^d)$ holds,  the assumptions  $D_{H^{1/2}}(\mathcal{P}_{\mathcal{A}}), D_{H^{1/2}}(\mathcal{P}_{\mathcal{A}^*})$ in \eqref{eq.factorization} or \eqref{eq.factorization.half} can be simply replaced by $H^1 (\R^d)$.

}
\end{rem}

As an application of our result, we will consider 
$L^2$ solvability of the Dirichlet problem 
\begin{equation}\label{eq.dirichlet}
\tag{D}
\begin{cases}
& \mathcal{A} u  = F~~~~~~~~{\rm in}~~~\R^{d+1}_+,\\
& ~~u  = g~~~~~~~~~{\rm on}~~\partial\R^{d+1}_+,
\end{cases}
\end{equation}
and the Neumann problem
\begin{equation}
\label{eq.neumann}
\tag{N}
\begin{cases}
& ~~~~~~~~~~~\hspace{0.3cm}~ \mathcal{A} v  = F~~~~~~~{\rm in}~~~\R^{d+1}_+,\\
& -\langle {\bf e}_{d+1}, A \nabla v \rangle  = g~~~~~~~{\rm on}~~\partial\R^{d+1}_+.
\end{cases}
\end{equation}
There is a lot of literature on the solvability of these problems
under weak assumptions for the matrix $A$ 
and the data $F$, $g$. 
In the case $F \in \dot{H}^{-1}(\R^{d+1}_+)$ and  $g\in H^{1/2}(\R^d)$, 
the Lax-Milgram theorem is directly applied to find the unique solution in $\dot{H}^1 (\R^{d+1}_+)$. However, in other cases 
the abstract theory from functional analysis is not always applicable and the problem becomes subtle.
The factorization formula \eqref{eq.factorization.strong.1} has an important application to the boundary value problems, for it leads to the 
formal representations of solutions  to \eqref{eq.dirichlet};
\begin{align}
u (t) = e^{-t\mathcal{P}_{\mathcal{A}}} g +  \int_0^t e^{-(t-s)\mathcal{P}_\mathcal{A}}\int_s^\infty e^{-(\tau-s)\mathcal{Q}_\mathcal{A}} M_{1/b} F (\tau ) \dd\tau \dd s, \label{eq.thm.representation.1}
\end{align}
and to \eqref{eq.neumann};
\begin{align}
& v (t) =  e^{-t \mathcal{P}_{\mathcal{A}}} \Lambda_{\mathcal{A}}^{-1} \big ( g + M_b   \int_0^\infty   e^{-s \mathcal{Q}_{\mathcal{A}}} M_{1/b} F (s) \dd s \big )   + \int_0^t e^{-(t-s) \mathcal{P}_{\mathcal{A}}} \int_s^\infty  e^{-(\tau-s) \mathcal{Q}_{\mathcal{A}}} M_{1/b} F (\tau)  \dd \tau \dd s. 
\label{eq.thm.representation.3}
\end{align}
Here the semigroup $\{e^{-t\mathcal{Q}_{\mathcal{A}}}\}_{t\geq 0}$ is given by
\begin{align}
e^{-t{Q}_{\mathcal{A}}} = M_{1/b} (e^{-t \mathcal{P}_{\mathcal{A}^*}} )^* M_b,\label{def.e^Q.intro}
\end{align}
which is well-defined  if the Poisson semigroup $\{E_{\mathcal{A}^*} (t) \}_{t\geq 0} = \{e^{-t \mathcal{P}_{\mathcal{A}^*}}\}_{t\geq 0}$ in $H^{1/2}(\R^d)$ (which is always well-defined) is extended as a semigroup 
in $L^2 (\R^d)$. These formula reduce the inhomogeneous problems \eqref{eq.dirichlet} and
 \eqref{eq.neumann} to the analysis of the semigroups $\{e^{-t {\mathcal{P}
}_{\mathcal{A}}}\}_{t\geq 0}$ and $\{e^{-t {Q}_{\mathcal{A}} }\}_{t\geq 0}$. 
We will call $u$ and $v$ the mild solutions, if \eqref{eq.thm.representation.1} and \eqref{eq.thm.representation.3} are well-defined; see Definition \ref{def.mild.solution.pre}. In order to clarify the relation between weak solutions and mild solutions
it is important to study the domain of the Poisson operators, which will be discussed in Section \ref{subsec.mild.weak}.

Now let us state some results on $L^2$ solvability 
of \eqref{eq.dirichlet} and \eqref{eq.neumann} in the simplest form. We set $\overline{\R_+}=[0,\infty)$, and for a Banach space $X$ we write $f\in C(\overline{\R_+}; X)$ if and only if $f\in C([0,T); X)$ for all $T>0$.
For the homogeneous problems (i.e., $F=0$ in \eqref{eq.dirichlet} or \eqref{eq.neumann}), Theorem 
\ref{thm.factorization.strong.1} implies the following result:

\begin{cor}
\label{thm.solvability.intro1}
Under the assumptions in Theorem 
\ref{thm.factorization.strong.1}, there exists a unique weak solution $u$ to 
\eqref{eq.dirichlet} with $F=0$ and $g\in L^2 (\R^d)$ such that $u\in C(\overline{\R_+};L^2(\R^{d}))$\ $\cap$\ $
\dot{H}^1(\R^d \times (\delta,\infty))$
for any $\delta>0$. 
If in addition $g$ belongs to the range of 
$\Lambda_{\mathcal{A}}$, then 
there exists a unique weak solution $v$ to \eqref{eq.neumann} 
with $F=0$ such that $v\in C(\overline{\R_+};H^{1/2}(\R^{d})) \cap \dot{H}^1(\R^{d+1}_+)$.
\end{cor}

\begin{rem}{\rm As we mentioned before, if $A$ is Hermite then $D_{L^2}(\Lambda_{\mathcal{A}})=H^1 (\R^d)$ holds. In this case the weak solution to \eqref{eq.neumann} obtained in Corollary \ref{thm.solvability.intro1} possesses further regularity such as $C(\overline{\R_+};H^{1}(\R^{d}))$.

}
\end{rem}
\begin{rem}{\rm
It is well-known that solvability of the homogeneous boundary value problems in $\R^{d+1}_+$
can be extended to that in the domain above a Lipschitz graph. The
$L^2$ solvability of the Laplace equation (i.e., $A=I$) in Lipschitz domains was shown in \cite{Dahlberg1, JK1, Verchota}. 
This result was extended by \cite{JK2, KP,AAAHK}  to the case when $A$ is real symmetric, and by \cite{AAM} to the case when $A$ is Hermite. 
In view of $L^2$ solvability of the  homogeneous boundary value problems, Corollary \ref{thm.solvability.intro1} gives a new contribution under the conditions {\rm (iia) }- {\rm (iib)} in Theorem \ref{thm.factorization.strong.1}.  
}
\end{rem}

When $A$ is not Hermite and nonsmooth, the boundary value problems are not always solvable for $L^2$ boundary data. 
If $A$ is a typical block matrix,  ${\bf r_1}={\bf r_2}={\bf 0}$ and $b=1$, then the homogenous Dirichlet problem 
is easily solved by using the semigroup theory, while the homogeneous Neumann problem in this case is  
equivalent with the Kato square root problem. 
Recently the authors 
in \cite{AAH} showed $L^2$ solvability of the homogeneous 
Dirichlet and Neumann problems when $A$ is a small $L^\infty$ 
perturbation of a block matrix; see also \cite{AAM,AAAHK,AA,AR} for related stability result.
In fact, Theorem \ref{thm.solvability.intro1} 
under the conditions (iia) - (iib) can be
regarded as another stability result for the block matrix case. 
Note that the principal scaling  of $\nabla_x\cdot {\bf r_j}$ in (iia), i.e., $\|\nabla_x \cdot {\bf r}_j\|_{L^{d,\infty}(\R^d)}$
for $d\ge 2$ or 
$\|\nabla_x \cdot {\bf r}_j\|_{\mathcal{M}(\R^d)}$ for $d=1$,
is same as that of 
$\|a_{i,j}\|_{L^{\infty}(\R^d)}$. Thus, in view of the scaling, the condition (iia) in Theorem \ref{thm.solvability.intro1}  is comparable to $L^\infty$ perturbations discussed in \cite{AAH,AAM,AAAHK}. 
On the other hand, for general real nonsymmetric $A$, $L^p$ solvability 
of the homogeneous Dirichlet problems in the half plane $\R^{2}_+$
was obtained in \cite{KKPT} for sufficiently large $p$ depending on $A$.
They also constructed an example of  the matrix $A$ such that the homogeneous Dirichlet problem in $\R^2_+$ 
is not solvable for the boundary data in $L^2(\R)$. 
In their example, $A$ is real but nonsymmetric, and $\nabla_x \cdot {\bf r_j}$ ($j=1,2$) is a Dirac measure whose mass is not small. This example 
shows the optimality of  our condition (iia) 
for the case of real nonsymmetric matrices when $d=1$. 
For further results on solvability of the homogeneous problems, 
see \cite{Kenig} and references therein.

The next result concerns $L^2$ solvability of the 
inhomogeneous problems. For simplicity of the presentation,
we will assume the boundary data are zero.

\begin{thm}
\label{thm.solvability.intro2}
\ Suppose that either

\vspace{0.3cm}

\noindent {\rm (i)} $A$ is Hermite or both

\vspace{0.3cm}

\noindent {\rm (iia')} $\nabla_x\cdot {\bf r_1}=0$ and  $\nabla_x \cdot {\bf r}_2 \in 
L^{d,\infty} (\R^d) + L^\infty (\R^d)$ if $d\ge 2$ 
(or $\nabla_x \cdot {\bf r}_2 \in 
\mathcal{M}(\R) + L^\infty(\R)$ if $d=1$) 
with small $L^{d,\infty}(\R^d)$ parts
(or small $\mathcal{M}(\R)$ parts resp.) and 

\vspace{0.1cm}

\noindent {\rm (iib')} ${\bf r_1}$, ${\bf r_2}$, and $b$ are real-valued.

\vspace{0.3cm}

\noindent Then  for given $F \in L^1 (\R_+; L^2 (\R^d))$ 
there exists a weak solution $u$ to \eqref{eq.dirichlet}
with  $g=0$ satisfying 
$$
u \in C(\overline{\R_+};L^2(\R^{d})) \quad \textit{and} \quad  
\nabla u  \in L^p_{loc}(\overline{\R_+};L^2(\R^d)) \ \ 
\textit{for any} \ \, p\in [1,\infty).
$$
If in addition $h =M_{b}  \int_0^\infty   e^{-s \mathcal{Q}_{\mathcal{A}}} M_{1/b} F (s) \dd s$ belongs to the range of $\Lambda_{\mathcal{A}}$, 
then there exists a weak solution $v$ to \eqref{eq.neumann} 
with $g=0$ satisfying
$$
v \in C(\overline{\R_+};L^2(\R^{d})) \quad \textit{and} \quad 
\nabla v  \in L^p_{loc}(\overline{\R_+};L^2(\R^d))
\ \ \textit{for any} \ \, p \in [1,2).
$$
\end{thm}
\begin{rem}{\rm Under the assumptions of Theorem \ref{thm.solvability.intro2} the semigroups $\{e^{-t\mathcal{P}_{\mathcal{A}}}\}_{t\geq 0}$ and  $\{e^{-t\mathcal{P}_{\mathcal{A}^*}}\}_{t\geq 0}$ acting on $H^{1/2}(\R^d)$ are extended as strongly continuous semigroups acting on $L^2 (\R^d)$  thanks to the results of Theorem \ref{thm.factorization.strong.1}.  

}
\end{rem}
\begin{rem}{\rm 
There is a lot of literature for the inhomogeneous boundary value
problems in bounded Lipschitz domains; see, e.g.,
\cite{Dahlberg, JK3, FMM, Mitrea1, Mitrea2, MiMi,NS}
and references therein.
As well as the case for the homogeneous problem, 
Theorem \ref{thm.solvability.intro2} for Hermite matrices
yields $L^2$ solvability of the inhomogeneous problems 
for matrices of the same type in domains 
above Lipschitz graphs. 
For the Laplace equation, $L^p$ solvability of 
the inhomogenous problems in bounded Lipschitz domains 
was proved in \cite{Dahlberg, JK3, FMM}. 
Our result also shows the gradient of the Dirichlet Green operator 
(i.e., the solution map for \eqref{eq.dirichlet} 
with the zero boundary data: $F \mapsto \nabla u$) 
maps $L^1(\R_+;L^2(\R^d))$ continuously to 
$L^p_{loc}(\R_+,L^2(\R^d))$.
Results of this type go back to \cite{Dahlberg} where
the author showed that the gradient of the Dirichlet 
Green operator for $\mathcal{A}=-\Delta$ 
in the bounded Lipschitz domain is a continuous map from 
$L^1(\Omega)$ to $L^{n/(n-1),\infty}(\Omega)$.
Recently, it was generalized in \cite{Mitrea1}
for the Neumann Green operator by using potential technique; 
see also \cite{Mitrea2, MiMi} for further results.}
\end{rem}

As is well-known in the spectral theory, 
it is a subtle problem to determine sufficient conditions for 
$F$ to solve the problems \eqref{eq.dirichlet} or \eqref{eq.neumann}.
Indeed, due to the lack of the Poincar{\'e} inequality, the origin belongs to the continuous spectrum of $\mathcal{A}$ (with the zero boundary condition) in $L^2(\R^{d+1}_+)$. 
Hence the inhomogeneous problem is not always solvable for $F\in L^2 (\R^{d+1}_+)$, even if $A$ is real symmetric 
and smooth. 
Therefore some additional conditions related to the spatial 
decay have to be imposed on $F$ to find the solution. 
Furthermore, the solution may fail to 
decay at spatial infinity even if it exists. 
To show Theorem \ref{thm.solvability.intro2} we will make use of the representation formulas 
\eqref{eq.thm.representation.1} and \eqref{eq.thm.representation.3}. Then  it is clear that
the temporal decay of $e^{-t\mathcal{Q}_{\mathcal{A}}}$ 
is crucial for solving our problems.
In fact,  the conditions in Theorem \ref{thm.solvability.intro2} guarantee the boundedness of  the semigroup $\{e^{-t\mathcal{\mathcal{Q}}_{\mathcal{A}}}\}_{t\geq 0}$ in $L^2 (\R^d)$, and hence,  the integrals in 
\eqref{eq.thm.representation.1} and  \eqref{eq.thm.representation.3} converge absolutely if $F\in L^1(\R_+; L^2 (\R^d))$. 
By a simple observation of the scaling, it is easy to see 
that the space $L^1(\R_+;L^2(\R^d))$ includes some functions 
decaying more slowly at (time) infinity than those in 
$\dot{H}^{-1}(\R^{d+1}_+)$.
In this sense, our result 
generalizes the class of the inhomogeneous terms 
for the solvability 
in terms of the decay at infinity. In should be emphasized here that the factorization in Theorem \ref{thm.factorization.strong.1} plays an essential role behind the proof of Theorem \ref{thm.solvability.intro2},
for the representation formulas  such as \eqref{eq.thm.representation.1} and \eqref{eq.thm.representation.3} are nothing but a result of \eqref{eq.factorization.strong.1}.
In Section \ref{subsec.solvability}, we will state 
a detailed version of Theorem \ref{thm.solvability.intro2}. 

We finally note that our approach  provides a unified view for the Dirichlet and 
Neumann  problems to both homogeneous and inhomogeneous  equations.  
In \cite{MM2} the results of the present paper  is applied to the  study of the Helmholtz decomposition of vector fields in domains with noncompact boundary.

The rest of this paper is organized as follows. Section \ref{sec.preliminary} collects some preliminarily facts on the Poisson semigroup. In particular, in Section \ref{subsec.def.weak.sol} we give the notions of weak solutions, and in Section \ref{subsec.poisson} the analyticity of the Poisson semigroup in $H^{1/2}(\R^d)$ is proved. The relation $D_{H^{1/2}}(\mathcal{P}_{\mathcal{A}})\hookrightarrow H^1 (\R^d)\cap D_{L^2} (\Lambda_{\mathcal{A}})$ will be shown in Section \ref{subsec.relation.Lambda.P}. Using this embedding property, we prove the Rellich type identity and Theorem \ref{thm.factorization} in Section \ref{subsec.thm.factorization}. Then the relation between the domain of the Poisson operator and $H^1 (\R^d)$ is investigated in Section \ref{subsec.embedding.poisson.H^1}. The factorization of $\mathcal{A}$ in $L^2 (\R^d)$ is established in an abstract setting in Section \ref{subsec.lem.factorization}. Based on this result, we analyze in Section \ref{sec.poisson} two specific classes of $A$ as described in Theorem \ref{thm.factorization.strong.1}. In Section \ref{sec.application} we apply Theorem \ref{thm.factorization.strong.1} to the boundary value problem. In particular,  the relation between mild solutions and weak solutions is discussed in Section \ref{subsec.mild.weak}, and Theorem \ref{thm.solvability.intro2} and its extension are proved in Section \ref{subsec.solvability}. Finally, we state the counterpart of Theorem \ref{thm.factorization.strong.1} for $\mathcal{A}+\lambda$, $\lambda\in \C$ with ${\rm Re} ~\lambda > 0$ in the appendix without proof.

%
%
%
%

\section{Preliminaries}\label{sec.preliminary}

In this preliminary section we first define the notion of weak solutions and mild solutions, and next collect fundamental  results on  Poisson semigroups and Dirichlet-Neumann maps. We note that most of the results stated in this section are valid only under the ellipticity condition \eqref{ellipticity}, without any other assumptions on the matrix $A$.

\subsection{Definitions of weak solution and mild solution}
\label{subsec.def.weak.sol}

We first define the notion of weak solutions treated in this paper. 
Set $\R^{d+1}_{\delta,+} = \{ (x,t) \in \R^{d+1}~|~t > \delta\}$. 

\begin{df} Let $F\in L^1_{loc}(\R^{d+1}_+)$ and $g\in L^1_{loc} (\R^d)$. 

\noindent (i) We say that $u\in L^1_{loc} (\R^{d+1}_+)$ is a weak solution to \eqref{eq.dirichlet} if $u$ belongs to $\displaystyle {\large\cap}_{\delta>0}  W^{1,1}_{loc}(\R^{d+1}_{\delta,+})$ and satisfies 
\begin{align}
\langle A\nabla u, \nabla \varphi \rangle _{L^2 (\R^{d+1}_+)} = \langle F,\varphi \rangle _{L^2(\R^{d+1}_+)}~~~~~~~~{\rm for~all}~~\varphi\in C_0^\infty (\R^{d+1}_+) \label{def.weak.dirichlet}
\end{align}
and ~$u(t) \rightarrow  g$ as $t\rightarrow + 0$ in the sense of distributions.

\noindent (ii)  We say that $u\in L^1_{loc} (\R^{d+1}_+)$ is a weak solution to \eqref{eq.neumann} if $u$ belongs to $W^{1,1}_{loc}(\R^{d+1}_+)$ and satisfies 
\begin{align}
\langle A\nabla u, \nabla \varphi \rangle _{L^2 (\R^{d+1}_+)} = \langle F,\varphi \rangle _{L^2(\R^{d+1}_+)} + \langle g, \gamma \varphi\rangle _{L^2 (\partial\R^{d+1}_+)}~~~~~~~~{\rm for~all}~~\varphi \in  C_0^\infty (\overline{\R^{d+1}_+}).\label{def.weak.neumann}
\end{align}
Here $\langle g, \gamma \varphi\rangle _{L^2 (\partial\R^{d+1}_+)} = \int_{\R^d} g(x) \overline{\varphi (x,0)}\dd x$.  
\end{df}

As stated in the introduction, it is natural to introduce the notion of mild solutions for elliptic boundary value problems in terms of the Poisson semigroups. 
\begin{df}[Mild solution]\label{def.mild.solution.pre} Assume that the Poisson semigroups $\{e^{-t\mathcal{P}_{\mathcal{A}}}\}_{t\geq 0}$ and  $\{e^{-t\mathcal{P}_{\mathcal{A}^*}}\}_{t\geq 0}$ acting on $H^{1/2}(\R^d)$ are extended as  strongly continuous semigroups acting on $L^2 (\R^d)$. We set the one-parameter family  $\{e^{-t\mathcal{Q}_{\mathcal{A}}}\}_{t\geq 0}$ by 
\begin{align}
e^{-t \mathcal{Q}_{\mathcal{A}}} = M_{1/b} (e^{-t \mathcal{P}_{\mathcal{A}^*}} )^* M_b,\label{def.mild.solution.Q}
\end{align}
which is again a strongly continuous semigroup acting on $L^2 (\R^d)$. Let $F \in L^1_{loc}(\R_+; L^2 (\R^d))$ and $g\in L^2(\R^d)$.  If  the function $u\in L^1_{loc}(\R^{d+1}_+)$ has the well-defined representation 
\begin{align}
u (t) = e^{-t\mathcal{P}_{\mathcal{A}}} g +  \int_0^t e^{-(t-s)\mathcal{P}_\mathcal{A}}\int_s^\infty e^{-(\tau-s)\mathcal{Q}_\mathcal{A}} M_{1/b} F (\tau ) \dd\tau \dd s, \label{eq.thm.representation.1}
\end{align}
then we call $u$ a mild solution to \eqref{eq.dirichlet}. Similarly,  if the function  $v\in L^1_{loc}(\R^{d+1}_+)$ has the well-defined representation  
\begin{align}
& v (t) =  e^{-t \mathcal{P}_{\mathcal{A}}} \Lambda_{\mathcal{A}}^{-1} \big ( g + M_b   \int_0^\infty   e^{-s \mathcal{Q}_{\mathcal{A}}} M_{1/b} F (s) \dd s \big )   + \int_0^t e^{-(t-s) \mathcal{P}_{\mathcal{A}}} \int_s^\infty  e^{-(\tau-s) \mathcal{Q}_{\mathcal{A}}} M_{1/b} F (\tau)  \dd \tau \dd s, \label{eq.thm.representation.3}
\end{align}
then we call $v$ a mild solution to \eqref{eq.neumann}. 

\end{df}

\begin{rem}{\rm It is easy to see that the generator of the semigroup $\{e^{-t\mathcal{Q}_{\mathcal{A}}}\}_{t\geq 0}$ defined by \eqref{def.mild.solution.Q} is $-M_{1/b} ( M_{\bar{b}} \mathcal{P}_{\mathcal{A}^*})^*$, where  $-\mathcal{P}_{\mathcal{A}^*}$ is the generator of $\{e^{-t\mathcal{P}_{\mathcal{A}^*}}\}_{t\geq 0}$ in $L^2 (\R^d)$.
}
\end{rem}

\subsection{Construction of $\mathcal{A}$-extension operator and Dirichlet-Neumann map}\label{subsec.extension.operator}

In this section we give a realization of the $\mathcal{A}$-extension operator and the Dirichlet-Neumann map in Definition \ref{def.intro}. The argument is rather standard. For $f \in \dot{H}^{1/2} (\R^d)$ let $E_{-\Delta} f\in \dot{H}^1(\R^{d+1}_+)$ be the classical harmonic extension of $f$ to $\R^{d+1}_+$, i.e., $E_{-\Delta} f$ is the solution to \eqref{eq.dirichlet0} for $\mathcal{A}=-\Delta$. Then $E_{\mathcal{A}} f$ is constructed in the form 
\begin{align}
E_{\mathcal{A}} f = E_{-\Delta} f + \mathcal{A}^{-1}_D \nabla\cdot (A-I) \nabla E_{-\Delta} f,
\end{align}
where $\mathcal{A}^{-1}_D \nabla\cdot (A-I) \nabla E_{-\Delta} f$ is the solution to \eqref{eq.dirichlet} with $F=\nabla\cdot (A-I) \nabla E_{-\Delta} f$ and $g=0$.  To be precise,  for  $w\in \dot{H}^1_0 (\R^{d+1}_+)$ we  define the distribution $\mathcal{A}_D w\in \dot{H}^{-1}(\R^{d+1}_+)$   by $\langle \mathcal{A}_D w, v\rangle _{\dot{H}^{-1},\dot{H}^1_0} = \langle A \nabla w, \nabla v \rangle _{L^2(\R^{d+1}_+)}$. We note that $\mathcal{A}_D$ is realized as an {\it m}-sectorial operator $\mathcal{A}_D: D_{L^2}(\mathcal{A}_D)\subset H^1_0(\R^{d+1}_+) \rightarrow L^2 (\R^{d+1}_+)$  defined by $\langle \mathcal{A}_D w, v\rangle _{L^2(\R^{d+1}_+)} = \langle A \nabla w, \nabla v \rangle _{L^2(\R^{d+1}_+)}$ for $w\in D_{L^2}(\mathcal{A}_D)$ and $v\in {H}^1_0(\R^{d+1}_+)$; \cite[Chapter VI-2]{Kato}.  Then  the function $w =\mathcal{A}^{-1}_D \nabla\cdot (A-I) \nabla E_{-\Delta} u \in \dot{H}^1_0 (\R^{d+1}_+)$ is the solution to  $\langle \mathcal{A}_D w, v\rangle _{\dot{H}^{-1},\dot{H}^1_0} = - \langle  (A-I)\nabla E_{-\Delta} f, \nabla v \rangle _{L^2(\R^{d+1}_+)}$ for all $v\in \dot{H}^1_0(\R^{d+1}_+)$. The unique existence of such $w\in \dot{H}^1_0 (\R^{d+1}_+)$ is ensured by the Lax-Milgram theorem due to \eqref{ellipticity}, and it follows that $\| \nabla w \|_{L^2(\R^{d+1}_+)}\leq C \| \nabla E_{-\Delta} f\|_{L^2(\R^{d+1}_+)}$. By the definition we have $\| \nabla E_{-\Delta} f \|_{L^2(\R^{d+1}_+)} = \| f \|_{\dot{H}^{1/2}(\R^d)}<\infty$. In particular, $E_{\mathcal{A}}$ is extended as a bounded operator from $\dot{H}^{1/2}(\R^d)$ to $\dot{H}^1(\R^{d+1}_+)$. Moreover, from the construction of $E_{\mathcal{A}}$ we have $\langle A\nabla E_{\mathcal{A}} f, \nabla v \rangle_{L^2(\R^{d+1}_+)}=0$ for any $v\in \dot{H}^1_0(\R^{d+1}_+)$.


\vspace{0.3cm}

Next we construct the Dirichlet-Neumann map $\Lambda_{\mathcal{A}}$ associated with the elliptic operator $\mathcal{A}$. As in Definition \ref{def.intro} we define the operator $\Lambda_{\mathcal{A}}: H^{1/2} (\R^d)\rightarrow \dot{H}^{-1/2} (\R^d)$  by 
\begin{equation}
\langle  \Lambda_{\mathcal{A}} f, g\rangle_{\dot{H}^{-\frac12},\dot{H}^{\frac12}}  = \langle  A\nabla E_{\mathcal{A}}f, \nabla E_{\mathcal{A}} g \rangle_{L^2(\R^{d+1}_+)}.\label{def.Lambda'}
\end{equation}
Then the ellipticity \eqref{ellipticity} implies ${\rm Re}~\langle \Lambda_{\mathcal{A}} f, g\rangle _{\dot{H}^{-1/2},\dot{H}^{1/2}} \geq \nu_1 \| \nabla E_{\mathcal{A}} f \|_{L^2(\R^{d+1}_+)}^2 \geq \nu_1 \| f \|_{\dot{H}^{1/2}(\R^d)}^2$, while $|\langle   \Lambda_{\mathcal{A}} f, g \rangle _{\dot{H}^{-1/2},\dot{H}^{1/2}}|\leq \nu_2 \| \nabla E_{\mathcal{A}} f \|_{L^2(\R^{d+1}_+)} \| \nabla E_{\mathcal{A}} g \|_{L^2(\R^{d+1}_+)} \leq C \| f \|_{\dot{H}^{1/2}(\R^d)} \| g \|_{\dot{H}^{1/2}(\R^d)}$ from the definition and the estimate of $E_{\mathcal{A}}$ in Section \ref{subsec.extension.operator}. Thus the theory of the sesquilinear forms \cite[Chapter VI-2]{Kato} shows that there exists an {\it m}-sectorial operator, denoted again by $\Lambda_{\mathcal{A}}$, such that $D_{L^2} (\Lambda_{\mathcal{A}})\subset  H^{1/2}(\R^d)$ and $\langle \Lambda_{\mathcal{A}} f, g\rangle _{L^2(\R^d)} = \langle \Lambda_{\mathcal{A}} f, g\rangle _{\dot{H}^{-1/2},\dot{H}^{1/2}}$ for $f\in D_{L^2} (\Lambda_{\mathcal{A}})$ and $g\in H^{1/2}(\R^d)$. It is clear that $\Lambda_{\mathcal{A}}$ is injective.  Note that 
\begin{align}
\langle  A\nabla E_{\mathcal{A}}f, \nabla E_{\mathcal{A}} g \rangle_{L^2(\R^{d+1}_+)} = \langle A \nabla E_{\mathcal{A}} f, \nabla  v\rangle _{L^2 (\R^{d+1}_+)}\label{dn.hosoku}
\end{align}
for any $v\in \dot{H}^1(\R^{d+1}_+)$ whose trace on $\partial\R^{d+1}_+=\R^d$ is $g$. Thus we have $\langle  A\nabla E_{\mathcal{A}}f, \nabla E_{\mathcal{A}} g \rangle_{L^2(\R^{d+1}_+)}  = \langle  \nabla E_{\mathcal{A}^*} f, A^* \nabla E_{\mathcal{A}^*} g \rangle_{L^2(\R^{d+1}_+)}$. This identity implies that the adjoint operator $\Lambda_{\mathcal{A}}^*$ of $\Lambda_{\mathcal{A}}$ in $L^2(\R^d)$ is given by $\Lambda_{\mathcal{A}^*}$, where $\mathcal{A}^*=\nabla\cdot A^*\nabla$ and $A^*$ is the adjoint matrix of $A$. In particular, $\Lambda_{\mathcal{A}}^{**} = \Lambda_{\mathcal{A}}$ holds.

\subsection{Construction of  Poisson operator in $H^{1/2}(\R^d)$}
\label{subsec.poisson}

In this section we consider the Poisson semigroup $\{E_{\mathcal{A}}(t)\}_{t\geq 0}$, where $E_{\mathcal{A}} (0)=I$ is the identity map and 
\begin{equation}
E_{\mathcal{A}} (t) f = \big (E_{\mathcal{A}} f \big )(\cdot, t), ~~~~~~~~~ t>0,~~~f \in {H}^\frac12 (\R^d).
\label{def.poisson.semigroup}
\end{equation}

\noindent In \eqref{def.poisson.semigroup} it is important to take $f$ from ${H}^{1/2}(\R^d)$ rather than $L^2(\R^d)$, since $E_{\mathcal{A}} f$ 
does not make sense for $f\in L^2(\R^d)$ in general when $A$ is nonsymmetric and nonsmooth, as stated in the introduction. The construction of $E_{\mathcal{A}}$ in Section \ref{subsec.extension.operator} shows $E_{\mathcal{A}}f\in \dot{H}^1(\R^{d+1}_+)$, so we have $\partial_t E_{\mathcal{A}} f\in L^2 (\R_+; L^2(\R^d))$. Hence, $E_{\mathcal{A}} f\in C (\overline{\R_+}; L^2(\R^d))$ when $u\in H^{1/2}(\R^d)$. In particular, if $f\in {H}^{1/2} (\R^d)$ then \eqref{def.poisson.semigroup} is well defined for all $t>0$. The main observation of this section is the analyticity of  $\{E_{\mathcal{A}}(t)\}_{t\geq 0}$  in ${H}^{1/2} (\R^d)$.

\begin{prop}\label{prop.poisson.analytic} The one-parameter family $\{E_{\mathcal{A}} (t)\}_{t\geq 0}$ defined by \eqref{def.poisson.semigroup}  is a strong continuous and analytic semigroup in ${H}^{1/2} (\R^d)$. Moreover, we have 
\[
\displaystyle \lim_{t\rightarrow \infty} t^{-\frac{1-k}{2}} \| E_{\mathcal{A}} (t) f \|_{\dot{H}^{\frac{k}{2}}(\R^d)} =0,~~~~~~~~~k=0,1,2.
\]
Here $\dot{H}^{0}(\R^d)$ is interpreted  as $L^2 (\R^d)$. 

\end{prop}

\begin{rem}{\rm From the proof below we will see that  $\| E_{\mathcal{A}} (t) f \|_{\dot{H}^{1/2}}\leq C \| f \|_{\dot{H}^{1/2}}$ holds for $t>0$, i.e., $\{E_{\mathcal{A}} (t)\}_{t\geq 0}$ is  bounded in the homogeneous space $\dot{H}^{1/2} (\R^d)$.

}
\end{rem}

\noindent {\it Proof of Proposition \ref{prop.poisson.analytic}.} (i) Boundedness: Let $s>0$ and set $v(\cdot,t) = \big ( E_{\mathcal{A}}f \big ) (\cdot, t+s)$ for $t\geq 0$. Then the variational characterization of $\|\cdot \|_{\dot{H}^{1/2}(\R^d)}$ implies
\begin{align}
\| E_{\mathcal{A}} (s) f \|_{\dot{H}^{\frac12}(\R^d)}^2  \leq \| \nabla v \|_{L^2(0,\infty; L^2 (\R^d))} ^2  &  = \int_s^\infty \|\nabla  E_{\mathcal{A}} f (\tau )\|_{L^2(\R^d)}^2 \dd \tau \label{est.proof.prop.poisson.analytic.1'} \\
& \leq \int_0^\infty \| \nabla E_{\mathcal{A}} f (\tau )\|_{L^2(\R^d)}^2 \dd \tau \leq C \| f \|_{\dot{H}^{\frac12} (\R^d)}^2.\label{est.proof.prop.poisson.analytic.1}
\end{align}

\noindent On the other hand, we have 
\begin{align}
\| E_{\mathcal{A}} (s) f -f \|_{L^2(\R^d)}  \leq \int_0^s \|\partial_t E_{\mathcal{A}} f (t) \|_{L^2(\R^d)}\dd t \leq  C s^\frac12 \| f \|_{\dot{H}^{\frac12}(\R^d)}.\label{est.proof.prop.poisson.analytic.2}
\end{align}

\noindent  (ii) Semigroup property: Let $s>0$. Since $E_{\mathcal{A}} (s) f\in H^{1/2}(\R^d)$ we can define $E_{\mathcal{A}} E_{\mathcal{A}} (s) f$. Set $w(t) = E_{\mathcal{A}} (t+s) f -   E_{\mathcal{A}} (t) E_{\mathcal{A}} (s) f$. Then, by recalling that $f$ and  $E_{\mathcal{A}} (s)f$ belong to $H^{1/2}(\R^d)$ it is easy to see that $w\in \dot{H}^1(\R^{d+1}_+)$, $\lim_{t\downarrow 0} w(t)=0$ in $L^2 (\R^d)$, and $\langle A\nabla w, \nabla \varphi \rangle _{L^2(\R^{d+1}_+)}=0$ for all $\varphi\in \dot{H}^1_0 (\R^{d+1}_+)$. In particular, we may take $\varphi = w$, which leads to $w=0$ by the coercive condition. Hence $E_{\mathcal{A}} (t+s) f = E_{\mathcal{A}} (t) E_{\mathcal{A}} (s) f$ holds for all $t,s>0$ if $f\in H^{1/2}(\R^d)$.

\noindent (iii) Time derivative estimate: Our aim is to show $\partial_t E_{\mathcal{A}} f\in H^1 (\R^{d+1}_{\delta,+})$ for all $\delta>0$ and 
\begin{align}
\| \partial_t E_{\mathcal{A}}  (t) f  \|_{H^\frac12 (\R^d)}\leq C t^{-1} (1+t^\frac12) \| f \|_{\dot{H}^\frac12 (\R^d )},~~~~~~~~~~t>0.\label{est.analytic.H.12}
\end{align}
\noindent But \eqref{est.analytic.H.12} immediately follows once we obtain  $\partial_t E_{\mathcal{A}} f\in H^1 (\R^{d+1}_{\delta,+})$ for all $\delta>0$. Indeed, from this regularity we see $\langle A \nabla \partial_t E_{\mathcal{A}} f, \nabla \varphi \rangle _{L^2(\R^{d+1}_+)} =0$ holds for all $\varphi\in C_0^\infty (\R^{d+1}_+)$ and thus we can take $\varphi = \phi_{t_0} (t) \partial_t E_{\mathcal{A}}f$ by the density argument, where $\phi_{r}(t)$, $r>0$, is a smooth nonnegative cut-off function such that $\phi_{r}(t)=0$ if $0\leq t\leq r/2$ and $\phi_{r} (t) = 1$ if $t\geq r$.  We can take $\phi_{r}$ so that $|\phi_{r}'| \leq C r^{-1} \phi_{r}^{1/2}$. Set $\phi_r^c=1-\phi_r$.  Then the term $|\langle A\nabla \partial_t E_{\mathcal{A}} f, \phi_{t_0} \nabla \partial_t E_{\mathcal{A}} f\rangle |$ is bounded from above by $C t_0^{-1} \| \phi_{t_0}^{1/2} \nabla \partial_t E_{\mathcal{A}} f \|_{L^2 (\R^{d+1}_+)} \| \partial_t E_{\mathcal{A}} f \|_{L^2 (\R^{d+1}_+)}$, which leads to 
\begin{align}
\| \partial_t E_{\mathcal{A}} (t_0) f\|_{\dot{H}^{\frac12} (\R^d)} \leq \| \nabla \partial_t E_{\mathcal{A}} f \|_{L^2 (\R^{d+1}_{t_0,+})} \leq C t_0^{-1} \| f \|_{\dot{H}^\frac12 (\R^d) }.\label{est.proof.prop.poisson.analytic.3}
\end{align}
The $L^2$ norm of $ \partial_t E_{\mathcal{A}} f (t_0)$ is estimated from the equality 
\begin{align*}
\langle \partial_t E_{\mathcal{A}} (t_0)  f , \varphi \rangle _{L^2(\R^d)} = - \langle \partial_t^2 E_{\mathcal{A}} f , \phi_{2 t_0}^c \varphi \rangle _{L^2 (\R^{d+1}_{t_0,+})} - \langle \partial_t E_{\mathcal{A}} f, \partial_t \phi_{2 t_0}^c\varphi \rangle _{L^2(\R^{d+1}_{t_0,+})}
\end{align*}
for $\varphi\in L^2(\R^d)$. By the duality we obtain 
\begin{align}
\| \partial_t E_{\mathcal{A}} (t_0) f  \|_{L^2 (\R^d)} \leq C t_0^\frac12 \| \partial_t^2 E_{\mathcal{A}} f \| _{L^2 (\R^{d+1}_{t_0,+})} + Ct_0^{-\frac12} \| \partial_t E_{\mathcal{A}} f \|_{L^2(\R^{d+1}_+)} \leq C t_0^{-\frac12} \| f \|_{\dot{H}^{\frac12}(\R^d)}. \label{est.proof.prop.poisson.analytic.4}
\end{align}
In the same manner we also have 
\begin{align}
\| \nabla_x E_{\mathcal{A}} (t_0) f  \|_{L^2 (\R^d)} \leq C t_0^\frac12 \| \partial_t \nabla_x E_{\mathcal{A}} f \| _{L^2 (\R^{d+1}_{t_0,+})} + Ct_0^{-\frac12} \| \nabla_x E_{\mathcal{A}} f \|_{L^2(\R^{d+1}_+)} \leq C t_0^{-\frac12} \| f \|_{\dot{H}^{\frac12}(\R^d)}. \label{est.proof.prop.poisson.analytic.4'}
\end{align}

Now the fact $\partial_t E_{\mathcal{A}} f\in H^1 (\R^{d+1}_{\delta,+})$ for all $\delta>0$ is shown by the standard argument, that is,  the  estimate of  the difference $h^{-1} \big ( E_{\mathcal{A}} (t+h) f  - E_{\mathcal{A}} (t) f   \big )$ and the limiting procedure as $h\rightarrow 0$. The assumption that $A$ is $t$-independent is essential here. We omit the details. By the bootstrap argument we also have $\partial_t^k E_{\mathcal{A}} f \in H^1 (\R^{d+1}_{\delta,+})$ for all $\delta>0$ and $k\in \N$.

\noindent (iv) Strong continuity: It suffices to show that $\lim_{t\downarrow  0} E_{\mathcal{A}} (t) f = f$ in $\dot{H}^{1/2}(\R^d)$, for the convergence in $L^2 (\R^d)$ follows from \eqref{est.proof.prop.poisson.analytic.2}. Let us recall that $E_{\mathcal{A}} f$ is constructed in the form $E_{\mathcal{A}} f= E_{-\Delta} f + w$ with $w = \mathcal{A}_D^{-1} \nabla \cdot \big ( (A-I) \nabla E_{-\Delta} f \big )\in \dot{H}^1_0 (\R^{d+1}_+)$. Since it is straightforward to see $\lim_{t\downarrow 0} E_{-\Delta} (t) f  = f$ in ${H}^{1/2}(\R^d)$, it suffices to prove $\| w (t_0) \|_{\dot{H}^{1/2}(\R^d)}\rightarrow 0$ as $t_0\rightarrow + 0$. But from the inequality 
\begin{align}
\| w (t_0)\|_{\dot{H}^\frac12(\R^d)} \leq \| \nabla \phi_{2 t_0}^c w \|_{L^2 (\R^{d+1}_{t_0,+})} & \leq  \| \phi_{2 t_0}^c \nabla w \|_{L^2 (\R^{d+1}_{t_0,+})} +  \| w  \partial_t \phi_{2 t_0}  \|_{L^2 (\R^{d+1}_{t_0,+})} \nonumber \\
& \leq C \| \nabla w \|_{L^2(0,2 t_0; L^2 (\R^d ) )},\label{strong.continuity}
\end{align}
we get the desired convergence.

The properties (i) - (iv) are sufficient conditions so that $\{ E_{\mathcal{A}} (t)\}_{t\geq 0}$ has an analytic extension in the sector $S=\{\lambda\in \C~|~\lambda\ne 0,~|{\rm arg}~\lambda |<\theta-\pi/2\}$ for some $\theta\in (\pi/2,\pi)$ (see \cite[Proposition 2.1.9]{Lunardi}), and thus, $\{ E_{\mathcal{A}} (t)\}_{t\geq 0}$  is an analytic semigroup in $H^{1/2}(\R^d)$.

Finally we observe that \eqref{est.proof.prop.poisson.analytic.1'} yields $\displaystyle \lim_{t\rightarrow \infty} \| E_{\mathcal{A}} (t) f  \|_{\dot{H}^{1/2} (\R^d)} =0$. Since \eqref{est.proof.prop.poisson.analytic.1} implies $\| E_{\mathcal{A}} (t+t_0) f \|_{L^2 (\R^d)}\leq \| E_{\mathcal{A}} (t_0) f \|_{L^2 (\R^d)} + C t^{1/2} \| E_{\mathcal{A}} (t_0) f \|_{\dot{H}^{1/2}(\R^d)}$, the above convergence in $\dot{H}^{1/2}(\R^d)$ leads to $\displaystyle \lim_{t\rightarrow \infty} t^{-1/2} \| E_{\mathcal{A}} (t)  f  \|_{L^2 (\R^d)} =0$.    Similarly, \eqref{est.proof.prop.poisson.analytic.4'} implies the bound  $\|\nabla_x E_{\mathcal{A}} (2 t_0) f  \|_{L^2(\R^d)} \leq C t_0^{-1/2} \|E_{\mathcal{A}} (t_0) f  \|_{\dot{H}^{1/2} (\R^d)}$, which proves $\displaystyle \lim_{t_0\rightarrow \infty} t_0^{1/2} \| \nabla_x E_{\mathcal{A}} (t_0) f \|_{L^2 (\R^d)} =0$.  The proof is complete.

\begin{rem}{\rm The generator of  $\{E_{\mathcal{A}}(t)\}_{t\geq 0}$ in $H^{1/2}(\R^d)$ is denoted by $-\mathcal{P}_{\mathcal{A}}$, and $\mathcal{P}_{\mathcal{A}}$ is called the {\it Poisson operator} associated with $\mathcal{A}$. Since $\{e^{-t\mathcal{P}_{\mathcal{A}}}\}_{t\geq 0}$ is strongly continuous in $H^{1/2}(\R^d)$, $D_{H^{1/2}}(\mathcal{P}_{\mathcal{A}})$ is dense in $H^{1/2}(\R^d)$. Moreover, the above proof implies $\|e^{-t\mathcal{P}_{\mathcal{A}}}  f \|_{H^{1/2}(\R^d)}\leq C(1+t^{1/2}) \| f \|_{H^{1/2}(\R^d)}$. Thus,  the spectrum of $\mathcal{P}_{\mathcal{A}}$ in $H^{1/2}(\R^d)$, denoted by $\sigma (\mathcal{P}_{\mathcal{A}})$, is included in  the closure of the sector  $\{ \lambda \in \C~|~\lambda\ne 0,~ |{\rm arg}~\lambda |<\theta\}$ for some $\theta\in (0,\pi/2)$. In the sequel we will freely use the maximal regularity for sectorial operators (with sectorial angles less than $\pi/2$) in the Hilbert space setting;  cf. \cite[Chapter 9.3.3]{Haase}.
}
\end{rem}

\begin{rem}\label{rem.prop.poisson.analytic.1}{\rm By taking  $v= \phi_{t_0} \phi_M^c E_{\mathcal{A}} f $, $M\gg t_0>0$, in $\langle A\nabla E_{\mathcal{A}} f, \nabla v\rangle_{L^2(\R^{d+1}_+)} =0$ and then by letting $M\rightarrow \infty$, it is not difficult to see $\|\nabla E_{\mathcal{A}} f \|_{L^2(t_0,\infty; L^2(\R^{d}))} \leq Ct_0^{-1}\| E_{\mathcal{A}} f \|_{L^2 (t_0/2,t_0; L^2(\R^{d}))}$. Since the variational characterization of the norm  $\|\cdot \|_{\dot{H}^{1/2}(\R^d)}$ yields the estimate $\| E_{\mathcal{A}} (t_0) f  \|_{\dot{H}^{1/2} (\R^d)} \leq  \|\nabla E_{\mathcal{A}} f \|_{L^2 (t_0,\infty; L^2(\R^{d}))}$ we have 
\begin{align}
\| e^{-t_0\mathcal{P}_{\mathcal{A}}}  f \|_{\dot{H}^{\frac12}(\R^d)}\leq C t_0^{-1} \| e^{-\cdot \mathcal{P}_{\mathcal{A}}}  f \|_{L^2 (\frac{t_0}{2},t_0; L^2 (\R^d ) )}.\label{est.rem.prop.poisson.analytic.1}
\end{align}
Here the constant $C$ depends only on $\nu_1$ and $\nu_2$. This estimate will be used later.
}
\end{rem}

The following proposition will be used in Section \ref{subsec.poisson.regular}.
\begin{prop}\label{prop.key.estimate.poisson.analytic} Assume that  $\beta\in [0,1)$ and $f\in H^{1/2}(\R^d)$, and let  $t>0$. Then there exists $C_\beta>0$ depending only on $d$, $\beta$, $\nu_1$, and $\nu_2$ such that
\begin{align}
& ~~~ \int_0^t \| e^{-s\mathcal{P}_{\mathcal{A}}} f \|_{\dot{H}^{\frac{1+\beta}{2}} (\R^d)} \| e^{-s\mathcal{P}_{\mathcal{A}}} f \|_{\dot{H}^{\frac{1-\beta}{2}}(\R^d)}  \dd s  + \int_0^t s^{-\beta} \| e^{-s\mathcal{P}_{\mathcal{A}}} f \|_{\dot{H}^\frac{1-\beta}{2} (\R^d)}^2 \dd s \nonumber \\
& \leq C_\beta \bigg ( t^{1-\beta}  \|  e^{-t\mathcal{P}_{\mathcal{A}}} f \|_{L^2 (\R^d)}^{2\beta} \|  e^{-t\mathcal{P}_{\mathcal{A}}} f \|_{\dot{H}^\frac{1}{2} (\R^d)}^{2(1-\beta)}  + \int_0^t  \|  e^{-s\mathcal{P}_{\mathcal{A}}} f \|_{\dot{H}^\frac12 (\R^d)}^{2} \dd s \bigg ). \label{est.prop.key.estimate.poisson.analytic}  
\end{align}

\end{prop}

\noindent {\it Proof.} From \eqref{est.proof.prop.poisson.analytic.1} and \eqref{est.proof.prop.poisson.analytic.4'} the interpolation inequality implies that 
\[
\| e^{-s\mathcal{P}_{\mathcal{A}}} f\|_{\dot{H}^{\frac{1+\beta}{2}} (\R^d)} \leq C s^{-\frac{\beta}{2}} \| e^{-\frac{s}{2} \mathcal{P}_{\mathcal{A}}} f\|_{\dot{H}^{\frac12} (\R^d)}.
\]
Then the Young inequality and the integration by parts lead to 
\begin{align}
& ~~~ 2 \int_0^t s^{-\frac{\beta}{2}}\| e^{-\frac{s}{2}\mathcal{P}_{\mathcal{A}}} f\|_{\dot{H}^\frac12 (\R^d)} \| e^{-s\mathcal{P}_{\mathcal{A}}} f\|_{\dot{H}^{\frac{1-\beta}{2}} (\R^d)} \dd s \nonumber \\
& \leq \int_0^t s^{-\beta}  \| e^{-s\mathcal{P}_{\mathcal{A}}} f\|_{\dot{H}^{\frac{1-\beta}{2}} (\R^d)}^2 \dd s + \int_0^{t} \| e^{-\frac{s}{2} \mathcal{P}_{\mathcal{A}}} f\|_{\dot{H}^{\frac12} (\R^d)}^2 \dd s \label{proof.prop.key.estimate.poisson.analytic.1} \\
& =   \frac{1}{1-\beta} \int_0^t (s^{1-\beta} )^{'}  \| e^{-s\mathcal{P}_{\mathcal{A}}} f\|_{\dot{H}^{\frac{1-\beta}{2}} (\R^d)}^2 \dd s + \int_0^{\frac{t}{2}} \| e^{-s \mathcal{P}_{\mathcal{A}}} f\|_{\dot{H}^{\frac12} (\R^d)}^2 \dd s \nonumber  \\
& = \frac{t^{1-\beta}}{1-\beta} \| e^{-t \mathcal{P}_{\mathcal{A}}} f \|_{\dot{H}^{\frac{1-\beta}{2}} (\R^d)}^2 - \frac{2}{1-\beta} \int_0^t s^{1-\beta}  {\rm Re} \langle (-\Delta_x)^{\frac{1-\beta}{2}} e^{-s\mathcal{P}_{\mathcal{A}}} f, \frac{\dd}{\dd s} e^{-s\mathcal{P}_{\mathcal{A}}} f \rangle _{L^2 (\R^d)} \dd s \nonumber \\
& ~~~~~ + \int_0^{\frac{t}{2}} \| e^{-s \mathcal{P}_{\mathcal{A}}} f\|_{\dot{H}^{\frac12} (\R^d)}^2 \dd s \nonumber \\
& \leq  \frac{t^{1-\beta}}{1-\beta} \| e^{-t \mathcal{P}_{\mathcal{A}}} f \|_{L^2 (\R)}^{2\beta} \| e^{-t \mathcal{P}_{\mathcal{A}}} f \|_{\dot{H}^\frac12 (\R)}^{2 (1-\beta)} \nonumber \\
& ~~~~~~~ +  \frac{C}{1-\beta} \int_0^t s^{\frac12 -\beta} \| e^{-s \mathcal{P}_{\mathcal{A}}} f \|_{\dot{H}^{1-\beta} (\R^d)}\|  e^{- \frac{s}{2} \mathcal{P}_{\mathcal{A}}} f \|_{\dot{H}^\frac12 (\R^d)}\dd s + \int_0^{t} \| e^{-s \mathcal{P}_{\mathcal{A}}} f\|_{\dot{H}^{\frac12} (\R)}^2 \dd s.\label{proof.prop.key.estimate.poisson.analytic.2}  
\end{align}
Here we have used \eqref{est.proof.prop.poisson.analytic.4}. By the interpolation 
\[
\| e^{-t \mathcal{P}_{\mathcal{A}}} f \|_{\dot{H}^{1-\beta} (\R^d)}\leq \| e^{-t \mathcal{P}_{\mathcal{A}}} f \|_{\dot{H}^{\frac{1-\beta}{2}} (\R^d)}^{\frac{2\beta}{1+\beta}} \| e^{-t \mathcal{P}_{\mathcal{A}}} f \|_{\dot{H}^{1} (\R^d)}^{\frac{1-\beta}{1+\beta}}
\]
and \eqref{est.proof.prop.poisson.analytic.4'} we see that 
\begin{align}
& ~~~ \frac{C}{1-\beta} \int_0^t s^{\frac12 -\beta} \| e^{-s \mathcal{P}_{\mathcal{A}}} f \|_{\dot{H}^{1-\beta} (\R^d)}\|  e^{- \frac{s}{2} \mathcal{P}_{\mathcal{A}}} f \|_{\dot{H}^\frac12 (\R^d)}\dd s  \nonumber \\
&  \leq \int_0^t s^{-\frac{\beta}{2}}\| e^{-s \mathcal{P}_{\mathcal{A}}} f\|_{\dot{H}^\frac12 (\R^d)} \| e^{-s \mathcal{P}_{\mathcal{A}}} f\|_{\dot{H}^{\frac{1-\beta}{2}} (\R^d)} \dd s  + C_\beta \int_0^t  s^\frac12   \| e^{-s \mathcal{P}_{\mathcal{A}}} f \|_{\dot{H}^{1} (\R^d)}\|  e^{- \frac{s}{2} \mathcal{P}_{\mathcal{A}}} f \|_{\dot{H}^\frac12 (\R^d)}\dd s \nonumber \\
& \leq  \int_0^t s^{-\frac{\beta}{2}}\| e^{-s\mathcal{P}_{\mathcal{A}}} f\|_{\dot{H}^\frac12 (\R^d)} \| e^{-s\mathcal{P}_{\mathcal{A}}} f\|_{\dot{H}^{\frac{1-\beta}{2}} (\R^d)} \dd s  + C_\beta \int_0^t   \|  e^{- \frac{s}{2} \mathcal{P}_{\mathcal{A}}} f \|_{\dot{H}^\frac12 (\R^d)}^2\dd s.\label{proof.prop.key.estimate.poisson.analytic.3}
\end{align}
Here $C_\beta$ depends only on $d$, $\nu_1$, $\nu_2$, and $\beta$. Now it is easy to derive  \eqref{est.prop.key.estimate.poisson.analytic} from \eqref{proof.prop.key.estimate.poisson.analytic.1} - \eqref{proof.prop.key.estimate.poisson.analytic.3}. The proof is complete.

\subsection{Poisson semigroup in $L^2 (\R^d)$}

In general, the Poisson semigroup $\{e^{-t\mathcal{P}_{\mathcal{A}}}\}_{t\geq 0}$ acting on $H^{1/2}(\R^d)$ is not extended as a strongly continuous semigroup in $L^2 (\R^d)$. However, a simple Caccioppoli type inequality shows that once we have  the bound of $\{e^{-t\mathcal{P}_{\mathcal{A}}}\}_{t\geq 0}$ in $L^2 (\R^d)$, then the analyticity in $L^2 (\R^d)$ automatically follows.

\begin{prop}\label{prop.poisson.analytic.L^2} Assume that  the Poisson semigroup $\{e^{-t\mathcal{P}_{\mathcal{A}}}\}_{t\geq 0}$ in $H^{1/2}(\R^d)$ is extended as a strongly continuous semigroup in $L^2 (\R^d)$. Then it is an analytic semigroup acting on $L^2 (\R^d)$, and it follows that $\| e^{-t\mathcal{P}_{\mathcal{A}}} \|_{L^2\rightarrow L^2} \leq C  (1+t )^{1/2}$ and 
\begin{equation}
t \| \frac{\dd }{\dd t} e^{-t\mathcal{P}_{\mathcal{A}}} \|_{L^2\rightarrow L^2} + t^\beta \| e^{-t\mathcal{P}_{\mathcal{A}}} \|_{L^2\rightarrow \dot{H}^{\beta}} \leq C' \sup_{\frac{t}{4}< \tau <t} \| e^{-\tau\mathcal{P}_{\mathcal{A}}}  \|_{L^2\rightarrow L^2}~~~~~~~~~t>0,~~~\beta\in [0,1]. \label{est.prop.poisson.analytic.L^2}
\end{equation}
Moreover, we have $\displaystyle \lim_{t\rightarrow \infty} t^{\beta-1/2} \| e^{-t\mathcal{P}_{\mathcal{A}}} f \|_{\dot{H}^{\beta}(\R^d)} =0$, $\beta\in [0,1]$, for any $f\in L^2 (\R^d)$.
The constant $C'$ in \eqref{est.prop.poisson.analytic.L^2} depends only on $\nu_1$ and $\nu_2$.
\end{prop}

\begin{rem}\label{rem.prop.poisson.analytic.L^2.1}{\rm If $\{e^{-t\mathcal{P}_{\mathcal{A}}}  \}_{t\geq 0}$ is a bounded semigroup in $L^2(\R^d)$ then \eqref{est.prop.poisson.analytic.L^2} implies the estimate  $t \|\dd/\dd t~e^{-t\mathcal{P}_{\mathcal{A}}} \|_{L^2\rightarrow L^2} + t^\beta \|  e^{-t\mathcal{P}_{\mathcal{A}}} \|_{L^2\rightarrow \dot{H}^{\beta}} \leq C$ for $t>0$ and $\beta\in [0,1]$. 

}
\end{rem}

\begin{rem}\label{rem.prop.poisson.analytic.L^2}{\rm By the density argument,  $\{ e^{-t\mathcal{P}_{\mathcal{A}}} \}_{t\geq 0}$  in $H^{1/2}(\R^d)$ is extended as a strongly continuous semigroup in $L^2 (\R^d)$ if we have the estimate $\|e^{-t\mathcal{P}_{\mathcal{A}}} f\|_{L^2(\R^d)} \leq C \| f \|_{L^2 (\R^d)}$ for  $t\in (0,1]$ and $f\in H^{1/2}(\R^d)$. 

}
\end{rem}

\noindent {\it Proof of Proposition \ref{prop.poisson.analytic.L^2}.}  By the assumption of the proposition, for $f \in L^2(\R^d)$ we have $e^{-t\mathcal{P}_{\mathcal{A}}} f = \displaystyle \lim_{n\rightarrow \infty} e^{-t\mathcal{P}_{\mathcal{A}}} f_n$ in $L^2 (\R^d)$, where $\{f_n\}_{n\in \N}\subset H^{1/2}(\R^d)$ and $f_n\rightarrow f$ in $L^2 (\R^d)$ as $n\rightarrow \infty$. Hence, to prove the analyticity it suffices to show the estimate
\begin{equation}
\| \frac{\dd}{\dd t} e^{-t\mathcal{P}_{\mathcal{A}}} f \|_{L^2 (\R^d)}  \leq C t^{-1} \| f \|_{L^2 (\R^d)},~~~~~~~~~t \in (0,1],~~~f\in H^{\frac12} (\R^d).\label{est.proof.prop.poisson.analytic.L^2.1}
\end{equation}
Since $\{e^{-t\mathcal{P}_{\mathcal{A}}}\}_{t\geq 0}$ is assumed to be strongly continuous in $L^2(\R^d)$, we have from \eqref{est.rem.prop.poisson.analytic.1}
\begin{align}
\| e^{-t_0\mathcal{P}_{\mathcal{A}}} f \|_{\dot{H}^{\frac12} (\R^d)} \leq  C t_0^{-\frac12}  \sup_{\frac{t_0}{2}<\tau <t_0} \| e^{-\tau\mathcal{P}_{\mathcal{A}}}\|_{L^2\rightarrow L^2} \| f \|_{L^2 (\R^d)}.\label{est.proof.prop.poisson.analytic.L^2.2}
\end{align}
Combining \eqref{est.proof.prop.poisson.analytic.4} with \eqref{est.proof.prop.poisson.analytic.L^2.2}, we get 
\begin{align*}
\|\frac{\dd}{\dd t} e^{-2 t_0 \mathcal{P}_{\mathcal{A}}} f  \|_{L^2 (\R^d)} \leq Ct_0^{-\frac12} \| e^{-t_0 \mathcal{P}_{\mathcal{A}}} f   \|_{\dot{H}^\frac12 (\R^d)} \leq C t_0^{-1}  \sup_{\frac{t_0}{2}<\tau < t_0} \| e^{-\tau\mathcal{P}_{\mathcal{A}}}  \|_{L^2\rightarrow L^2} \| f \|_{L^2 (\R^d)}.
\end{align*}
This is nothing but \eqref{est.proof.prop.poisson.analytic.L^2.1}, as desired. The estimate for $\| e^{-t\mathcal{P}_{\mathcal{A}}} \|_{L^2\rightarrow L^2}$ follows from \eqref{est.proof.prop.poisson.analytic.2} and \eqref{est.proof.prop.poisson.analytic.L^2.2}. Indeed, we have $\| e^{-(t+\tau )\mathcal{P}_{\mathcal{A}}} f \|_{L^2 (\R^d) } \leq  \| e^{-\tau\mathcal{P}_{\mathcal{A}}} f \|_{L^2 (\R^d) } +  C t^{1/2} \| e^{-\tau\mathcal{P}_{\mathcal{A}}} f \|_{\dot{H}^{1/2} (\R^d)}$ for $t,\tau>0$. This proves $\|e^{-t\mathcal{P}_{\mathcal{A}}} \|_{L^2\rightarrow L^2} \leq C(1+t^{1/2})$ by taking $\tau=1$. The estimate \eqref{est.prop.poisson.analytic.L^2} with $\beta=1$ follows from \eqref{est.proof.prop.poisson.analytic.4'} and  \eqref{est.proof.prop.poisson.analytic.L^2.2}. The case $\beta\in (0,1)$ is a consequence of the interpolation inequality. The convergence as $t\rightarrow \infty$ is already shown in Proposition \ref{prop.poisson.analytic}. The proof is complete.

\vspace{0.5cm}

If the Poisson semigroup  in $H^{1/2}(\R^d)$ is extended as a semigroup in $L^2 (\R^d)$ it is important to determine the domain of $D_{L^2}(\mathcal{P}_{\mathcal{A}})$, although it is a difficult problem in general. In the proceeding section we will study the relation between $D_{L^2}(\mathcal{P}_{\mathcal{A}})$ and $H^1 (\R^d)$. We note that, once the characterization  $D_{L^2}(\mathcal{P}_{\mathcal{A}}) = H^1 (\R^d)$ is verified,  the Poisson semigroup  is analytic in $H^s (\R^d)$ for each $s\in (0,1]$ by the interpolation.

\subsection{General relation between Dirichlet-Neumann map and Poisson operator}\label{subsec.relation.Lambda.P}

In this section we consider the relation between $\Lambda_{\mathcal{A}}$ and $\mathcal{P}_{\mathcal{A}}$, which is formally described as $\mathcal{P}_{\mathcal{A}} =  M_{1/b}\Lambda_{\mathcal{A}} + M_{{\bf r_2}/b}\cdot \nabla_x$ if one does not take into account the relation among the domains of each operator. For $f\in H^{1/2}(\R^d)$ set $\tilde \Lambda_{\mathcal{A}} f := - M_b \partial_t E_{\mathcal{A}} f - M_{{\bf r_2}}\cdot \nabla_x E_{\mathcal{A}} f \in L^2(\R^{d+1}_+)$.  Then the identity $\langle A\nabla E_{\mathcal{A}} f, \nabla \varphi\rangle =0$ for $\varphi\in C_0^\infty (\R^{d+1}_+)$ yields 
\begin{align}
\langle \partial_t \tilde \Lambda_{\mathcal{A}} f , \varphi \rangle _{L^2(\R^{d+1}_+)} = - \sum_{1\leq i\leq d, 1\leq j\leq d+1} \langle a_{i,j} \partial_j E_{\mathcal{A}} f, \partial_i \varphi\rangle _{L^2(\R^{d+1}_+)}, ~~~~~~~\varphi\in C_0^\infty (\R^{d+1}_+), \label{eq.subsec.relation.Lambda.P.1}
\end{align}
where we have used that $\partial_t^k E_{\mathcal{A}} f \in H^1(\R^{d+1}_{\delta,+})$ for all $\delta>0$ and $k\in \N$ when $f \in H^{1/2}(\R^d)$. Note that the right-hand side of \eqref{eq.subsec.relation.Lambda.P.1} makes sense for all $\varphi\in L^2 (\R_+; \dot{H}^1(\R^d))$. Hence $\partial_t \tilde \Lambda_{\mathcal{A}} f$ is extended as a distribution belonging to $L^2 (\R_+; \dot{H}^{-1}(\R^d))$. In particular, the function $\tilde \Lambda_{\mathcal{A}} f$ belongs to $C([0,\infty); H^{-1}(\R^d))$, and  we have from \eqref{eq.subsec.relation.Lambda.P.1} with $\varphi\in C_0^\infty (\overline{\R^{d+1}_+})$,
\begin{align}
\langle \tilde \Lambda_{\mathcal{A}}  f |_{t=0}, \varphi|_{t=0}\rangle _{\dot{H}^{-1} (\R^d), \dot{H}^1(\R^d)} & = - \int_0^\infty \langle \partial_t\tilde \Lambda_{\mathcal{A}}  f, \varphi\rangle _{\dot{H}^{-1} (\R^d), \dot{H}^1(\R^d) }\dd t - \int_0^\infty \langle \tilde \Lambda_{\mathcal{A}} f , \partial_t\varphi\rangle _{L^2} \dd t\nonumber \\
& =  \sum_{1\leq i\leq d, 1\leq j\leq d+1} \langle a_{i,j} \partial_j E_{\mathcal{A}} f, \partial_i \varphi\rangle _{L^2(\R^{d+1}_+)} -   \langle \tilde \Lambda_{\mathcal{A}} f , \partial_t\varphi\rangle _{L^2(\R^{d+1}_+)} \nonumber \\
& = \langle \Lambda_{\mathcal{A}}  f, \varphi |_{t=0}\rangle _{\dot{H}^{-\frac12}(\R^d), \dot{H}^{\frac12}(\R^d)}.\label{eq.subsec.relation.Lambda.P.2}
\end{align}
Here we have used \eqref{dn.hosoku} in the last line. Since $E_{\mathcal{A}} (t) f  = e^{-t\mathcal{P}_{\mathcal{A}}} f$,  the equality \eqref{eq.subsec.relation.Lambda.P.2}  implies
\begin{prop}\label{prop.subsec.relation.1} Let $f\in H^{1/2}(\R^d)$. Then we have
\begin{align}
\lim_{t\downarrow 0} \big ( M_b \mathcal{P}_{\mathcal{A}}  e^{-t\mathcal{P}_{\mathcal{A}} } f  - M_{{\bf r_2}}\cdot \nabla_x e^{-t\mathcal{P}_{\mathcal{A}} } f \big ) = \Lambda_{\mathcal{A}} f ~~~~~~~~~{\rm in}~~H^{-1}(\R^d).\label{eq.prop.subsec.relation.1}
\end{align}
\end{prop}

When $f$ belongs to $D_{H^{1/2}}(\mathcal{P}_{\mathcal{A}})$  we have the desired representation as follows.
\begin{prop}\label{prop.subsec.relation.2} It follows that $D_{H^{1/2}} (\mathcal{P}_{\mathcal{A}}) \hookrightarrow  D_{L^2} (\Lambda_{\mathcal{A}})\cap H^1 (\R^d)$ and 
\begin{align}
\mathcal{P}_{\mathcal{A}}  f = M_{1/b} \Lambda_{\mathcal{A}}  f + M_{{\bf r_2}/b} \cdot \nabla_x f, ~~~~~~~~f\in D_{H^{\frac12}} (\mathcal{P}_{\mathcal{A}}). \label{eq.prop.subsec.relation.2}
\end{align}

\end{prop}

\noindent {\it Proof.} When $f \in D_{H^{1/2}} (\mathcal{P}_{\mathcal{A}})$ we have $\dd /\dd t ~ e^{-t\mathcal{P}_{\mathcal{A}}}f = - e^{-t \mathcal{P}_{\mathcal{A}}} \mathcal{P}_{\mathcal{A}} f = - (E_{\mathcal{A}} \mathcal{P}_{\mathcal{A}} f ) (\cdot,t)$, which belongs to $L^2 (\R_+; \dot{H}^1 (\R^{d}) )$. Hence, $e^{-t\mathcal{P}_{\mathcal{A}} } f$ belongs to $C(\overline{\R_+}; H^1 (\R^d))$,  and in particular, $f\in H^1(\R^d)$.  It is easy to see $\|f\|_{H^1(\R^d)}\leq C \| f\|_{D_{H^{1/2}}(\mathcal{P}_{\mathcal{A}})}$. Moreover, the left-hand side of \eqref{eq.prop.subsec.relation.1} converges to $M_b \mathcal{P}_{\mathcal{A}} f - M_{{\bf r_2}}\cdot \nabla_x f$ in $L^2(\R^d)$. This implies $f\in D_{L^2} (\Lambda_{\mathcal{A}} )$ and \eqref{eq.prop.subsec.relation.2} holds. The proof is complete.

\section{Identity of Rellich type and factorization of elliptic operator}\label{sec.proof.thm.factorization}

This section is devoted to derive the Rellich type identity with the aid of the class $D_{H^{1/2}}(\mathcal{P}_{\mathcal{A}})$, which is used in studying the relation among the domains of Poisson operators, Dirichlet-Neumann maps, and $H^1(\R^d)$. 

\subsection{Identity in weak form for general case}\label{subsec.thm.factorization}

In this section we give a proof of Theorem \ref{thm.factorization}. The following proposition plays central roles.
\begin{prop}\label{prop.thm.factorization} Let $f\in D_{H^{1/2}} (\mathcal{P}_{\mathcal{A}})$ and $g\in H^1 (\R^d)$. Then it follows that
\begin{align}
\langle \mathcal{A}' f, g\rangle _{\dot{H}^{-1},\dot{H}^1} = \langle  \mathcal{P}_{\mathcal{A}} f,  \Lambda_{\mathcal{A}^*} g \rangle _{\dot{H}^{\frac12}, \dot{H}^{-\frac12}} + \langle  \mathcal{P}_{\mathcal{A}} f,  M_{{\bf \bar{r}_1}} \cdot \nabla_x  g \rangle _{L^2(\R^d)}.\label{eq.prop.thm.factorization}
\end{align}
Here $\langle \mathcal{A}'f,g\rangle  _{\dot{H}^{-1} ,\dot{H}^1}= \langle A' \nabla_x f, \nabla_x g\rangle _{L^2(\R^d)}$ with $A'(x) = \big ( a_{i,j} (x) \big ) _{1\leq i,j \leq d}$. 

\end{prop}

\begin{rem}{\rm When $A$ is Hermite and $f=g$ the equality \eqref{eq.prop.thm.factorization}  is formally written as 
$$\langle A'\nabla_x f, \nabla_x f\rangle _{L^2(\R^d)} = \| M_{\sqrt{b}} \mathcal{P}_{\mathcal{A}} f \| _{L^2(\R^d)}^2,$$
which is a variant of the classical Rellich identity \cite{Rellich}.

}

\end{rem}

\noindent {\it Proof of Proposition \ref{prop.thm.factorization}.}  Set $w (t) = A' \nabla_x e^{-t\mathcal{P}_{\mathcal{A}}}  f  - M_{\bf r_1} \mathcal{P}_{\mathcal{A}}e^{-t\mathcal{P}_{\mathcal{A}}}  f$ for $f\in D_{H^{1/2}}(\mathcal{P}_{\mathcal{A}})$. Now let us recall the equality \eqref{eq.subsec.relation.Lambda.P.1}, which implies 
\begin{align*}
- \int_0^\infty \langle w(t), \nabla_x \varphi (t)  \rangle _{L^2(\R^d)} \dd t = \int_0^\infty \langle M_b \frac{\dd}{\dd t} e^{-t\mathcal{P}_{\mathcal{A}}} \mathcal{P}_{\mathcal{A}} f+ M_{\bf r_2} \cdot \nabla_x e^{-t\mathcal{P}_{\mathcal{A}}}  \mathcal{P}_{\mathcal{A}} f, \varphi (t) \rangle _{L^2(\R^d)} \dd t
\end{align*}
for all $\varphi\in C_0^\infty (\R^{d+1}_+)$. Since each of $w(t)$, $\frac{\dd}{\dd t} e^{-t\mathcal{P}_{\mathcal{A}}}  \mathcal{P}_{\mathcal{A}} f$, and $M_{\bf r_2} \cdot \nabla_x e^{-t\mathcal{P}_{\mathcal{A}}}  \mathcal{P}_{\mathcal{A}} f$ belongs to $C([\delta,\infty); L^2 (\R^d))$ for $\delta>0$, we have from Proposition \ref{prop.subsec.relation.2},
\begin{align}
- \langle w(t), \nabla_x g  \rangle _{L^2(\R^d)}  =  \langle M_b \frac{\dd}{\dd t} e^{-t\mathcal{P}_{\mathcal{A}}}  \mathcal{P}_{\mathcal{A}} f + M_{\bf r_2} \cdot \nabla_x e^{-t\mathcal{P}_{\mathcal{A}}} \mathcal{P}_{\mathcal{A}} f, g  \rangle _{L^2(\R^d)} =  - \langle \Lambda_{\mathcal{A}} e^{-t\mathcal{P}_{\mathcal{A}}} \mathcal{P}_{\mathcal{A}} f, g \rangle _{L^2(\R^d)} \label{proof.prop.thm.factorization.1}
\end{align}
for all $g \in C_0^\infty (\R^{d})$ and $t>0$. It is clear that \eqref{proof.prop.thm.factorization.1} holds for all $g\in H^1 (\R^d)$. As stated in the proof of Proposition \ref{prop.subsec.relation.2}, $e^{-t\mathcal{P}_{\mathcal{A}}} f\rightarrow f$ in $H^1(\R^d)$  as $t\rightarrow 0$ if $f\in D_{H^{1/2}}(\mathcal{P}_{\mathcal{A}})$. Moreover, we see $\langle \Lambda_{\mathcal{A}} e^{-t\mathcal{P}_{\mathcal{A}}} \mathcal{P}_{\mathcal{A}} f, g \rangle _{L^2(\R^d)}= \langle \Lambda_{\mathcal{A}} e^{-t\mathcal{P}_{\mathcal{A}}} \mathcal{P}_{\mathcal{A}} f, g \rangle _{\dot{H}^{-1/2}, \dot{H}^{1/2}}\rightarrow \langle \Lambda_{\mathcal{A}} \mathcal{P}_{\mathcal{A}} f, g \rangle _{\dot{H}^{-1/2},\dot{H}^{1/2}}$ as $t\rightarrow 0$ for $f\in D_{H^{1/2}}(\mathcal{P}_{\mathcal{A}})$ and $g\in H^{1/2}(\R^d)$. Hence, letting $t\rightarrow 0$ in \eqref{proof.prop.thm.factorization.1}, we get 
\begin{align}
- \langle  A' \nabla_x  f  - M_{\bf r_1} \mathcal{P}_{\mathcal{A}}  f, \nabla_x g  \rangle _{L^2(\R^d)}  =  - \langle \Lambda_{\mathcal{A}} \mathcal{P}_{\mathcal{A}} f, g \rangle _{\dot{H}^{-\frac12}, \dot{H}^\frac12} =- \langle \mathcal{P}_{\mathcal{A}} f, \Lambda_{\mathcal{A}^*} g \rangle _{\dot{H}^{\frac12}, \dot{H}^{-\frac12}} \label{proof.prop.thm.factorization.2}
\end{align}
if $f\in D_{H^{1/2}}(\mathcal{P}_{\mathcal{A}})$ and $g\in H^1(\R^d)$. This proves \eqref{eq.prop.thm.factorization}. The proof is complete.

\vspace{0.5cm}

\noindent {\it Proof of Theorem \ref{thm.factorization}.} Let us prove \eqref{eq.factorization}. Assume that $u\in L^2 (\R; D_{H^{1/2}}(\mathcal{P}_{\mathcal{A}}))\cap \dot{H}^1 (\R; L^2 (\R^d))$ and
$v\in L^2 (\R; D_{H^{1/2}}(\mathcal{P}_{\mathcal{A}^*}))\cap \dot{H}^1 (\R; L^2 (\R^d))$. Then  the definition of $\mathcal{A}$ and Proposition \ref{prop.thm.factorization} imply
\begin{align}
 \langle A\nabla u, \nabla v\rangle_{L^2 (\R^{d+1})}  & = \langle A' \nabla _x u, \nabla_x v\rangle _{L^2(\R^{d+1})} + \langle \partial_t u, M_{\bar{b}} \partial_t v \rangle_{L^2(\R^{d+1})}\nonumber \\
& ~~~ + \langle M_{{\bf r_2}}\cdot \nabla_x u, \partial_t v\rangle _{L^2(\R^{d+1})} + \langle \partial_t u, M_{{\bf \bar{r}_1}} \cdot \nabla_x v\rangle _{L^2(\R^{d+1})}\nonumber \\
& = \langle  \mathcal{P}_{\mathcal{A}} u, ~ M_{\bar{b}} \mathcal{P}_{\mathcal{A}^*} v  \rangle _{L^2(\R^{d+1})} +  \langle \partial_t u, M_{\bar{b}} \partial_t v \rangle_{L^2(\R^{d+1})} \nonumber \\
& ~~~ +  \langle M_{{\bf r_2}}\cdot \nabla_x u, \partial_t v\rangle _{L^2(\R^{d+1})} +  \langle \partial_t u, M_{{\bf \bar{r}_1}} \cdot \nabla_x v\rangle _{L^2(\R^{d+1})} \label{proof.thm.factorization.1} \\
& = \langle ~ \big ( \partial_t  + \mathcal{P}_{\mathcal{A}} \big ) u, ~M_{\bar{b}} \big ( \partial_t + \mathcal{P}_{\mathcal{A}^*} \big ) v ~ \rangle _{L^2(\R^{d+1})}\nonumber \\
& ~~~ -\langle \partial_t u,  \Lambda_{\mathcal{A}^*} v \rangle _{L^2(\R^{d+1})} - \langle M_{1/b} \Lambda_{\mathcal{A}} u, M_{\bar{b}} \partial_t v  \rangle_{L^2(\R^{d+1})}.\nonumber
\end{align}
Thus we get \eqref{eq.factorization} since $\langle M_{1/b} \Lambda_{\mathcal{A}} u, M_{\bar{b}} \partial_t v\rangle_{L^2(\R^{d+1})} = - \langle \partial_t u,  \Lambda_{\mathcal{A}^*} v\rangle _{L^2(\R^{d+1})}$ holds. 

Next we consider \eqref{eq.factorization'}. Assume that  $u\in L^2 (\R; D_{H^{1/2}} (\mathcal{P}_{\mathcal{A}} ))\cap \dot{H}^1(\R; H^{1/2}(\R^d))$ and $v\in H^1 (\R^{d+1})$. In this case, instead of \eqref{proof.thm.factorization.1}, we have from \eqref{eq.prop.thm.factorization},
\begin{align}
\langle A\nabla u, \nabla v\rangle_{L^2 (\R^{d+1})} & = \int_\R \langle  \mathcal{P}_{\mathcal{A}} u (t) , ~\Lambda_{\mathcal{A}^*} v (t)  \rangle _{\dot{H}^\frac12, \dot{H}^{-\frac12} } \dd t +  \langle  \mathcal{P}_{\mathcal{A}} u , M_{{\bf \bar{r}_1}}\cdot \nabla_x v \rangle _{L^2(\R^{d+1})} \nonumber \\
&~~ + \langle \partial_t u, M_{\bar{b}} \partial_t v \rangle_{L^2(\R^{d+1})}  +  \langle M_{{\bf r_2}}\cdot \nabla_x u, \partial_t v\rangle _{L^2(\R^{d+1})} +  \langle \partial_t u, M_{{\bf \bar{r}_1}} \cdot \nabla_x v\rangle _{L^2(\R^{d+1})} \nonumber \\
& =  \int_\R \langle(  \partial_t +  \mathcal{P}_{\mathcal{A}} ) u (t) , ~\Lambda_{\mathcal{A}^*} v (t)  \rangle _{\dot{H}^\frac12, \dot{H}^{-\frac12} } \dd t +  \langle  ( \partial_t + \mathcal{P}_{\mathcal{A}}) u , ~ (M_{\bar{b}} \partial_t  + M_{{\bf \bar{r}_1}}\cdot \nabla_x ) v \rangle _{L^2(\R^{d+1})} \nonumber \\
&~~-  \int_\R \langle  \partial_t u (t) , ~\Lambda_{\mathcal{A}^*} v (t)  \rangle _{\dot{H}^\frac12, \dot{H}^{-\frac12} } \dd t  - \int_\R \langle M_{1/b}\Lambda_{\mathcal{A}} u (t) , ~ M_{\bar{b}} \partial_t v (t) \rangle _{L^2(\R^d)} \dd t.\nonumber 
\end{align}
But the density argument implies 
\[
 \int_\R \langle  \partial_t u (t) , ~\Lambda_{\mathcal{A}^*} v (t)  \rangle _{\dot{H}^\frac12, \dot{H}^{-\frac12} } \dd t = -\int_\R \langle M_{1/b}\Lambda_{\mathcal{A}} u (t) , ~ M_{\bar{b}} \partial_t v (t) \rangle _{L^2(\R^d)} \dd t,
\]
which gives \eqref{eq.factorization'}. It remains to show \eqref{eq.factorization.half} and   \eqref{eq.thm.factorization.half}. We will give a proof only for  \eqref{eq.thm.factorization.half} since \eqref{eq.factorization.half}  is proved in the same manner. As in the proof of \eqref{eq.factorization'}, we have 
\begin{align}
\langle A \nabla u, \nabla v\rangle _{L^2 (\R^{d+1}_+)}  & =  \int_0^\infty \langle(  \partial_t +  \mathcal{P}_{\mathcal{A}} ) u (t) , ~\Lambda_{\mathcal{A}^*} v (t)  \rangle _{\dot{H}^\frac12, \dot{H}^{-\frac12} } \dd t \nonumber \\
& ~~~ +  \langle  ( \partial_t + \mathcal{P}_{\mathcal{A}}) u , ~ (M_{\bar{b}} \partial_t  + M_{{\bf \bar{r}_1}}\cdot \nabla_x ) v \rangle _{L^2(\R^{d+1}_+)} \nonumber \\
&~~~~ -  \int_0^\infty \langle  \partial_t u (t) , ~\Lambda_{\mathcal{A}^*} v (t)  \rangle _{\dot{H}^\frac12, \dot{H}^{-\frac12} } \dd t  - \int_0^\infty \langle \Lambda_{\mathcal{A}} u (t) ,  \partial_t v (t) \rangle _{L^2(\R^d)} \dd t.\nonumber 
\end{align}
Let $v_\epsilon = j_\epsilon *v$ be the mollification of $v$ with respect to the $x$ variables, which converges to $v$ in $H^1 (\R^{d+1}_+)$. Then we see from the integration by parts,
\begin{align*}
\int_0^\infty \langle  \partial_t u (t) , ~\Lambda_{\mathcal{A}^*} v (t)  \rangle _{\dot{H}^\frac12, \dot{H}^{-\frac12} } \dd t & =  \lim_{\epsilon\rightarrow 0} \int_0^\infty \langle  \partial_t u (t) , ~\Lambda_{\mathcal{A}^*} v_\epsilon (t)  \rangle _{\dot{H}^\frac12, \dot{H}^{-\frac12} } \dd t\\
& = \lim_{\epsilon\rightarrow 0} \langle \gamma u, \Lambda_{\mathcal{A}^*}  \gamma v_\epsilon  \rangle_{\dot{H}^\frac12, \dot{H}^{-\frac12} } - \lim_{\epsilon\rightarrow 0} \int_0^\infty \langle  \Lambda_{\mathcal{A}}u (t) , ~ \partial_t v_\epsilon (t)  \rangle _{L^2 (\R^d)} \dd t\\
& = \langle \gamma u, \Lambda_{\mathcal{A}^*} \gamma  v   \rangle_{\dot{H}^\frac12, \dot{H}^{-\frac12} } -  \int_0^\infty \langle  \Lambda_{\mathcal{A}}u (t) , ~ \partial_t v (t)  \rangle _{L^2 (\R^d)} \dd t.
\end{align*}
Here $\gamma: H^{1}(\R^{d+1}_+)\rightarrow H^{1/2}(\partial\R^{d+1}_+)$ is the trace operator. Hence  \eqref{eq.thm.factorization.half} holds. The proof is complete.

\subsection{Criterion on the embedding between $H^1(\R^d)$ and $D_{L^2}(\mathcal{P}_{\mathcal{A}})$}\label{subsec.embedding.poisson.H^1}

In this section we show various relations among $D_{L^2}(\mathcal{P}_{\mathcal{A}})$, $D_{L^2}(\Lambda_{\mathcal{A}})$, and $H^1 (\R^d)$.

\begin{prop}\label{prop.embedding.section3.1} The following two statements are equivalent.

\vspace{0.1cm}

\noindent {\rm (i)} $D_{H^{1/2}} (\mathcal{P}_{\mathcal{A}})\subset D_{L^2}(\Lambda_{\mathcal{A}^*})$ and $\| \Lambda_{\mathcal{A}^*} f\|_{L^2 (\R^d)} \leq C  \| f \|_{H^1 (\R^d)}$ holds for $f\in D_{H^{1/2}}(\mathcal{P}_{\mathcal{A}})$.

\noindent {\rm (ii)} $\{e^{-t\mathcal{P}_{\mathcal{A}}}\}_{t\geq 0}$ is extended as a strongly continuous semigroup in $L^2 (\R^d)$ and  $D_{L^2}(\mathcal{P}_{\mathcal{A}})$ is continuously embedded in $H^1 (\R^d)$. 

\vspace{0.1cm}

\noindent Moreover, if the condition {\rm (ii)} (and hence, {\rm (i)}) holds then $D_{L^2}(\mathcal{P}_{\mathcal{A}})$ is continuously embedded in $D_{L^2}(\Lambda_{\mathcal{A}})$,  $H^1 (\R^d)$ is continuously embedded in $D_{L^2}(\Lambda_{\mathcal{A}^*})$, and it follows that
\begin{align}
\mathcal{P}_{\mathcal{A}}  f   & = M_{1/b} \Lambda_{\mathcal{A}} f + M_{{\bf r_2}/b} \cdot \nabla_x f,\label{assume.prop.embedding.section3.1.1}\\
\langle \mathcal{A}' f, g\rangle _{\dot{H}^{-1},\dot{H}^1} & = \langle \mathcal{P}_{\mathcal{A}} f, \Lambda_{\mathcal{A}^*} g + M_{\bf \bar{r}_1}\cdot \nabla_x g\rangle _{L^2(\R^d)}\label{assume.prop.embedding.section3.1.2}
\end{align}
for $f\in D_{L^2}(\mathcal{P}_{\mathcal{A}})$ and $g\in H^1(\R^d)$.

\end{prop}

\noindent {\it Proof.} Assume that the statement (ii) holds. Then, from Proposition \ref{prop.subsec.relation.2}  we have $\|\Lambda_{\mathcal{A}} f\|_{L^2(\R^d)} \leq C (\| \mathcal{P}_{\mathcal{A}}f \|_{L^2(\R^d)} + \|f \|_{L^2(\R^d)} )$ for $f\in D_{H^{1/2}}(\mathcal{P}_{\mathcal{A}})$. Thus, since $D_{H^{1/2}}(\mathcal{P}_{\mathcal{A}})$ is dense in $D_{L^2}(\mathcal{P}_{\mathcal{A}})$ and $\Lambda_{\mathcal{A}}$ is closed in $L^2 (\R^d)$, the space $D_{L^2}(\mathcal{P}_{\mathcal{A}})$ is continuously embedded in $H^1(\R^d)\cap D_{L^2}(\Lambda_{\mathcal{A}})$.  Next we take $f=(1+\mathcal{P}_{\mathcal{A}})^{-1}h$ with $h\in H^{1/2}(\R^d)$ in \eqref{eq.prop.thm.factorization}. Then we have from $(\Lambda_{\mathcal{A}^*})^* =\Lambda_{\mathcal{A}}$,
\begin{align}
\langle  h, \Lambda_{\mathcal{A}^*} g\rangle _{\dot{H}^{\frac12}, \dot{H}^{-\frac12}} & = \langle \Lambda_{\mathcal{A}} (1+\mathcal{P}_{\mathcal{A}})^{-1}h, g\rangle _{L^2(\R^d)}+ \langle A'\nabla_x (1+\mathcal{P}_{\mathcal{A}})^{-1}h, \nabla_x g\rangle _{L^2(\R^d)} \nonumber \\
& ~~~ - \langle \mathcal{P}_{\mathcal{A}} (1 + \mathcal{P}_{\mathcal{A}})^{-1} h, M_{{\bf \bar{r}_1}}\cdot \nabla_x  g\rangle _{L^2(\R^d)}\label{proof.prop.embedding.section3.1.1}
\end{align}
for any $h\in H^{1/2}(\R^d)$ and $g\in H^1(\R^d)$. Since the right-hand side of \eqref{proof.prop.embedding.section3.1.1} can be extended to $h\in L^2(\R^d)$, we conclude that $H^1(\R^d)$ is continuously embedded in $D_{L^2}(\Lambda_{\mathcal{A}^*})$. The identities \eqref{assume.prop.embedding.section3.1.1} and \eqref{assume.prop.embedding.section3.1.2} follow from Proposition \ref{prop.thm.factorization} and the above embeddings  by using the fact that $D_{H^{1/2}}(\mathcal{P}_{\mathcal{A}})$ is dense in $D_{L^2}(\mathcal{P}_{\mathcal{A}})$.

Next we assume that the statement (i) holds. We first show that for any $f\in L^2(\R_+; D_{H^{1/2}}(\mathcal{P}_{\mathcal{A}}))$ the function $\Psi_{\mathcal{P}_{\mathcal{A}}}[f] (t) = \int_0^t e^{-(t-s)\mathcal{P}_{\mathcal{A}}} f(s)\dd s$ satisfies 
\begin{align}
\| \mathcal{P}_{\mathcal{A}} \Psi_{\mathcal{P}_{\mathcal{A}}}[f] \|_{L^2(0,T; L^2 (\R^d ) )}\leq C \| f\|_{L^2 (0,T; L^2 (\R^d ) )},~~~~~~~~T\in (0,1].\label{proof.prop.embedding.section3.1.2}
\end{align}
Fix $T\in (0,1]$. We may assume that $f(t)=0$ for $t> T$. Let $\phi_M^c$, $M>2$, be the cut-off function in the proof of Proposition \ref{prop.poisson.analytic}. Then by extending $\Psi_{\mathcal{P}_{\mathcal{A}}}[f] (t)$ by zero for $t<0$ it is easy to see that $\phi_M^c \Psi_{\mathcal{P}_{\mathcal{A}}}[f]\in L^2 (\R; D_{H^{1/2}}(\mathcal{P}_{\mathcal{A}}))\cap H^1 (\R; H^{1/2}(\R^d))\hookrightarrow H^1 (\R^{d+1})$. Thus from \eqref{eq.factorization'} we have 
\begin{align*}
& ~~~~ \langle ~A\nabla (\phi_M^c \Psi_{\mathcal{P}_{\mathcal{A}}}[f]),~ \nabla (\phi_M^c \Psi_{\mathcal{P}_{\mathcal{A}}}[f])~\rangle _{L^2 (\R^{d+1})} \\
& = \int_0^\infty \langle (\partial_t + \mathcal{P}_{\mathcal{A}}) \phi_M^c \Psi_{\mathcal{P}_{\mathcal{A}}}[f] (t) , (M_{\bar{b}} \partial_t + M_{{\bf \bar{r}_1}}\cdot \nabla_x ) \phi_M^c \Psi_{\mathcal{P}_{\mathcal{A}}}[f] (t) \rangle _{L^2 (\R^d)} \dd t\\
&~~~ + \int_0^\infty  \langle (\partial_t + \mathcal{P}_{\mathcal{A}}) \phi_M^c \Psi_{\mathcal{P}_{\mathcal{A}}}[f] (t) ,\Lambda_{\mathcal{A}^*} \phi_M^c \Psi_{\mathcal{P}_{\mathcal{A}}}[f] (t) \rangle _{\dot{H}^\frac12, \dot{H}^{-\frac12}} \dd t.
\end{align*}
Since $(\partial_t + \mathcal{P}_{\mathcal{A}}) \Psi_{\mathcal{P}_{\mathcal{A}}}[f] = f$ the elliptic condition on $A$ yields 
\begin{align*}
\| \nabla ( \phi_M^c \Psi_{\mathcal{P}_{\mathcal{A}}}[f] )\|^2_{L^2 (\R^{d+1})} &\leq  C \int_0^T \| f(t) \|_{L^2 (\R^d)}^2 \dd t + C M^{-2} \int_{\frac{M}{2}}^M \| \Psi_{\mathcal{P}_{\mathcal{A}}}[f] (t)\| _{L^2 (\R^d)}^2 \dd t\\
&~~~ + C {\rm Re} \int_0^T   \langle f, \Lambda_{\mathcal{A}^*} \Psi_{\mathcal{P}_{\mathcal{A}}}[f] (t) \rangle _{\dot{H}^\frac12, \dot{H}^{-\frac12}} \dd t \\
& ~~~~~~+ C {\rm Re}   \int_0^\infty (\phi_M^c)' \phi_M^c \langle \Psi_{\mathcal{P}_{\mathcal{A}}}[f] (t),  \Lambda_{\mathcal{A}^*} \Psi_{\mathcal{P}_{\mathcal{A}}}[f] (t) \rangle _{\dot{H}^\frac12, \dot{H}^{-\frac12}} \dd t\\
& \leq  C \int_0^T \| f(t) \|_{L^2 (\R^d)}^2 \dd t + C M^{-2} \int_{\frac{M}{2}}^M \| \Psi_{\mathcal{P}_{\mathcal{A}}}[f] (t)\| _{L^2 (\R^d)}^2 \dd t\\
&~~~ + C {\rm Re} \int_0^T    \langle f, \Lambda_{\mathcal{A}^*} \Psi_{\mathcal{P}_{\mathcal{A}}}[f] (t) \rangle _{\dot{H}^\frac12, \dot{H}^{-\frac12}} \dd t. 
\end{align*}
Here we have used the fact $(\phi_M^c)'\leq 0$. Since $\Psi_{\mathcal{P}_{\mathcal{A}}}[f] (t)= e^{-(t-T)\mathcal{P}_{\mathcal{A}}} \Psi_{\mathcal{P}_{\mathcal{A}}}[f] (T)$ for $t>T$,  by using Proposition \ref{prop.poisson.analytic} we take $M\rightarrow \infty$ in the above inequality and arrive at 
\begin{align}
\| \nabla  \Psi_{\mathcal{P}_{\mathcal{A}}}[f] \|^2_{L^2 (\R^{d+1}_+)} &\leq  C \int_0^T \| f(t) \|_{L^2 (\R^d)}^2 \dd t  + C {\rm Re} \int_0^T    \langle f (t) , \Lambda_{\mathcal{A}^*} \Psi_{\mathcal{P}_{\mathcal{A}}}[f] (t) \rangle _{\dot{H}^\frac12, \dot{H}^{-\frac12}} \dd t. \label{proof.prop.embedding.section3.1.3'}
\end{align}
Now  the assumption of (i) leads to 
\[
{\rm Re} \langle f (t) , \Lambda_{\mathcal{A}^*} \Psi_{\mathcal{P}_{\mathcal{A}}}[f] (t) \rangle _{\dot{H}^\frac12, \dot{H}^{-\frac12}} \leq C \| f (t) \|_{L^2 (\R^d)} \big ( \| \nabla \Psi_{\mathcal{P}_{\mathcal{A}}}[f] (t)\|_{L^2 (\R^d)} + \| \Psi_{\mathcal{P}_{\mathcal{A}}}[f] (t)\|_{L^2 (\R^d)}  \big ).
\]
Then, applying the estimate $\|\Psi_{\mathcal{P}_{\mathcal{A}}}[f] (t)\|_{L^2 (\R^d)} \leq Ct^{1/2} \| \partial_t \Psi_{\mathcal{P}_{\mathcal{A}}}[f]\|_{L^2 (0,t; L^2 (\R^d))}$, the left-hand side of \eqref{proof.prop.embedding.section3.1.3'} is bounded from above by $C \| f\|_{L^2 (0,T; L^2(\R^{d}) )}^2$. Thus  we get \eqref{proof.prop.embedding.section3.1.2} from $\mathcal{P}_{\mathcal{A}} \Psi_{\mathcal{P}_{\mathcal{A}}} [f] =-\partial_t\Psi_{\mathcal{P}_{\mathcal{A}}} [f] + f$, as desired. Set $u(t) = e^{-t\mathcal{P}_{\mathcal{A}}} g$ with $g\in D_{H^{1/2}}(\mathcal{P}_{\mathcal{A}})$. Then the identity $u(t) = g - \mathcal{P}_{\mathcal{A}}\Psi _{\mathcal{P}_{\mathcal{A}}} [g](t)$ and \eqref{proof.prop.embedding.section3.1.2} imply $\| u\|_{L^2 (0,T; L^2 (\R^d))} \leq C T^{1/2} \| g\|_{L^2(\R^d)}$ for any $T\in (0,1]$.  Thus, for any $t\in (0,1]$ there is $t_0\in [t/2,t]$ such that $\|u(t_0) \|_{L^2(\R^d)} \leq C \| g \|_{L^2(\R^d)}$ with $C$ independent of $f$. Then we have from \eqref{est.proof.prop.poisson.analytic.2} and  \eqref{est.rem.prop.poisson.analytic.1},
\begin{align}
\| u(t) \|_{L^2 (\R^d)}  \leq \| u (t_0) \|_{L^2 (\R^d)}+ C t^\frac12 \| u(t_0 )\|_{\dot{H}^\frac12(\R^d)} & \leq C \| g \|_{L^2(\R^d)} + C t^\frac12 t_0^{-1} \| u \|_{L^2(\frac{t_0}{2}, t_0; L^2 (\R^d))} \nonumber \\
& \leq C\| g \|_{L^2(\R^d)} + C t^\frac12 t_0^{-\frac12} \| g \|_{L^2 (\R^d)}\leq C \| g \|_{L^2 (\R^d)}.\label{proof.prop.embedding.section3.1.3}
\end{align} 
Here we have also used $t_0\in [t/2,t]$. By  Remark \ref{rem.prop.poisson.analytic.L^2} the estimate \eqref{proof.prop.embedding.section3.1.3} implies that $\{ e^{-\mathcal{P}_{\mathcal{A}}}\}_{t\geq 0}$ is extended as a strongly continuous semigroup in $L^2(\R^d)$. Finally, taking $g=f$ in \eqref{eq.prop.thm.factorization} and using (i), we have $\| \nabla_x f \|_{L^2(\R^d)}\leq C (\|\mathcal{P}_{\mathcal{A}} f \|_{L^2(\R^d)} + \| f \|_{L^2(\R^d)})$ for any $f\in D_{H^{1/2}}(\mathcal{P}_{\mathcal{A}})$. Since $D_{H^{1/2}}(\mathcal{P}_{\mathcal{A}})$ is dense in $D_{L^2}(\mathcal{P}_{\mathcal{A}})$, this embedding inequality is valid for all $f\in D_{L^2}(\mathcal{P}_{\mathcal{A}})$. Hence (ii) follows.  The proof of Proposition \ref{prop.embedding.section3.1} is now complete.

\begin{cor}\label{cor.prop.embedding.section3.1'} Assume that $\{e^{-t\mathcal{P}_{\mathcal{A}^*}}\}_{t\geq 0}$ is extended as a strongly continuous semigroup in $L^2 (\R^d)$. Assume further that  $H^1 (\R^d)$ is continuously embedded in $D_{L^2}(\mathcal{P}_{\mathcal{A}^*})$ and 
$\displaystyle \liminf_{t\rightarrow 0} \| \nabla_x e^{-t\mathcal{P}_{\mathcal{A}^*}}  f  \|_{L^2 (\R^d)} <\infty$ holds for any 
$f \in H^1 (\R^d)$. Then the condition {\rm (i)} in Proposition \ref{prop.embedding.section3.1} holds and 
\begin{align}
\mathcal{P}_{\mathcal{A}^*} f = M_{1/\bar{b}} \Lambda_{\mathcal{A}^*} f + M_{{\bf \bar{r}_1}/\bar{b}} \cdot \nabla_x f,~~~~~~\quad f\in H^1 (\R^d). \label{eq.cor.prop.embedding.section3.1'} 
\end{align}
\end{cor}

\noindent {\it Proof.}  Under the assumptions of the corollary we have  $\mathcal{P}_{\mathcal{A}^*} e^{-t\mathcal{P}_{\mathcal{A}^*}} f \rightarrow \mathcal{P}_{\mathcal{A}^*} f$ in $L^2 (\R^d)$ for all $f\in H^1 (\R^d)\hookrightarrow D_{L^2} (\mathcal{P}_{\mathcal{A}^*})$ and we can also take a sequence $\{t_n\}$, $t_n\rightarrow 0$, such that $\nabla_x e^{-t_n\mathcal{P}_{\mathcal{A}^*}} f$  converges to $\nabla_x f$ weakly in $L^2 (\R^d)$. Then \eqref{eq.prop.subsec.relation.2}  (with $\mathcal{A}$ replaced by $\mathcal{A}^*$) implies that 
\begin{align*}
\langle \Lambda_{\mathcal{A}} g, f\rangle_{L^2 (\R^d)} = \lim_{n\rightarrow \infty} \langle \Lambda_{\mathcal{A}} g, e^{-t_n\mathcal{P}_{\mathcal{A}^*}} f\rangle_{L^2 (\R^d)} &  = \lim_{n\rightarrow \infty} \langle g,  \big ( M_{\bar{b}} \mathcal{P}_{\mathcal{A}^*}  - M_{\bf \bar{r}_1} \cdot \nabla_x \big )e^{-t_n\mathcal{P}_{\mathcal{A}^*}} f\rangle _{L^2 (\R^d)}\\
& = \langle g,   M_{\bar{b}} \mathcal{P}_{\mathcal{A}^*} f  - M_{\bf \bar{r}_1} \cdot \nabla_x f\rangle _{L^2 (\R^d)}.
\end{align*}
for all $g\in D_{L^2}(\Lambda_{\mathcal{A}})$. Hence we have  $f \in D_{L^2}(\Lambda_{\mathcal{A}^*})$ and \eqref{eq.cor.prop.embedding.section3.1'}. Since $D_{H^{1/2}}(\mathcal{P}_{\mathcal{A}})\hookrightarrow H^1 (\R^d)$  the condition {\rm (i)} in Proposition \ref{prop.embedding.section3.1} is also valid.  The proof is complete.

\begin{cor}\label{cor.prop.embedding.section3.1}  Assume that $\{e^{-t\mathcal{P}_{\mathcal{A}}}\}_{t\geq 0}$ and $\{e^{-t\mathcal{P}_{\mathcal{A}^*}}\}_{t\geq 0}$ are extended as  strongly continuous semigroups in $L^2 (\R^d)$ and that $D_{L^2}(\mathcal{P}_{\mathcal{A}})$ and $D_{L^2}(\mathcal{P}_{\mathcal{A}^*})$ are continuously embedded in $H^1 (\R^d)$. Then we have 
\begin{align}
C' \| f \|_{H^1 (\R^d)} \leq \| \mathcal{P}_{\mathcal{A}} f\|_{L^2 (\R^d)} + \| f\|_{L^2(\R^d)} \leq C  \| f\|_{H^1 (\R^d)},~~~~~~~~~f\in D_{L^2}(\mathcal{P}_{\mathcal{A}}).
\end{align}
Moreover, $\langle A'\nabla_x f, \nabla_x g\rangle _{L^2 (\R^d)} = \langle \mathcal{P}_{\mathcal{A}}f, M_{\bar{b}} \mathcal{P}_{\mathcal{A}^*} g\rangle _{L^2 (\R^d)}$ holds for all $f\in D_{L^2}(\mathcal{P}_{\mathcal{A}})$ and $g\in D_{L^2}(\mathcal{P}_{\mathcal{A}^*})$.  

\end{cor}

\noindent {\it Proof.} By Proposition \ref{prop.embedding.section3.1} the space $H^1(\R^d)$ is continuously embedded in  $D_{L^2}(\Lambda_{\mathcal{A}})$, and \eqref{assume.prop.embedding.section3.1.1} holds  for $f\in D_{L^2}(\mathcal{P}_{\mathcal{A}})$. 
Hence we have $\| \mathcal{P}_{\mathcal{A}} f \|_{L^2 (\R^d)} \leq C (\| \Lambda_{\mathcal{A}} f \|_{L^2(\R^d)} + \| \nabla_x f \|_{L^2 (\R^d )} ) \leq C \| f \|_{H^1 (\R^d)}$ for all $f\in D_{L^2}(\mathcal{P}_{\mathcal{A}})$. The last assertion on the identity is obtained by applying  Proposition \ref{prop.embedding.section3.1} to $\mathcal{P}_{\mathcal{A}}$ and $\mathcal{P}_{\mathcal{A}^*}$. The proof is complete.

\vspace{0.5cm}

Although Corollary  \ref{cor.prop.embedding.section3.1} gives the comparability of $\|f\|_{H^1(\R^d)}$ with the graph norm of $D_{L^2}(\mathcal{P}_{\mathcal{A}})$ in the case of  $f\in D_{L^2}(\mathcal{P}_{\mathcal{A}})$, this does not mean $H^1(\R^d)=D_{L^2}(\mathcal{P}_{\mathcal{A}})$ because we do not know whether $D_{L^2}(\mathcal{P}_{\mathcal{A}})$ is dense in $H^1 (\R^d)$ or not in general. In order to obtain the exact characterization of $D_{L^2}(\mathcal{P}_{\mathcal{A}})$ we need an additional estimate as follows.
\begin{prop}\label{prop.embedding.section3.2} {\rm (i)} Assume that  $\{ e^{-t\mathcal{P}_{\mathcal{A}}}\}_{t\geq 0}$ is  extended as  strongly continuous semigroups in $L^2 (\R^d)$ and that  
\begin{align}
\liminf_{t\rightarrow 0} \| \frac{\dd}{\dd t} e^{-t\mathcal{P}_{\mathcal{A}}} f \|_{L^2(\R^d )}<\infty ~~~~~~{\rm for~all}~~f\in H^1 (\R^d).\label{assume.prop.embedding.section3.2}
\end{align}
Then $H^1 (\R^d) \subset D_{L^2}(\mathcal{P}_{\mathcal{A}})$.

\noindent {\rm (ii)} Suppose that the assumptions of Corollary \ref{cor.prop.embedding.section3.1} holds. Assume further that 
\begin{align}
\liminf_{t\rightarrow 0} \| \frac{\dd}{\dd t} e^{-t\mathcal{P}_{\mathcal{A}}} f \|_{L^2(\R^d )}<\infty ~~~~~~{\rm for~all}~~f\in C_0^\infty (\R^d).\label{assume.prop.embedding.section3.3}
\end{align}
Then $D_{L^2}(\mathcal{P}_{\mathcal{A}}) = H^1 (\R^d)$ with equivalent norms.
\end{prop}

\noindent {\it Proof.} (i) We note that $e^{-t \mathcal{P}_{\mathcal{A}}}f \rightarrow f$ in $L^2 (\R^d)$ as $t\rightarrow 0$ for any $f\in H^1 (\R^d)$. Then for $f\in H^1 (\R^d)$ and $g\in D_{L^2}((\mathcal{P}_{\mathcal{A}})^*)$ we have 
\begin{align*}
|\langle f, (\mathcal{P}_{\mathcal{A}})^*g\rangle _{L^2(\R^d)}| = \lim_{t \rightarrow 0}| \langle e^{-t \mathcal{P}_{\mathcal{A}}}f, (\mathcal{P}_{\mathcal{A}})^*g\rangle _{L^2(\R^d)}|= \liminf_{t\rightarrow 0}| \langle \mathcal{P}_{\mathcal{A}}e^{-t \mathcal{P}_{\mathcal{A}}}f, g\rangle _{L^2(\R^d)}|\leq C_f \| g\|_{L^2 (\R^d)},
\end{align*}
where $C_f = \liminf_{t\rightarrow 0} \| \dd /\dd t ~ e^{-t\mathcal{P}_{\mathcal{A}}} g \|_{L^2(\R^d )}<\infty$ by the assumption. Hence the sesquilinear functional $S_f: D_{L^2}((\mathcal{P}_{\mathcal{A}})^*) \ni g \mapsto S_f (g)= \langle f, (\mathcal{P}_{\mathcal{A}})^*g\rangle _{L^2(\R^d)}\in \C$ is extended as a bounded sesquilinear functional acting on $L^2 (\R^d)$. Then the representation theorem implies the existence of $h\in L^2 (\R^d)$ such that $S_f (g) = \langle h, g\rangle _{L^2(\R^d)}$ for all $g\in L^2 (\R^d)$. In particular, $ \langle f, (\mathcal{P}_{\mathcal{A}})^*h\rangle _{L^2(\R^d)} = \langle h, g\rangle _{L^2(\R^d)}$ holds for all $g\in D_{L^2}((\mathcal{P}_{\mathcal{A}})^*)$, which gives $f\in D_{L^2}((\mathcal{P}_{\mathcal{A}})^{**})$. Since $\mathcal{P}_{\mathcal{A}}$ is a densely defined closed operator and $L^2(\R^d)$ is reflexive, we have $(\mathcal{P}_{\mathcal{A}})^{**}=\mathcal{P}_{\mathcal{A}}$, i.e., $f\in D_{L^2}(\mathcal{P}_{\mathcal{A}})$. The proof of (i) is complete.

\noindent (ii) It suffices to show $H^1 (\R^d)\subset D_{L^2}(\mathcal{P}_{\mathcal{A}})$. The argument of (i) implies that if \eqref{assume.prop.embedding.section3.3} holds then each $f\in C_0^\infty (\R^d)$ belongs to $D_{L^2}(\mathcal{P}_{\mathcal{A}})$.  For any $f\in H^1 (\R^d)$ let us take $\{f_n\}\subset C_0^\infty (\R^d)$ such that $f_n\rightarrow f$ in $H^1 (\R^d)$. Then for each $t\in (0,1)$ there is $n(t)\in \N$ such that $\| f_{n(t)} -f \|_{L^2(\R^d)}\leq t$. Thus we have for $t\in (0,1)$,
\begin{align*}
\| \frac{\dd}{\dd t}e^{-t\mathcal{P}_{\mathcal{A}}} f\|_{L^2(\R^d)} & \leq \|\frac{\dd}{\dd t} e^{-t\mathcal{P}_{\mathcal{A}}} f_{n(t)} \|_{L^2 (\R^d)} + \| \frac{\dd}{\dd t} e^{-t\mathcal{P}_{\mathcal{A}}} (f_{n(t)} - f) \|_{L^2 (\R^d)} \\
& \leq C \| \mathcal{P}_{\mathcal{A}} f_{n(t)} \|_{L^2 (\R^d)}  + C t^{-1} \| f_{n(t)} - f \|_{L^2 (\R^d)} \\
& \leq C \| f_{n(t)} \|_{H^1 (\R^d)}  + C  \leq C ( 1+ \sup_n \| f_n \|_{H^1 (\R^d)} )<\infty.
\end{align*}
Here we have used Corollary \ref{cor.prop.embedding.section3.1}. Hence \eqref{assume.prop.embedding.section3.2} holds, i.e., we have $H^1 (\R^d)\subset D_{L^2}(\mathcal{P}_{\mathcal{A}})$ by the assertion (i). The proof is complete.

\subsection{Characterization of $D_{L^2}(\mathcal{P}_{\mathcal{A}})$ and factorization of $\mathcal{A}$ in $L^2(\R^d)$}\label{subsec.lem.factorization}

Next lemma shows that the factorization of $\mathcal{A}$ as in Theorem \ref{thm.factorization.strong.1} is obtained as a consequence of the $H^1(\R^d)$ characterization of $D_{L^2}(\mathcal{P}_{\mathcal{A}})$ and $D_{L^2}(\mathcal{P}_{\mathcal{A}^*})$.

\begin{lem}\label{lem.factorization.strong} Assume that the semigroups $\{ e^{-t\mathcal{P}_{\mathcal{A}}}\}_{t\geq 0}$  and $\{ e^{-t\mathcal{P}_{\mathcal{A}^*}}\}_{t\geq 0}$ in $H^{1/2}(\R^d)$ are extended as  strongly continuous semigroups in $L^2 (\R^d)$ and that $D_{L^2}(\mathcal{P}_{\mathcal{A}})= D_{L^2}(\mathcal{P}_{\mathcal{A}^*}) = H^1 (\R^d)$ holds with equivalent norms. Then $H^1 (\R^d)$ is continuously embedded in $D_{L^2}(\Lambda_{\mathcal{A}})\cap D_{L^2}(\Lambda_{\mathcal{A}^*})$ and 
\begin{align}
\mathcal{P}_{\mathcal{A}} f & = M_{1/b} \Lambda_{\mathcal{A}} f + M_{{\bf r_2}/b}\cdot \nabla_x f,~~~~~~f\in H^1 (\R^d),\label{eq.lem.factorization.strong.1} \\
\mathcal{P}_{\mathcal{A}^*} g &  = M_{1/\bar{b}} \Lambda_{\mathcal{A}^*} g + M_{{\bf \bar{r}_1}/\bar{b}}\cdot \nabla_x g,~~~~~~ g\in H^1 (\R^d).\label{eq.lem.factorization.strong.2}
\end{align}
Moreover, the realization of $\mathcal{A}'$ in $L^2 (\R^d)$ and the realization of $\mathcal{A}$ in $L^2 (\R^{d+1})$ are respectively factorized as
\begin{align}
\mathcal{A}'& = M_b \mathcal{Q}_{\mathcal{A}} \mathcal{P}_{\mathcal{A}},~~~~~~~~~~~~~~\quad \mathcal{Q}_{\mathcal{A}} = M_{1/b} ( M_{\bar{b}} \mathcal{P}_{\mathcal{A}^*} )^*,\label{eq.lem.factorization.strong.4}\\
\mathcal{A} & = - M_b (\partial _t - \mathcal{Q}_{\mathcal{A}} ) ( \partial_t + \mathcal{P}_{\mathcal{A}}).\label{eq.lem.factorization.strong.3}
\end{align}
Here $( M_{\bar{b}} \mathcal{P}_{\mathcal{A}^*} )^*$ is the adjoint of $ M_{\bar{b}} \mathcal{P}_{\mathcal{A}^*}$ in $L^2 (\R^d)$. 

\end{lem}

\noindent {\it Proof.}  Proposition \ref{prop.embedding.section3.1} shows that $H^1 (\R^d)\hookrightarrow D_{L^2}(\Lambda_{\mathcal{A}})\cap D_{L^2}(\Lambda_{\mathcal{A}^*})$ and \eqref{eq.lem.factorization.strong.1}-\eqref{eq.lem.factorization.strong.2} hold. Moreover, we have $\langle A' \nabla_x f, \nabla_x g \rangle _{L^2 (\R^d)} = \langle \mathcal{P}_{\mathcal{A}} f, M_{\bar{b}} \mathcal{P}_{\mathcal{A}^*} g\rangle _{L^2 (\R^d)}$ for all $f,g \in H^1 (\R^d)$. If $f\in D_{L^2}(\mathcal{A}')\subset H^1 (\R^d)$ then this identity implies $\mathcal{P}_{\mathcal{A}} f\in D_{L^2}\big ( (M_{\bar{b}} \mathcal{P}_{\mathcal{A}^*})^* \big )$ and thus $\mathcal{A}'\subset M_b \mathcal{Q}_{\mathcal{A}}\mathcal{P}_{\mathcal{A}}$ holds. The converse relation is proved similarly, and \eqref{eq.lem.factorization.strong.4} is proved. Next we consider \eqref{eq.lem.factorization.strong.3}. Theorem \ref{thm.factorization} and Remark \ref{rem.thm.factorization} verify the identity 
\begin{align}
\langle A\nabla u, \nabla v\rangle_{L^2 (\R^{d+1})} = \langle (\partial_t + \mathcal{P}_{\mathcal{A}} ) u, M_{\bar{b}} (\partial_t + \mathcal{P}_{\mathcal{A}^*} ) v \rangle _{L^2 (\R^{d+1})},~~~~~~~~~u,~v \in H^1 (\R^{d+1}).\label{proof.lem.factorization.strong.1}
\end{align}
Set 
\begin{align*}
& S = \partial_t + \mathcal{P}_{\mathcal{A}},~~~~~~D_{L^2} (S) = \{ u\in L^2 (\R^{d+1})~|~ \partial_t u, \mathcal{P}_{\mathcal{A}} u\in L^2 (\R^{d+1})\} = H^1 (\R^{d+1}),\\
& T = \partial_t - \mathcal{Q}_{\mathcal{A}},~~~~~~D_{L^2} (T) = \{ u\in L^2 (\R^{d+1})~|~ \partial_t u, (M_{\bar{b}} \mathcal{P}_{\mathcal{A}^*})^* u\in L^2 (\R^{d+1})\}.
\end{align*} 
Assume that $u\in D_{L^2}(\mathcal{A})$, which is continuously embedded in $H^1 (\R^{d+1})= D_{L^2}(S)$. We observe the inequality $\| \partial_t \nabla u \|_{L^2 (\R^{d+1})} \leq C (\| \mathcal{A} u \|_{L^2 (\R^{d+1})} + \| u\|_{L^2 (\R^{d+1})} )$, which then  implies from \eqref{proof.lem.factorization.strong.1} that 
\begin{align*}
|\langle S u, M_{\bar{b}} \mathcal{P}_{\mathcal{A}^*}  v \rangle _{L^2 (\R^{d+1})} | &\leq | \langle  M_{b} S \partial_t u, v \rangle _{L^2 (\R^{d+1})}|  + | \langle \mathcal{A} u, v\rangle _{L^2 (\R^{d+1})}|\\
& \leq  C ( \| \partial_t u \|_{H^1 (\R^{d+1})} + \|\mathcal{A} u \|_{L^2 (\R^{d+1})} ) \| v \|_{L^2 (\R^{d+1})} \\
& \leq  C ( \| \mathcal{A} u \|_{L^2 (\R^{d+1})} + \| u \|_{L^2 (\R^{d+1})}) \| v \|_{L^2 (\R^{d+1})}.
\end{align*}
Hence $S u \in D_{L^2}( (M_{\bar{b}} \mathcal{P}_{\mathcal{A}^*})^* )= D_{L^2}( \mathcal{Q}_{\mathcal{A}} )$. Thus we have $S u \in D_{L^2}(T) = D_{L^2}(M_b T)$, that is, $u \in D_{L^2}(M_b TS)$. Conversely, assume that $u\in D_{L^2}(M_b TS)\subset H^1 (\R^{d+1})$. Then $S u\in D_{L^2}(T)$ and we have from \eqref{proof.lem.factorization.strong.1},
\begin{align*}
 \langle A\nabla u, \nabla v\rangle_{L^2 (\R^{d+1})} & = \langle S u, M_{\bar{b}} (\partial_t + \mathcal{P}_{\mathcal{A}^*} ) v \rangle _{L^2 (\R^{d+1})} \\
& = - \langle M_b \partial_t S u, v\rangle_{L^2 (\R^{d+1})} + \langle  (M_{\bar{b}} \mathcal{P}_{\mathcal{A}^*})^* S u, v \rangle _{L^2 (\R^{d+1})} = \langle - M_b TS u, v\rangle _{L^2 (\R^{d+1})}. 
\end{align*}
This shows $u\in D_{L^2}(\mathcal{A})$ and $\mathcal{A}u = - M_b TS u$. The proof is complete.

\section{Analysis of Poisson operator for specific cases}\label{sec.poisson}

For general elliptic operators the verification of $D_{L^2}(\mathcal{P}_{\mathcal{A}})=H^1 (\R^d)$ is a difficult problem in general. As far as the authors know,  at least  the following cases are settled so far; (I) $A$ is a constant matrix, (II) $A$ is an Hermitian, i.e., $A=A^*$, (III) $A$ is block-type, i.e., ${\bf r_1}={\bf r_2}={\bf 0}$, (IV) $A$ is a small perturbation of $B$ in $(L^\infty(\R^d))^{(d+1)\times (d+1)}$, where $B$ satisfies one of (I)-(III) above.  Roughly speaking, the relation $D_{L^2}(\mathcal{P}_{\mathcal{A}})=H^1 (\R^d)$ follows from $L^2$ solvability of both Dirichlet and Neumann problem. The case (I) is verified by the direct application of the Fourier transform. As stated in the introduction, the case (II) is classical and $L^2$ solvability is proved by \cite{Dahlberg1,JK1,JK2,Verchota,Dahlberg2,KP,AAM}. The case (III) is related with the Kato square root problem, which was finally solved by \cite{AHLMT}. The stability with respect to small $L^\infty$ perturbations,  the case (IV), is solved by \cite{FJK} when $B$ is a constant matrix, by \cite{H,AAAHK} when $B$ is a constant or  real symmetric matrix, by \cite{AAH} when $B$ is a constant, real symmetric, or block matrix, and by \cite{AAM} when $B$ is a constant, Hermite, or block matrix.  In this section we study the domain of the generator for the Poisson semigroup when $A$ is Hermite (Section \ref{subsec.poisson.hermite}) and when the off-block vectors ${\bf r_1}$, ${\bf r_2}$ have additional regularity (Section \ref{subsec.poisson.regular}).



\subsection{Analysis for Hermitian case}\label{subsec.poisson.hermite}

In this section we consider the case $A$ is Hermite.  This is a classical case, and the derivation of various estimates for Poisson semigroups here is based on the Rellich type identity \cite{Rellich} as already observed in \cite{PW,JK1,JK2}. Strictly speaking, in \cite{PW,JK1,JK2} the problem is studied when the domain is bounded and  $A$ is real symmetric, and one has to be careful about the noncompactness of the domain $\R^{d+1}_+$ in the present problem. In the proof of Theorem \ref{thm.poisson.hermite} below, to verify $D_{L^2}(\mathcal{P}_{\mathcal{A}})=D_{L^2}(\Lambda_{\mathcal{A}})=H^1 (\R^d)$  based on the Rellich identity we firstly approximate the nonsmooth $A$ by smooth ones $\{A_\epsilon\}_{\epsilon>0}$ and then take  the limit. In this approach  the results for the smooth $A_\epsilon$ have to be obtained in advance.  If  the domain is bounded then  the case of smooth $A$ is handled  by using the localization and perturbation argument from the  constant matrix case. This method is quite robust (also for $L^p$ framework),  but when the domain is unbounded  it does not always work successfully. Thus, although Theorem \ref{thm.poisson.hermite} itself is not essentially new, we will give an alternative proof using the result in our companion work \cite{MM1}, where $D_{L^2}(\mathcal{P}_{\mathcal{A}})=D_{L^2}(\Lambda_{\mathcal{A}})=H^1 (\R^d)$ is verified when $A$ is Lipschitz continuous. 

\begin{lem}\label{lem.poisson.hermite} Assume that $A$ is Hermite, i.e., $A^*=A$. Let $f,g \in D_{H^{1/2}} (\mathcal{P}_\mathcal{A})$ and set $\Psi_{\mathcal{P}_\mathcal{A}} [F] (t) = \int_0^t e^{-(t-s)\mathcal{P}_\mathcal{A}} F  (s) \dd s$ for $F \in C_0^\infty (\R^{d+1}_+)$. Then 
\begin{align}
\langle A'\nabla_x f, \nabla_x f \rangle_{L^2(\R^d)} & =\| M_{\sqrt{b}} \mathcal{P}_\mathcal{A} f \|_{L^2(\R^d)}^2,\label{est.prop.poisson.hermite.1} \\
C_4^{-1} \|\nabla_x f \|_{L^2 (\R^d)} \leq \| \Lambda_{\mathcal{A}} f \|_{L^2(\R^d)} & \leq C_4 \| \nabla_x f \|_{L^2(\R^d)}, \label{est.prop.poisson.hermite.2} \\ 
\sup_{t>0} \| e^{-t\mathcal{P}_\mathcal{A}} g \|_{L^2 (\R^d)} + \sup_{t>0} t \|  \mathcal{P}_\mathcal{A} e^{-t\mathcal{P}_\mathcal{A}} g \|_{L^2 (\R^d)} & \leq C_5 \| g \|_{L^2 (\R^d)},\label{est.prop.poisson.hermite.3}\\
\| \partial_t \Psi_{\mathcal{P}_\mathcal{A}} [F] \|_{L^2 (\R^{d+1}_+)} + \| \mathcal{P}_\mathcal{A}  \Psi_{\mathcal{P}_\mathcal{A}} [F]  \|_{L^2 (\R^{d+1}_+)} & \leq C_5 \| F \|_{L^2 (\R^{d+1}_+)}. \label{est.prop.poisson.hermite.4}
\end{align}
Here the constants $C_4,~C_5$ depend only on $\nu_1$ and $\nu_2$.  

\end{lem}

\noindent {\it Proof.} The equality \eqref{est.prop.poisson.hermite.1} directly follows from \eqref{eq.prop.thm.factorization}. The inequality $\| \Lambda_{\mathcal{A}} f \|_{L^2 (\R^d)} \leq C_4 \| \nabla_x f \|_{L^2 (\R^d)}$ follows from \eqref{eq.prop.subsec.relation.2} and \eqref{est.prop.poisson.hermite.1}, while $\| \Lambda_{\mathcal{A}} f \|_{L^2 (\R^d)} \geq C_4^{-1} \| \nabla_x f \|_{L^2 (\R^d)}$ follows by combining \eqref{eq.prop.subsec.relation.2}, \eqref{est.prop.poisson.hermite.1}, and $\langle A'\nabla_x f, \nabla_x f \rangle_{L^2(\R^d)} - \| M_{1/\sqrt{b}} M_{{\bf r_2}} \cdot \nabla_x f\|_{L^2 (\R^d)}^2 \geq c \| \nabla_x f \|_{L^2 (\R^d)}^2$ for some $c>0$ by the ellipticity condition \eqref{ellipticity}. The details are omitted here. We can now conclude that  $\{e^{-t\mathcal{P}_{\mathcal{A}}}\}_{t\geq 0}$ is extended as a strongly continuous analytic semigroup acting on $L^2 (\R^{d})$ by Propositions  \ref{prop.poisson.analytic.L^2}, \ref{prop.subsec.relation.2}, and \ref{prop.embedding.section3.1} since $\Lambda_{\mathcal{A}}=\Lambda_{\mathcal{A}}$ when $A$ is Hermite. As for \eqref{est.prop.poisson.hermite.4}, we note that $\partial_t \Psi_{\mathcal{A}}[F], ~\mathcal{P}_{\mathcal{A}} \Psi_{\mathcal{A}}[F]\in L^2 (\R_+; H^{1/2}(\R^{d}))$ from the maximal regularity in the Hilbert space and $\nabla \Psi_{\mathcal{A}}[F]$ also belongs to $L^2 (\R^{d+1}_+)$ since $F$ is compactly supported. Hence, by using $\displaystyle \lim_{t\rightarrow 0} \Psi_{\mathcal{A}}[F] (t)=0$ in $H^{1/2}(\R^d)$ we see from Theorem \ref{thm.factorization} that
\begin{align*}
\langle A\nabla \Psi_{\mathcal{A}} [F], \nabla \Psi_{\mathcal{A}} [F]\rangle _{L^2 (\R^{d+1}_+)} & = \langle  M_{b} (\partial_t + \mathcal{P}_{\mathcal{A}} ) \Psi_{\mathcal{A}} [F], (\partial_t + \mathcal{P}_{\mathcal{A}} ) \Psi_{\mathcal{A}} [F] \rangle _{L^2 (\R^{d+1}_+)} \\
& = \| M_{\sqrt{b}} F \|_{L^2(\R^{d+1}_+)}^2.
\end{align*}
This proves \eqref{est.prop.poisson.hermite.4}, which is now extended to any $F \in L^2(\R^{d+1}_+)$. To prove \eqref{est.prop.poisson.hermite.3} set $w (t) = e^{-t\mathcal{P}_{\mathcal{A}}} g= g - \mathcal{P}_{\mathcal{A}} \int_0^t e^{-(t-s)\mathcal{P}_{\mathcal{A}}} g \dd s$. Fix any $T>0$. Then, as a consequence of \eqref{est.prop.poisson.hermite.4}, we have 
\begin{align}
\| w \|_{L^2(0,T; L^2 (\R^d))} \leq \| g \|_{L^2 (0,T; L^2 (\R^d ))} +  \| \mathcal{P}_{\mathcal{A}} \Psi_{\mathcal{A}}[\chi_T g]\|_{L^2 (0,T; L^2 (\R^d))} \leq  C T^{1/2} \| g \|_{L^2(\R^d)}, \label{proof.prop.poisson.hermite.1}
\end{align}
where $\chi_T=\chi_T(t)$ is the characteristic function on the interval $[0,T]$ and  $C>0$ is independent of $T$ and $g$. Thus there is $T_0\in [T/2,T]$ such that $\| w (T_0)\|_{L^2 (\R^d)} \leq C \| g \|_{L^2 (\R^d)}$ with $C$ independent of $T$ and $g$. 
Set $w_1 (t) = t\mathcal{P}_{\mathcal{A}} w (t)$. Then $w_1$ satisfies $\partial_t w_1 + \mathcal{P}_{\mathcal{A}} w_1 =  \mathcal{P}_{\mathcal{A}} w$  with the zero initial data.  Hence we have the representation $w_1(t) =  \int_0^t e^{-(t-s)\mathcal{P}_{\mathcal{A}}} \mathcal{P}_{\mathcal{A}} w (s) \dd s =  \mathcal{P}_{\mathcal{A}} \int_0^t e^{-(t-s)\mathcal{P}_{\mathcal{A}}} w (s) \dd s$. Thus \eqref{est.prop.poisson.hermite.4} yields 
\begin{align}
\| w_1 \| _{L^2(0,T; L^2 (\R^d))} \leq C \| w \|_{L^2 (0,T; L^2 (\R^d ))}\leq C T^\frac12 \| g\|_{L^2(\R^d)},~~~~~~~~~~T>0. \label{proof.prop.poisson.hermite.2}
\end{align}
By the identity $w (T) = e^{-(t-T_0)\mathcal{P}_{\mathcal{A}}} w (T_0) = w (T_0) -  \int_{T_0}^T \mathcal{P}_{\mathcal{A}}w (t) \dd t$ we have from \eqref{proof.prop.poisson.hermite.1} and \eqref{proof.prop.poisson.hermite.2},
\begin{align*}
\| w (T) \|_{L^2 (\R^d)} & \leq \| w (T_0) \|_{L^2 (\R^d)} + \int_{T_0}^T t^{-1} \| w_1 (t) \|_{L^2 (\R^d)} \dd t\\
& \leq  C \| g\|_{L^2 (\R^d)} + T_0^{-\frac12} \| w_1 \| _{L^2(0,T; L^2 (\R^d))} \leq C \| g \|_{L^2 (\R^d)}.
\end{align*}
Here we have used $T_0\in [T/2,T]$. Since $T>0$ is arbitrary we have shown that $\{ e^{-t\mathcal{P}_{\mathcal{A}}}\}_{t\geq 0}$ is a bounded semigroup in $L^2 (\R^d)$. The uniform bound estimate for $\{t\mathcal{P}_{\mathcal{A}} e^{-t\mathcal{P}_{\mathcal{A}}}\}_{t> 0}$ then follows from Proposition \ref{prop.poisson.analytic.L^2}, and \eqref{est.prop.poisson.hermite.3} is now proved. The proof is complete.

\begin{thm}\label{thm.poisson.hermite} Assume that $A$ is Hermite. Then the Poisson semigroup $\{e^{-t\mathcal{P}_{\mathcal{A}}}\}_{t\geq 0}$ in $H^{1/2}(\R^d)$ is extended as a bounded strongly continuous analytic semigroup in $L^2 (\R^d)$, and  $D_{L^2}(\mathcal{P}_{\mathcal{A}}) = D_{L^2}(\Lambda_{\mathcal{A}})= H^1 (\R^d)$ with equivalent norms. Moreover, the estimates \eqref{est.prop.poisson.hermite.1} - \eqref{est.prop.poisson.hermite.4} hold for all $f\in H^1 (\R^d)$, $g\in L^2 (\R^d)$, and $F\in L^2 (\R^{d+1}_+)$.

\end{thm}

\noindent {\it Proof.} In the proof of Lemma \ref{lem.poisson.hermite}, as a result of the estimate \eqref{est.prop.poisson.hermite.2} and Proposition \ref{prop.embedding.section3.1}, we have already seen that  $\{e^{-t\mathcal{P}_{\mathcal{A}}}\}_{t\geq 0}$ is extended as a strongly continuous analytic semigroup in $L^2 (\R^d)$, which is bounded due to \eqref{est.prop.poisson.hermite.3}. Moreover, Proposition \ref{prop.embedding.section3.1} also gives the continuous embeddings $D_{L^2}(\mathcal{P}_{\mathcal{A}})\hookrightarrow H^1 (\R^d)\hookrightarrow D_{L^2}(\Lambda_{\mathcal{A}})$, and \eqref{est.prop.poisson.hermite.2} is valid at least for $u\in D_{L^2}(\mathcal{P}_{\mathcal{A}})$. In view of Proposition \ref{prop.embedding.section3.2}, to show $D_{L^2}(\mathcal{P}_{\mathcal{A}})=H^1 (\R^d)$ it suffices to show $\sup_{0<t<1} \| \dd/\dd t ~e^{-t\mathcal{P}_{\mathcal{A}}} f \|_{L^2 (\R^d)} <\infty$ for any $f\in H^1 (\R^d)$. For this purpose we first consider the mollified matrix $A_\epsilon = ( a_{i,j}^{(\epsilon)} )_{1\leq i,j,\leq d+1}$, where $a_{i,j}^{(\epsilon)} = j_\epsilon * a_{i,j}$ with a standard mollifier $j_\epsilon$, $\epsilon>0$. It is easy to see that $A_\epsilon$ is Hermite and also satisfies \eqref{ellipticity} with the same constants $\nu_1,\nu_2$. Hence by \cite[Theorem 1.2]{MM1} we have the characterization $D_{L^2}(\mathcal{P}_{\mathcal{A}_\epsilon}) = H^1 (\R^d)$, and \eqref{est.prop.poisson.hermite.1} implies the estimate
\begin{align}
\frac{\nu_1}{\nu_2} \| \nabla_x f \|_{L^2 (\R^d)}^2 \leq  \| \mathcal{P}_{\mathcal{A}_\epsilon} f \|_{L^2 (\R^d)}^2 \leq \frac{\nu_2}{\nu_1} \| \nabla_x f \|_{L^2 (\R^d)}^2,~~~~~~~~~f\in H^1 (\R^d).\label{proof.thm.poisson.hermite.1}
\end{align}
Then the function $u_\epsilon (t) = e^{-t\mathcal{P}_{\mathcal{A}_\epsilon}} f$ with $f\in H^1 (\R^d)$ satisfies 
\begin{align}
\|\nabla  u_\epsilon (t) \|_{L^2 (\R^d)} \leq C \| \mathcal{P}_{\mathcal{A}_\epsilon} u_\epsilon (t) \|_{L^2 (\R^d)}  = C \| e^{-t\mathcal{P}_{\mathcal{A}_\epsilon}} \mathcal{P}_{\mathcal{A}_\epsilon} f \|_{L^2 (\R^d)} \leq C \| \mathcal{P}_{\mathcal{A}_\epsilon} f \|_{L^2 (\R^d)} \leq C \|\nabla_x f \|_{L^2 (\R^d)},\label{proof.thm.poisson.hermite.2}
\end{align}
as well as $ \sup_{t>0} \| u_\epsilon (t) \|_{L^2 (\R^d)} \leq C \| f\|_{L^2 (\R^d)}$. Here the constant $C$ in these estimates depends only on $\nu_1$ and $\nu_2$. Hence we can take a sequence $\{\epsilon_n\}$, $\epsilon_n\rightarrow 0$, such that $u_{\epsilon_n}$  converges to some $u\in  \dot{H}^1 (\R^{d+1}_+)$ weakly in $\dot{H}^1 (\R^{d+1}_+)$ and strongly in $L^2_{loc} (\R^{d+1}_+)$, and $u$ also satisfies the estimate $\sup_{t>0}\| \nabla u (t) \|_{L^2 (\R^d)}<\infty$. Since $a_{i,j}^{(\epsilon)}$  converges to $a_{i,j}$ strongly in $L^2_{loc}(\R^{d+1}_+)$,  we have $0=\langle A_{\epsilon_n} \nabla u_{\epsilon_n}, \nabla \varphi \rangle_{L^2 (\R^{d+1}_+)} \rightarrow  \langle A \nabla u, \nabla \varphi \rangle_{L^2 (\R^{d+1}_+)}$ for any $\varphi\in C_0^\infty (\R^{d+1}_+)$, and it is not difficult to see $\lim_{t\rightarrow 0} u(t) \rightarrow f$ strongly in $L^2 (\R^d)$ from the above uniform bound. Hence the uniqueness of the weak solutions 
in $\dot{H}^1 (\R^{d+1}_+)$ to \eqref{eq.dirichlet} 
implies $u(t) = e^{-t\mathcal{P}_{\mathcal{A}}} f$, i.e., $\sup_{t>0}\| \nabla e^{-t\mathcal{P}_{\mathcal{A}}} f\|_{L^2 (\R^d)}<\infty$ holds. thus $D_{L^2}(\mathcal{P}_{\mathcal{A}})=H^1 (\R^d)$ holds with equivalent norms by Proposition \ref{prop.embedding.section3.2}. Finally we show $D_{L^2}(\Lambda_{\mathcal{A}})=H^1 (\R^d)$, which requires more technicality. Let us recall that we have already known that $H^1(\R^d) \subset D_{L^2}(\Lambda_{\mathcal{A}})$ and $\|\nabla_x f\|_{L^2(\R^d)} \simeq \|\Lambda_{\mathcal{A}} f \|_{L^2(\R^d)}$ holds for any $f\in H^1 (\R^d)$. 
Now  let us consider the approximation $A_\epsilon$ as above, and for any $f\in L^2 (\R^d)$ set  $g_\epsilon=(\lambda + \Lambda_{\mathcal{A}_\epsilon})^{-1} f$. Since $A_\epsilon$ is Lipschitz we see $g_\epsilon\in D_{L^2}(\Lambda_{\mathcal{A}_\epsilon})=H^1 (\R^d)$ by \cite[Theorem 1.3]{MM1}. Moreover, using the fact that the ellipticity constants of $A_\epsilon$ is the same as those of $A$, we have $\| g_\epsilon\|_{H^1(\R^d)}\leq C(1 +\lambda^{-1}) \| f \|_{L^2 (\R^d)}$, where $C>0$ is independent of $\epsilon$.  Hence by taking suitable subsequence if necessary, we may assume that $g_\epsilon$ converges to some $g\in H^1 (\R^d)$ weakly in $L^2 (\R^d)$ and $\|g\|_{H^1 (\R^d)}\leq C (1+\lambda^{-1}) \| f \|_{L^2 (\R^d)}$ holds.  Thus for $\varphi\in H^1 (\R^d)$ we have 
\begin{align}
\langle g_\epsilon, (\lambda + \Lambda_{\mathcal{A}}) \varphi \rangle _{L^2 (\R^d)} \rightarrow  \langle g, (\lambda + \Lambda_{\mathcal{A}}) \varphi \rangle _{L^2 (\R^d)} = \langle (\lambda + \Lambda_{\mathcal{A}})  g, \varphi \rangle _{L^2 (\R^d)}  ~~~~~~~~~~\epsilon\rightarrow 0.\label{proof.thm.poisson.hermite.3}
\end{align}
On the other hand, we have 
\begin{align*}
\langle g_\epsilon,  (\lambda + \Lambda_{\mathcal{A}}) \varphi \rangle _{L^2 (\R^d)} 
& = \lambda\langle g_\epsilon, \varphi \rangle _{L^2 (\R^d)}  + \int_0^\infty \langle \nabla e^{-t\mathcal{P}_{\mathcal{A}_\epsilon}} g_\epsilon, A\nabla e^{-t\mathcal{P}_{\mathcal{A}}} \varphi\rangle _{L^2 (\R^d)}\dd t\\
& =  \lambda\langle g_\epsilon, \varphi \rangle _{L^2 (\R^d)}  + \int_0^\infty \langle A_\epsilon \nabla e^{-t\mathcal{P}_{\mathcal{A}_\epsilon}} g_\epsilon, \nabla e^{-t\mathcal{P}_{\mathcal{A}}} \varphi \rangle _{L^2 (\R^d)}\dd t\\
& ~~~~~+ \int_0^\infty \langle (A-A_\epsilon) \nabla e^{-t\mathcal{P}_{\mathcal{A}_\epsilon}} g_\epsilon, \nabla e^{-t\mathcal{P}_{\mathcal{A}}} \varphi \rangle _{L^2 (\R^d)}\dd t\\
& = \lambda\langle g_\epsilon, \varphi \rangle _{L^2 (\R^d)}  +  \langle \Lambda_{\mathcal{A}_\epsilon} g_\epsilon,  \varphi \rangle _{L^2 (\R^d)}\\
& ~~~~~+ \int_0^\infty \langle (A-A_\epsilon) \nabla e^{-t\mathcal{P}_{\mathcal{A}_\epsilon}} g_\epsilon, \nabla e^{-t\mathcal{P}_{\mathcal{A}}} \varphi \rangle _{L^2 (\R^d)}\dd t\\
& = \langle f, \varphi  \rangle _{L^2 (\R^d)} + \int_0^\infty \langle (A-A_\epsilon) \nabla e^{-t\mathcal{P}_{\mathcal{A}_\epsilon}} g_\epsilon, \nabla e^{-t\mathcal{P}_{\mathcal{A}}} \varphi \rangle _{L^2 (\R^d)}\dd t.
\end{align*}
Now we observe that $A_\epsilon\rightarrow A$ strongly in $L^p(B_R)$ for any $R>0$ and $p\in [1,\infty)$. Here $B_R=\{x\in \R^d~|~|x|<R\}$. Fix $R>0$. Then for each $\eta>0$ we have $\displaystyle\lim_{\epsilon\rightarrow 0} |\{x\in B_R~|~|A-A_\epsilon|\geq \eta\} |=0$. Thus we see from $\| \nabla e^{-\cdot \mathcal{P}_{\mathcal{A}_\epsilon}} g_\epsilon \|_{L^2 (\R^{d+1}_+)}\leq  C \|g_\epsilon\|_{\dot{H}^{1/2} (\R^d)}\leq C_\lambda \| f \|_{L^2 (\R^d)}$,
\begin{align}
& ~~~ \limsup_{\epsilon\rightarrow 0}|\int_0^\infty \langle (A-A_\epsilon) \nabla e^{-t\mathcal{P}_{\mathcal{A}_\epsilon}} g_\epsilon, \nabla e^{-t\mathcal{P}_{\mathcal{A}}} \varphi \rangle _{L^2 (\R^d)}\dd t|\nonumber\\
& \leq \limsup_{\epsilon\rightarrow 0} \int_0^\infty | \langle (A-A_\epsilon) \nabla e^{-t\mathcal{P}_{\mathcal{A}_\epsilon}} g_\epsilon, \nabla e^{-t\mathcal{P}_{\mathcal{A}}} \varphi \rangle _{L^2 (B_R\cap \{|A-A_\epsilon|<\eta\})} |\dd t \nonumber \\
& ~~~+ \limsup_{\epsilon\rightarrow 0} \int_0^\infty | \langle (A-A_\epsilon) \nabla e^{-t\mathcal{P}_{\mathcal{A}_\epsilon}} g_\epsilon, \nabla e^{-t\mathcal{P}_{\mathcal{A}}} \varphi \rangle _{L^2 (B_R\cap \{|A-A_\epsilon|\geq \eta\})} |\dd t \nonumber \\
&~~~~~~~ + \limsup_{\epsilon\rightarrow 0} \int_0^\infty | \langle (A-A_\epsilon) \nabla e^{-t\mathcal{P}_{\mathcal{A}_\epsilon}} g_\epsilon, \nabla e^{-t\mathcal{P}_{\mathcal{A}}} \varphi \rangle _{L^2 (\R^d\setminus B_R))} |\dd t \nonumber  \\
& \leq C_\lambda \eta \| f \|_{L^2 (\R^d)}  \| \varphi \|_{\dot{H}^\frac12 (\R^d)} + C_\lambda \| f \|_{L^2 (\R^d)}  \limsup_{\epsilon\rightarrow 0} \big (\int_0^\infty  \| \nabla e^{-t\mathcal{P}_{\mathcal{A}}} \varphi \|_{L^2 (B_R\cap \{|A-A_\epsilon|\geq \eta\})}^2 \dd t\big )^\frac12  \nonumber \\
& ~~~ + C_\lambda \| f \|_{L^2 (\R^d)} \big (\int_0^\infty  \| \nabla e^{-t\mathcal{P}_{\mathcal{A}}} \varphi \|_{L^2 (\R^d\setminus B_R)}^2 \dd t\big )^\frac12 \nonumber \\
& =  C_\lambda \eta \| f \|_{L^2 (\R^d)}  \| \varphi \|_{\dot{H}^\frac12 (\R^d)}  +  C_\lambda \| f \|_{L^2 (\R^d)} \big (\int_0^\infty  \| \nabla e^{-t\mathcal{P}_{\mathcal{A}}} \varphi \|_{L^2 (\R^d\setminus B_R)}^2 \dd t\big )^\frac12\nonumber. 
\end{align}
Since $\eta>0$ and $R>0$ are arbitrary, we conclude that
\begin{align*}
\limsup_{\epsilon\rightarrow 0} \big |\int_0^\infty \langle (A-A_\epsilon) \nabla e^{-t\mathcal{P}_{\mathcal{A}_\epsilon}} g_\epsilon, \nabla e^{-t\mathcal{P}_{\mathcal{A}}} \varphi \rangle _{L^2 (\R^d)}\dd t \big | =0.
\end{align*}
Collecting these above, we finally obtain
\begin{align}
\langle g_\epsilon, (\lambda + \Lambda_{\mathcal{A}})\varphi \rangle _{L^2 (\R^d)} ~~~ \rightarrow  ~~~ \langle f, \varphi \rangle_{L^2 (\R^d)} ~~~~~~\epsilon \rightarrow 0.\label{proof.thm.poisson.hermite.4}
\end{align}
Hence, \eqref{proof.thm.poisson.hermite.1} and \eqref{proof.thm.poisson.hermite.4} imply
\begin{align}
\langle ( \lambda+  \Lambda_{\mathcal{A}}) g, \varphi\rangle _{L^2 (\R^d)} =\langle g,  (\lambda +  \Lambda_{\mathcal{A}} )  \varphi \rangle _{L^2 (\R^d)} = \langle f, \varphi \rangle _{L^2 (\R^d)}, ~~~~~~~~~~\varphi\in H^1 (\R^d),
\end{align}
which gives $(\lambda + \Lambda_{\mathcal{A}}) g = f$.  Since $f\in L^2 (\R^d)$ is arbitrary and $g\in H^1 (\R^d)$, we have $D_{L^2}(\Lambda_{\mathcal{A}})\subset H^1 (\R^d)$, as desired. The proof is complete.


\subsection{Analysis for  regular off-block case}\label{subsec.poisson.regular}

In this section we consider the case  when the off-block vectors ${\bf r_1}$, ${\bf r_2}$ possess additional properties. As stated in Theorem \ref{thm.factorization.strong.1}, the key assumption will be made for the divergence of ${\bf r_1}$, ${\bf r_2}$. 
\begin{prop}\label{thm.domain.embedding} Assume that  $\{e^{-t\mathcal{P}_{\mathcal{A}}}\}_{t\geq 0}$ acting on $H^{1/2}(\R^d)$ is extended as a  strongly continuous semigroups acting on $L^2 (\R^d)$ and that

\noindent {\rm (i)} there are $\lambda\geq 0$ and $C_1>0$ such that 
\begin{equation}
\displaystyle \int_0^\infty \| e^{-t\mathcal{P}_{\mathcal{A}} -\lambda t} f \|_{\dot{H}^{\frac12}(\R^d)}^2 \dd t \leq C_1 \| f \|_{L^2(\R^d)}^2, ~~~~~~~~\quad\quad  f\in H^\frac12 (\R^d),\label{cond.thm.domain.embedding.1}
\end{equation} 
\noindent {\rm (ii)} there are $\alpha\in [0,1),~ C_2>0$, and $C_3>0$ such that 
\begin{equation}
| \langle M_{\nabla_x\cdot {\bf r_1}} f, f \rangle _{L^2 (\R^d)} |\leq C_2 \| f \|_{\dot{H}^\frac{1+\alpha}{2} (\R^d )} \| f \|_{H^\frac{1-\alpha}{2} (\R^d )} + C_3 \|f \|_{L^2 (\R^d)}^2,~ ~~~~~\quad  f \in H^1 (\R^d).\label{cond.thm.domain.embedding.2}
\end{equation}
Then $D_{L^2}(\mathcal{P}_{\mathcal{A}}) \subset H^1(\R^d)$ and 
\begin{align}
\| \nabla_x f \|_{L^2(\R^d)}\leq C (\|\mathcal{P}_{\mathcal{A}} f \|_{L^2 (\R^d)} + \| f \|_{L^2 (\R^d)} ),~~~~~~~~~f\in D_{L^2} (\mathcal{P}_{\mathcal{A}}),\label{est.thm.domain.embedding}
\end{align}
where the constant $C$ depends only on $d$, $\nu_1$, $\nu_2$,  $\alpha$, $C_1$, $C_2$, and $C_3$. 
\end{prop}


\begin{rem}{\rm As is seen in the proof below, the constant $C$ in \eqref{est.thm.domain.embedding} does not depend on the value $\displaystyle \sup_{t>0}\|e^{-t\mathcal{P}_{\mathcal{A}}-t}\|_{L^2\rightarrow L^2}$.  

}
\end{rem}

\begin{rem}{\rm Let $d\geq 2$ and let $h=h_1 + h_2\in L^{d,\infty}(\R^d) + L^\infty(\R^d)$. Then we have 
\begin{align}
| \langle M_h f, f \rangle _{L^2(\R^d)} |\leq C \| h_1 \|_{L^{d,\infty} (\R^d)} \| f \|_{\dot{H}^\frac12 (\R^d )}^2 +  \| h_2 \|_{L^\infty (\R^d)} \| f \|_{L^2 (\R^d)}^2, ~~~~~\quad \quad f \in H^{\frac12} (\R^d).
\end{align}
Hence, if $\nabla_x\cdot {\bf r_1}$ and $\nabla_x\cdot {\bf r_2}$ belong to $L^{d,\infty}(\R^d) + L^\infty(\R^d)$ and $d\geq 2$ then the condition (ii) in Proposition \ref{thm.domain.embedding} is satisfied. Similarly, when $d=1$  the condition (ii) in Proposition \ref{thm.domain.embedding} holds if $\nabla_x\cdot {\bf r_1}$ and $\nabla_x\cdot {\bf r_2}$ belong to $\mathcal{M}(\R) + L^\infty(\R)$.
}

\end{rem}

To prove Proposition \ref{thm.domain.embedding}  let us recall the operator $\mathcal{A}': D(\mathcal{A}')\subset L^2 (\R^d)\rightarrow L^2 (\R^d)$, which is defined through the sesquilinear form $\langle \mathcal{A}' f,g\rangle _{L^2(\R^d)} = \langle A' \nabla_x f, \nabla_x g\rangle _{L^2(\R^d)}$ for $f\in D_{L^2} (\mathcal{A}')$ and $g\in H^1(\R^d)$. Since $A(x)$ satisfies the ellipticity condition \eqref{ellipticity}, so does $A'(x)$. Hence $\mathcal{A}'$ is realized as an {\it m}-sectorial operator in $L^2 (\R^d)$. Furthermore, due to the result  \cite{AHLMT}  on the Kato square root problem we have $D_{L^2}(\sqrt{\mathcal{A}'}) = H^1 (\R^d)$ and $\| \sqrt{\mathcal{A}'} f \| _{L^2 (\R^d)} \simeq  \|\nabla_x f \|_{L^2 (\R^d)}$. Reminding this fact, we start from the next lemma based on the Littlewood-Paley theory. Let $\psi\in C_0^\infty (\R^d)$ be a real-valued function with zero average such that 
\[
\int_0^\infty \| \psi_s * f \|_{L^2(\R^d)}^2 \frac{\dd s}{s}  = \| f\|_{L^2 (\R^d)}^2,~~~~~~~~~ f\in L^2 (\R^d).
\]
Here $\psi_s (x) = s^{-d} \psi (x/s)$. We may take $\psi=\Delta \tilde \psi$ so that $\| s^{-1} \nabla_x \tilde \psi_s * f \|_{L^2(\R^d)}\leq C \| f\|_{L^2 (\R^d)}$ holds. 
\begin{lem}\label{lem.thm.domain.embedding} Suppose that the assumptions in Proposition \ref{thm.domain.embedding} hold. Set 
\[
U_{\lambda'} (t) = t^\frac12 (-\Delta_x)^{-\frac14} (e^{-t\mathcal{P}_{\mathcal{A}}})^* e^{-\lambda' t} (\sqrt{\mathcal{A}'})^* ~~~~~~~~~~~~{\rm for}~~t>0,~ \lambda' \geq 0.
\]
Then there is $\lambda_0\geq \lambda$  such that 
\begin{align}
\int_0^\infty \| U_{\lambda_0} (t) g \|_{L^2(\R^d)}^2 \frac{\dd t}{t} \leq C \| g \|_{L^2(\R^d)}^2, ~~~~~~~~~~~~~~g\in L^2 (\R^d).\label{est.lem.thm.domain.embedding.3}
\end{align}
Here the constant $C$ depends only on $d$, $\nu_1$, $\nu_2$, $\alpha$, $C_2$, and $C_3$. 

\end{lem}

\noindent {\it Proof.} For simplicity we assume $\lambda=0$ in the condition of Proposition \ref{thm.domain.embedding} and take $\lambda_0=1$. We also assume, instead of \eqref{cond.thm.domain.embedding.2}, that 
\begin{equation}
| \langle M_{\nabla_x\cdot {\bf r_1}} f, f \rangle _{L^2 (\R^d)} |\leq C_2 \| f \|_{\dot{H}^\frac{1+\alpha}{2} (\R^d )} \| f \|_{\dot{H}^\frac{1-\alpha}{2} (\R^d )} + C_3 \|f \|_{L^2 (\R^d)}^2,~ ~~~~~ \quad\quad  f\in H^\frac12 (\R^d),\label{cond.thm.domain.embedding.2'}
\end{equation}
which makes the computation slightly simpler. The general case is treated just in the same manner. We write $U (t)$ for $U_1(t)$. By taking into account the Schur lemma (cf. see \cite[pp.643-644]{Grafakos}) it suffices to show
\begin{align}
\| U (t) g \|_{L^2(\R^d)} & \leq C \| g\|_{L^2(\R^d)},~~~~~~~~~~~~~~~~~~~~~~~~~~~~~~~~t>0,~~g\in H^1 (\R^d),\label{est.lem.thm.domain.embedding.1}\\
\| U (t) Q_s  g \|_{L^2 (\R^d)} & \leq C \big (\min\{ \frac{t}{s}, \frac{s}{t}\}  \big )^{\frac{1-\alpha}{2}} \| g \|_{L^2 (\R^d),}~~~~~~~~~~~~t,s>0,~~g\in L^2 (\R^d),\label{est.lem.thm.domain.embedding.2}
\end{align}
with the constant  $C$ depends only on $d$, $\nu_1$, $\nu_2$, $\alpha$,  $C_2$, and $C_3$. Here we have set $Q_s g = \psi_s * g$.
Let us prove the desired estimates by the duality. Take any $f\in L^2 (\R^d)$ and $g\in H^1 (\R^d)$ with $\|  f\|_{L^2 (\R^d)} = \| g\|_{L^2 (\R^d)} = 1$. Then we have 
\begin{align*}
|\langle f, U(t) g\rangle _{L^2(\R^d)}| & = t^\frac12  |\langle \sqrt{\mathcal{A}'} e^{-t\mathcal{P}_{\mathcal{A}}-t} (-\Delta_x)^{-\frac14} f, g \rangle _{L^2(\R^d)} | \\
& \leq C t^\frac12  \| \nabla_x    e^{-t\mathcal{P}_{\mathcal{A}}-t} (-\Delta_x)^{-\frac14} f \|_{L^2(\R^d)},
\end{align*}
for  $\| \sqrt{\mathcal{A}'} u \| _{L^2 (\R^d)} \simeq  \|\nabla_x u \|_{L^2 (\R^d)}$. Then from \eqref{est.proof.prop.poisson.analytic.4'}  we have  
\begin{align*}
\| \nabla_x    e^{-t\mathcal{P}_{\mathcal{A}}-t} (-\Delta_x)^{-\frac14} f \|_{L^2(\R^d)}\leq C t^{-\frac12}\| (-\Delta_x)^{-\frac14} f\|_{\dot{H}^\frac12 (\R^d)}  \leq C t^{-\frac12}.
\end{align*}
This proves \eqref{est.lem.thm.domain.embedding.1}. To show \eqref{est.lem.thm.domain.embedding.2} we recall the characterization 
\[
\| \sqrt{\mathcal{A}'} f \| _{L^2 (\R^d)} \simeq \| (\sqrt{\mathcal{A}'} )^* f \| _{L^2 (\R^d)} \simeq \|\nabla_x f \|_{L^2 (\R^d)}
\]
and the fact that the realization of $\mathcal{A}'$ in $L^2 (\R^d)$ has the bounded imaginary powers leads to the characterization 
\[
\| {\mathcal{A}'}^{\frac14} f \| _{L^2 (\R^d)} \simeq \| ({\mathcal{A}'}^{\frac14})^* f \| _{L^2 (\R^d)} \simeq \| f \|_{\dot{H}^{\frac12}(\R^d)}
\]
by the complex interpolation. By \eqref{est.proof.prop.poisson.analytic.1} we obtain
\begin{align*}
|\langle f, U(t) Q_s  g\rangle _{L^2(\R^d)}| & = t^\frac12  |\langle {\mathcal{A}'}^\frac14  e^{-t\mathcal{P}_{\mathcal{A}}-t} (-\Delta_x)^{-\frac14} f,  ({\mathcal{A}'}^\frac14)^* Q_s g \rangle _{L^2(\R^d)} | \\
& \leq C t^\frac12  \| e^{-t\mathcal{P}_{\mathcal{A}}-t} (-\Delta_x)^{-\frac14} f \|_{\dot{H}^\frac12 (\R^d)} \| Q_s  g\|_{\dot{H}^\frac12 (\R^d)} \leq C  t^\frac12 s^{-\frac12}.
\end{align*}
Therefore it remains to consider the case $t\geq s$. Proposition \ref{prop.thm.factorization} yields
\begin{align*}
|\langle f, U(t) Q_s g\rangle _{L^2(\R^d)}| & = t^\frac12  |\langle  e^{-t\mathcal{P}_{\mathcal{A}}-t} (-\Delta_x)^{-\frac14} f,  (\mathcal{A}')^* ({\mathcal{A}'}^{-\frac12})^* Q_s g \rangle _{L^2(\R^d)} | \\
& =   t^\frac12  | \langle  A' \nabla_x e^{-t\mathcal{P}_{\mathcal{A}}-t} (-\Delta_x)^{-\frac14} f,  \nabla_x  ({\mathcal{A}'}^{-\frac12})^* Q_s g \rangle _{L^2(\R^d)} | \\
& \leq  t^\frac12 | \langle \mathcal{P}_{\mathcal{A}}  e^{-t\mathcal{P}_{\mathcal{A}}-t} (-\Delta_x)^{-\frac14} f, \Lambda_{\mathcal{A}^*} ({\mathcal{A}'}^{-\frac12})^* Q_s g \rangle _{\dot{H}^{\frac12}, \dot{H}^{-\frac12}} |  \\
& ~~~~ +  t^\frac12  | \langle \mathcal{P}_{\mathcal{A}}  e^{-t\mathcal{P}_{\mathcal{A}}-t} (-\Delta_x)^{-\frac14} f, M_{\bf \bar{r}_1} \cdot \nabla_x  ({\mathcal{A}'}^{-\frac12})^* Q_s g \rangle _{L^2(\R^d)} |\\
& = I_1 + I_2.
\end{align*}
The term $I_1$ is estimated as 
\begin{align*}
I_1 & \leq C t^\frac12 \| \frac{\dd}{\dd t}  e^{-t\mathcal{P}_{\mathcal{A}}-t} (-\Delta_x)^{-\frac14} f \|_{\dot{H}^\frac12(\R^d)} \|   ({\mathcal{A}'}^{-\frac12})^* Q_s g \| _{\dot{H}^{\frac12}(\R^d)} \leq C t^{-\frac12} s^\frac12. 
 \end{align*}
 Here we have used \eqref{est.proof.prop.poisson.analytic.3} and $\|({\mathcal{A}'}^{-1/2})^* Q_s g \| _{\dot{H}^{1/2}(\R^d)}\leq C \| Q_s g \|_{\dot{H}^{-1/2} (\R^d)}\leq C s^{1/2}$. As for the term $I_2$, we perform the integration by parts and get
 \begin{align*}
 I_2 & \leq t^\frac12  \| {\bf r_1} \|_{L^\infty (\R^d)} \| \nabla_x e^{-\frac{t}{2} \mathcal{P}_{\mathcal{A}}} \mathcal{P}_{\mathcal{A}}  e^{-\frac{t}{2}\mathcal{P}_{\mathcal{A}}-t} (-\Delta_x)^{-\frac14} f \|_{L^2(\R^d)} \|  ({\mathcal{A}'}^{-\frac12})^* Q_s g \| _{L^2(\R^d)} \\
 & ~~~~~ +  t^\frac12  |  \langle \frac{\dd}{\dd t}  e^{-t\mathcal{P}_{\mathcal{A}}-t} (-\Delta_x)^{-\frac14} f, M_{\nabla_x\cdot {\bf \bar{r}_1}} ({\mathcal{A}'}^{-\frac12})^* Q_s g \rangle _{L^2(\R^d)}|\\
 & = I_{2,1} + I_{2,2}.
 \end{align*}
By using \eqref{est.proof.prop.poisson.analytic.3} and \eqref{est.proof.prop.poisson.analytic.4'} the term $I_{2,1}$ is bounded from above by  $C t^{-1} s$. To estimate the last term $I_{2,2}$ we recall the standard identity for a given inner product $\langle,\rangle$:
\begin{align}
\langle u,h v\rangle  & = \frac14 \big \{ \langle u+v, h (u + v)\rangle - \langle  u-v , h (u-v)\rangle  + i \langle u + i v, h ( u+iv ) \rangle - i \langle u - i v,  h ( u-i v )\rangle \big \}.\label{inner.product}
\end{align}
Then, applying \eqref{cond.thm.domain.embedding.2'}, we have for any $m>0$,
\begin{align*}
I_{2,2} & \leq C  t^\frac12 e^{-\frac{t}{2}} \big ( \| m  \frac{\dd}{\dd t} e^{-t\mathcal{P}_{\mathcal{A}}-\frac{t}{2}} (-\Delta_x)^{-\frac14} f\|_{\dot{H}^{\frac{1+\alpha}{2}} (\R^d)} + \| m^{-1} ({\mathcal{A}'}^{-\frac12})^* Q_s g \|_{\dot{H}^{\frac{1+\alpha}{2}} (\R^d)} \big )\\
& ~~~~~ \cdot \big ( \| m  \frac{\dd}{\dd t} e^{-t\mathcal{P}_{\mathcal{A}}-\frac{t}{2}} (-\Delta_x)^{-\frac14} f\|_{\dot{H}^{\frac{1-\alpha}{2}} (\R^d)} + \| m^{-1} ({\mathcal{A}'}^{-\frac12})^* Q_s g \|_{\dot{H}^{\frac{1-\alpha}{2}} (\R^d)} \big )\\
& ~~~~~~~ + C t^\frac12 e^{-\frac{t}{2}} \big ( \| m  \frac{\dd}{\dd t}  e^{-t\mathcal{P}_{\mathcal{A}}-\frac{t}{2}} (-\Delta_x)^{-\frac14} f\|_{L^2 (\R^d)} + \| m^{-1} ({\mathcal{A}'}^{-\frac12})^* Q_s g \|_{L^2 (\R^d)} \big )^2.
\end{align*}
The interpolation inequality with \eqref{est.proof.prop.poisson.analytic.3} and \eqref{est.proof.prop.poisson.analytic.4'} shows that 
\[
 \| \frac{\dd}{\dd t} e^{-t\mathcal{P}_{\mathcal{A}}-\frac{t}{2}} (-\Delta_x)^{-\frac14} f\|_{\dot{H}^{\beta} (\R^d)}\leq C t^{-\frac{1}{2}-\beta} \|f \|_{L^2(\R^d)},~~~~~ \| ({\mathcal{A}'}^{-\frac12})^* Q_s g \|_{\dot{H}^{\beta} (\R^d)}\leq C s^{1-\beta} \| g \|_{L^2(\R^d)},
\]
for $\beta\in [0,1]$. Thus we have 
\begin{align*}
I_{2,2} \leq C t^\frac12 e^{-\frac{t}{2}}\big ( m^2 t^{-2} + t^{-1+\frac{\alpha}{2}} s^{\frac{1-\alpha}{2}} + t^{-1 -\frac{\alpha}{2}} s^{\frac{1+\alpha}{2}} + m^{-2} s^2 \big ).
\end{align*} 
Hence, setting $m^4 = t^2 s$, we get $I_{2,,2}\leq C t^{-(1-\alpha)/2} s^{(1-\alpha)/2}$ for $t\geq s>0$, which proves \eqref{est.lem.thm.domain.embedding.2}. The proof is complete.

\vspace{0.5cm}


We are now in position to prove Proposition \ref{thm.domain.embedding}.

\noindent {\it Proof of Proposition \ref{thm.domain.embedding}.} As in the proof of  Lemma \ref{lem.thm.domain.embedding}, we assume that $\lambda=0$ and $\lambda_0=1$, and write $U(t)$ for $U_{\lambda_0}(t)$.  It suffices to show $\| \sqrt{\mathcal{A}'} (I + \mathcal{P}_{\mathcal{A}})^{-1} f \|_{L^2(\R^d)}\leq C \| f\|_{L^2 (\R^d)}$ for all  $f\in H^{1/2} (\R^d)$ by the density argument and by the fact $\| \sqrt{\mathcal{A}'} f \|_{L^2 (\R^d)} \simeq \| \nabla_x f \|_{L^2(\R^d)}$. We apply the duality argument. For any $g\in H^1 (\R^d)$ we have 
\begin{align*}
J & := \langle \sqrt{\mathcal{A}'} \int_0^\infty e^{-2 t\mathcal{P}_{\mathcal{A}} - 2 t } f \dd t, g\rangle_{L^2(\R^d)} \\
& = \int_0^\infty \langle (-\Delta_x)^\frac14  e^{-t\mathcal{P}_{\mathcal{A}} -t } f, (-\Delta_x)^{-\frac14} (e^{-t\mathcal{P}_{\mathcal{A}}})^* e^{-t} (\sqrt{A'})^* g\rangle _{L^2(\R^d)} \dd t.
\end{align*}
Thus the H{\"o}lder inequality implies  
\begin{align*}
|J| &\leq \big (\int_0^\infty \| e^{-t\mathcal{P}_{\mathcal{A}} -t } f\|_{\dot{H}^\frac12 (\R^d)}^2 \dd t \big )^\frac12 \big ( \int_0^\infty \| U (t) g \|_{L^2(\R^d)}^2 \frac{\dd t}{t} \big )^\frac12  \leq C \| f\|_{L^2(\R^d)} \| g\|_{L^2 (\R^d)}.
\end{align*}
Here we have used the condition (i) of Proposition \ref{thm.domain.embedding} and Lemma \ref{lem.thm.domain.embedding}. The proof is complete.

\begin{prop}\label{prop.poisson.regular.1} Let $\alpha\in [0,1)$. Assume that ${\bf r_2}$ and $b$ are real-valued. Then there is $\delta_0>0$ depending only on $\nu_1$, $\nu_2$, and $\alpha$ such that if
\begin{align}
| \langle M_{\nabla_x\cdot {\bf r_2}} f, f \rangle _{L^2 (\R^d)} |\leq \delta_0 \| f \|_{\dot{H}^\frac{1+\alpha}{2} (\R^d )} \| f \|_{H^\frac{1-\alpha}{2} (\R^d )} + C_4 \| f\|_{L^2 (\R^d)}^2,~ ~~~~~ \quad\quad  f\in H^1 (\R^d),\label{cond.prop.poisson.regular.1.1}
\end{align}
for some  $C_4>0$  then $\{ e^{-t\mathcal{P}_{\mathcal{A}}}\}_{t\geq 0}$ in $H^{1/2}(\R^d)$ is extended as a strongly continuous semigroup in $L^2 (\R^d)$ and it follows that
\begin{align}
 \| M_{\sqrt{b}} e^{-t\mathcal{P}_{\mathcal{A}}} f\|_{L^2(\R^d)}^2 + \nu_1 \int_0^t  \| e^{-s\mathcal{P}_{\mathcal{A}}} f \|_{\dot{H}^\frac12 (\R^d)}^2 \dd s \leq 2 e^{\frac{2 C_4 + \delta_0}{\nu_1}  t} \| M_{\sqrt{b}} f\|_{L^2(\R^d)}^2,~~~~~~~~~~~t>0. \label{est.prop.poisson.regular.1}
\end{align}

\end{prop}

\noindent {\it Proof.} As in the proof of Lemma \ref{lem.thm.domain.embedding}, we assume 
\begin{align}
| \langle M_{\nabla_x\cdot {\bf r_2}} f, f \rangle _{L^2 (\R^d)} |\leq \delta_0 \| f \|_{\dot{H}^\frac{1+\alpha}{2} (\R^d )} \| f \|_{\dot{H}^\frac{1-\alpha}{2} (\R^d )} + C_4 \| f\|_{L^2 (\R^d)}^2,~ ~~~~~ \quad\quad  f \in H^1 (\R^d),\label{cond.prop.poisson.regular.1.1'}
\end{align}
instead of \eqref{cond.prop.poisson.regular.1.1} for simplicity of computations.  Let $f\in H^{1/2}(\R^d)$.  Then by Proposition \ref{prop.poisson.analytic} and Proposition \ref{prop.subsec.relation.2}  the function $e^{-t\mathcal{P}_{\mathcal{A}}}f$ belongs to $D_{H^{1/2}}(\mathcal{P}_{\mathcal{A}})\subset D_{L^2}(\Lambda_{\mathcal{A}}) \cap H^1 (\R^d)$ for all $t>0$, and the following  calculation is valid.
\begin{align}
\frac{\dd }{\dd t} \| M_{\sqrt{b}} e^{-t\mathcal{P}_{\mathcal{A}}} f \|_{L^2 (\R^d) }^2 & =  -2 {\rm Re} \langle \Lambda_{\mathcal{A}} e^{-t\mathcal{P}_{\mathcal{A}}} f, e^{-t\mathcal{P}_{\mathcal{A}}} f \rangle _{L^2(\R^d)} -2{\rm Re}  \langle M_{{\bf r_2}} \cdot \nabla_x  e^{-t\mathcal{P}_{\mathcal{A}}} f, e^{-t\mathcal{P}_{\mathcal{A}}} f \rangle _{L^2(\R^d)} \nonumber \\
& = -  2 {\rm Re} \langle \Lambda_{\mathcal{A}} e^{-t\mathcal{P}_{\mathcal{A}}} f, e^{-t\mathcal{P}_{\mathcal{A}}} f \rangle _{L^2(\R^d)} +  \langle M_{\nabla_x\cdot {\bf r_2}} e^{-t\mathcal{P}_{\mathcal{A}}} f, e^{-t\mathcal{P}_{\mathcal{A}}} f \rangle_{L^2(\R^d)}.\label{proof.prop.poisson.regular.1}
\end{align}
Applying \eqref{cond.prop.poisson.regular.1.1'},  we then have 
\begin{align*}
\langle M_{\nabla_x\cdot {\bf r_2}} e^{-t\mathcal{P}_{\mathcal{A}}} f, e^{-t\mathcal{P}_{\mathcal{A}}} f \rangle_{L^2(\R^d)} \leq \delta \|e^{-t\mathcal{P}_{\mathcal{A}}} f \|_{\dot{H}^{\frac{1+\alpha}{2}} (\R^d)}  \|e^{-t\mathcal{P}_{\mathcal{A}}} f \|_{\dot{H}^{\frac{1-\alpha}{2}} (\R^d)}  + C_4 \| e^{-t\mathcal{P}_{\mathcal{A}}} f \|_{L^2(\R^d)}^2. 
\end{align*}
Fix $T>0$. Integrating \eqref{proof.prop.poisson.regular.1} from $t=0$ to $t=T$, we  get 
\begin{align}
& ~~~ \| M_{\sqrt{b}} e^{-T\mathcal{P}_{\mathcal{A}}} f \|_{L^2 (\R^d ) }^2 + 2 \int_0^T {\rm Re}  \langle \Lambda_{\mathcal{A}} e^{-t\mathcal{P}_{\mathcal{A}}} f, e^{-t\mathcal{P}_{\mathcal{A}}} f \rangle _{L^2(\R^d)} \dd t \nonumber \\
& \leq \| M_{\sqrt{b}} f \|_{L^2 (\R^d) }^2 + C \delta   \int_0^T t^{-\frac{\alpha}{2}} \| e^{-\frac{t}{2}\mathcal{P}_{\mathcal{A}}} f\|_{\dot{H}^{\frac12} (\R^d)} \| e^{-t\mathcal{P}_{\mathcal{A}}} f\|_{\dot{H}^{\frac{1-\alpha}{2}} (\R^d )} \dd t + C_4  \int_0^T \| e^{-t\mathcal{P}_{\mathcal{A}}} f\|_{L^2 (\R)}^2 \dd t.\label{proof.prop.poisson.regular.2}
\end{align}
Then we use Proposition \ref{prop.key.estimate.poisson.analytic} to conclude that  
\begin{align}
& ~~~ \| M_{\sqrt{b}} e^{-T\mathcal{P}_{\mathcal{A}}} f \|_{L^2 (\R^d ) }^2 + 2\nu_1  \int_0^T \| e^{-t\mathcal{P}_{\mathcal{A}}} f \| _{\dot{H}^\frac12 (\R^d)}^2 \dd t \nonumber \\
& \leq \| M_{\sqrt{b}} f \|_{L^2 (\R^d) }^2 +  \frac{C \delta T^{1-\alpha}}{1-\alpha} \| e^{-T \mathcal{P}_{\mathcal{A}}} f \|_{L^2 (\R^d)}^{2\alpha} \| e^{-T \mathcal{P}_{\mathcal{A}}} f \|_{\dot{H}^\frac12 (\R^d)}^{2 (1-\alpha)} +  C_\alpha \delta \int_0^{T} \| e^{-t \mathcal{P}_{\mathcal{A}}} f\|_{\dot{H}^{\frac12} (\R^d)}^2 \dd t \nonumber \\ 
& ~~~~~ + C_4  \int_0^T \| e^{-t\mathcal{P}_{\mathcal{A}}} f\|_{L^2 (\R^d)}^2 \dd t.\label{proof.prop.poisson.regular.3}
\end{align}
Now we note that there is $T_0\in [T/2, T]$ such that 
\[
T_0   \|  e^{-T_0\mathcal{P}_{\mathcal{A}}} f\|_{\dot{H}^{\frac12} (\R^d )}^2\leq  4 \int_0^T \|  e^{-t\mathcal{P}_{\mathcal{A}}} f\|_{\dot{H}^{\frac12} (\R^d )}^2 \dd t, ~~~~ T \| e^{-T \mathcal{P}_{\mathcal{A}}} f \|_{\dot{H}^{\frac12} (\R^d)}^2 \leq  C  T_0 \| e^{-T_0 \mathcal{P}_{\mathcal{A}}} f \|_{\dot{H}^{\frac12} (\R^d)}^2.
\]
Thus, taking $\delta$ small enough if necessary, we arrive at
\begin{align}
\| M_{\sqrt{b}} e^{-T\mathcal{P}_{\mathcal{A}}} f \|_{L^2 (\R^d) }^2 +  \nu_1 \int_0^T \|  e^{-t\mathcal{P}_{\mathcal{A}}} f\|_{\dot{H}^\frac12 (\R^d )}^2 \dd t &  \leq 2\| M_{\sqrt{b}} f \|_{L^2 (\R^d) }^2 + \frac{2 C_4 }{\nu_1}  \int_0^T \| M_{\sqrt{b}} e^{-t\mathcal{P}_{\mathcal{A}}} f\|_{L^2 (\R^d)}^2 \dd t.\label{proof.prop.poisson.regular.5}
\end{align}
Since $T>0$ is arbitrary the Gronwall inequality  leads to \eqref{est.prop.poisson.regular.1}. The proof is complete.


\begin{thm}\label{thm.poisson.regular.off} Let $\alpha\in [0,1)$ and let $\delta_0$ be the number in Proposition \ref{prop.poisson.regular.1}. Then there is $\delta\in (0,\delta_0]$ depending only on $d$, $\alpha$, $\nu_1$, and $\nu_2$ such that the following statement holds: Suppose that

\vspace{0.3cm}

\noindent  ${\rm (1)} ~ |\langle M_{\nabla_x\cdot {\bf r_1}} f, f \rangle _{L^2 (\R^d)} | +  |\langle M_{\nabla_x\cdot  {\bf r_2}} f, f \rangle _{L^2 (\R^d)} | \leq  \delta \| f \|_{\dot{H}^\frac{1+\alpha}{2} (\R^d )} \| f\|_{H^\frac{1-\alpha}{2} (\R^d )} + C_5 \| f \|_{L^2 (\R^d)}^2,$

\noindent  ${\rm (2)} ~  {\rm Im} ~({\bf r_1 + r_2}) ={\bf 0}~{\rm and}~ {\rm Im} ~b =0.$

\vspace{0.3cm}

\noindent Then the semigroups $\{ e^{-t\mathcal{P}_{\mathcal{A}}}\}_{t\geq 0}$ and $\{ e^{-t\mathcal{P}_{\mathcal{A}^*}}\}_{t\geq 0}$   in $H^{1/2}(\R^d)$ are extended as strongly continuous semigroups in $L^2 (\R^d)$, and $D_{L^2}(\mathcal{P}_{\mathcal{A}}) = D_{L^2} (\mathcal{P}_{\mathcal{A}^*}) = H^1 (\R^d)$ holds with equivalent norms.
 
\end{thm}


\noindent {\it Proof.} To simplify the computation we assume, instead of (1), that
\begin{align}
|\langle M_{\nabla_x\cdot {\bf r_1}} f, f \rangle _{L^2 (\R^d)} | +  |\langle M_{\nabla_x\cdot  {\bf r_2}} f, f \rangle _{L^2 (\R^d)} | \leq  \delta \| f \|_{\dot{H}^\frac{1+\alpha}{2} (\R^d )} \| f \|_{\dot{H}^\frac{1-\alpha}{2} (\R^d )} + C_5 \| f \|_{L^2 (\R^d)}^2,\label{proof.thm.poisson.regular.off.1}
\end{align}
where $\delta\in (0,\delta_0]$ will be chosen later.  We denote by $A_1$ the matrix defined as
\begin{equation}\label{proof.prop.poisson.regular.off.2}
A_1 =
\begin{pmatrix}
\mbox{} & \mbox{} & \mbox{} & \mbox{} \\
\mbox{} & \large{A'} & \mbox{} & {\rm Re} ~{\bf r_1} \\
\mbox{} & \mbox{} & \mbox{} & \mbox{}\\
\mbox{} & {\rm Re} ~{\bf r_2}^\top & \mbox{} & b\end{pmatrix}~,~~~~~~~~~A' = (a_{i,j})_{1\leq i,j\leq d}.
\end{equation}
It is easy to see that $A_1$ satisfies \eqref{ellipticity} with the same numbers $\nu_1,\nu_2$. Let $t_0$ be a positive number satisfying $t_0\leq \min\{1, \nu_1/ (2 C_5 + \delta_0)\}$. Then \eqref{est.prop.poisson.regular.1} and \eqref{est.prop.poisson.analytic.L^2} yield the estimate
\begin{align}
t^\beta \big ( \| e^{-t\mathcal{P}_{\mathcal{A}_1}} f \|_{\dot{H}^\beta (\R^d)} + \| e^{-t\mathcal{P}_{\mathcal{A}_1^*}} f \|_{\dot{H}^\beta (\R^d)}\big ) \leq C \| f \|_{L^2 (\R^d)}, ~~~~~~~0<t\leq \frac{t_0}{2} \label{proof.thm.poisson.regular.off.3}
\end{align}
for $\beta \in [0,1]$, where $C$ depends only on $d$, $\alpha$, $\nu_1$, and $\nu_2$. 

For the moment we consider the case when each $a_{i,j}$ belongs to ${\rm Lip}(\R^d)$. Then $D_{L^2}(\mathcal{P}_{\mathcal{A}})=D_{L^2}(\mathcal{P}_{\mathcal{A}^*})=H^1 (\R^d)$ holds with equivalent norms; see \cite[Theorem 1.2]{MM1}. Let $f\in H^1 (\R^d)$. Then $u(t) = e^{-t\mathcal{P}_{\mathcal{A}}} f$ satisfies 
\begin{equation}
\mathcal{A}_1 u =  i \big ( \nabla_x \cdot M_{{\rm Im} {\bf r_1}} \partial_t u  + \partial_t   M_{{\rm Im} {\bf r_2}} \cdot \nabla_x u \big )\label{proof.thm.poisson.regular.off.4}
\end{equation}
in the weak sense and $u|_{t=0} = f$, where $\mathcal{A}_1 u=  - \nabla\cdot A_1 \nabla u$. Since  the right-hand side of \eqref{proof.thm.poisson.regular.off.4} satisfies the condition in Theorem \ref{thm.representation} we have the integral representation of $u$ such as 
\[ 
u(t) = e^{-t\mathcal{P}_{\mathcal{A}_1}} f +  i \sum_{k=1}^2  F_k (t).
\]
Here, for $0<t\leq t_0/2$ and $k=1,2$, each $F_k$ is defined by 
\begin{align*}
F_k (t) & = \int_0^t e^{-(t-s)\mathcal{P}_{\mathcal{A}_1}} G_k (s) \dd s,\\
G_1 (s) & = \int_{t_0}^\infty   e^{-(\tau-s)\mathcal{Q}_{\mathcal{A}_1}}  M_{1/b}   \big ( \nabla_x \cdot M_{{\rm Im} {\bf r_1}} \partial_\tau  u  + \partial_\tau   M_{{\rm Im} {\bf r_2}} \cdot \nabla_x u \big ) \dd \tau, \\
G_2 (s) & =-\int_s^{t_0}  e^{-(\tau-s)\mathcal{Q}_{\mathcal{A}_1}}  M_{1/b} \cdot M_{\nabla_x\cdot {\rm Im} {\bf r_2}} \partial_\tau u \dd\tau.
\end{align*}
Note that we have used the condition (2). Recalling that each $F_k$ satisfies $\partial_t F_k + \mathcal{P}_{\mathcal{A}_1}F_k = G_k$, we have from the energy calculation and \eqref{proof.thm.poisson.regular.off.1},
\begin{align}
& ~~~ \frac{\dd }{\dd t} \| M_{\sqrt{b}} F_k (t) \|_{L^2 (\R^d)}^2 \nonumber \\
& = - 2 {\rm Re} \langle \Lambda_{\mathcal{A}_1} F_k (t), F_k (t) \rangle _{L^2(\R^d)} + {\rm Re} \langle M_{\nabla_x\cdot {\bf r_2}} F_k (t), F_k (t)\rangle _{L^2(\R^d)}  +  2 {\rm Re} \langle M_b G_k (t), F_k (t) \rangle _{L^2(\R^d)}\nonumber \\
& \leq - 2\nu_1  {\rm Re} \langle \Lambda_{\mathcal{A}_1} F_k (t), F_k (t) \rangle _{L^2(\R^d)}  + \delta  \| F_k (t) \|_{\dot{H}^\frac{1+\alpha}{2} (\R^d)} \| F_k (t) \|_{\dot{H}^\frac{1-\alpha}{2} (\R^d)} \nonumber \\
& ~~~~ + C_5 \nu_1^{-1}\| M_{\sqrt{b}} F_k (t) \|_{L^2 (\R^d)}^2 +  2 {\rm Re} \langle M_b G_k (t), F_k (t) \rangle _{L^2(\R^d)}.\nonumber 
\end{align}
Thus, for $0<T\leq \min\{t_0/2, C_5\nu_1^{-1}\}$ the Gronwall inequality yields
\begin{align}
& ~~~~ \| M_{\sqrt{b}} F_k (T) \|_{L^2 (\R^d)}^2 + 2\nu_1 \int_0^T \|  F_k (t) \|_{\dot{H}^\frac12 (\R^d)}^2 \dd t \nonumber \\
& \leq C \bigg ( \delta \int_0^T \| F_k (t) \|_{\dot{H}^\frac{1+\alpha}{2} (\R^d)} \| F_k (t) \|_{\dot{H}^\frac{1-\alpha}{2} (\R^d)} \dd t  +  \int_0^T | \langle M_b G_k (t), F_k (t) \rangle _{L^2(\R^d)} | \dd t \bigg ),\label{proof.thm.poisson.regular.off.7}
\end{align}
where $C$ is a numerical constant.  Since $t$ is less than or equal to $t_0/2$, the term $F_1$ and $G_1$ are of lower order and easy to estimate. Indeed, by the definition $e^{-t\mathcal{Q}_{\mathcal{A}_1}} = M_{1/b} ( e^{-t\mathcal{P}_{\mathcal{A}_1^*}} )^* M_b$ the duality argument using the results of Proposition \ref{prop.poisson.analytic} and \eqref{proof.thm.poisson.regular.off.3} imply 
\begin{align}
\| G_1 (t) \|_{L^2 (\R^d)} \leq C t_0^{-\frac12} \| u (\frac{t_0}{2}) \|_{\dot{H}^\frac12 (\R^d)},~~~~~~~0<t\leq \frac{t_0}{2}, \label{proof.thm.poisson.regular.off.5}
\end{align}
where $C$ depends only on $d$, $\nu_1$, and $\nu_2$. The details are omitted. Then \eqref{proof.thm.poisson.regular.off.5} combined with \eqref{proof.thm.poisson.regular.off.3} yields
\begin{align}
\| F_1 (t) \|_{\dot{H}^\beta (\R^d)} \leq C\int_0^t (t-s)^{-\beta} \| G_1 (s) \|_{L^2 (\R^d)} \dd s \leq C t^{1-\beta} t_0^{-\frac12} \| u (t)\|_{\dot{H}^\frac12 (\R^d)},~~~~~ 0<t\leq \frac{t_0}{2} \label{proof.thm.poisson.regular.off.6}
\end{align}
for $\beta\in [0,1)$. Here the constant $C$ depends only on $d$, $\beta$, $\nu_1$, and $\nu_2$. From \eqref{proof.thm.poisson.regular.off.5} and \eqref{proof.thm.poisson.regular.off.6} we see 
\begin{align}
 \| F_1 (t) \|_{\dot{H}^\frac{1+\alpha}{2} (\R^d)} \| F_1 (t) \|_{\dot{H}^\frac{1-\alpha}{2} (\R^d)} + | \langle M_b G_1 (t), F_1 (t) \rangle _{L^2(\R^d)} |  + \| F_1 (t) \|_{L^2 (\R^d)}^2 \leq \frac{C t}{t_0} \| u (t) \|_{\dot{H}^\frac12 (\R^d)}^2\label{proof.thm.poisson.regular.off.8}
\end{align}
for $0<t\leq t_0/2$, with $C$ depending only on $d$, $\alpha$, $\nu_1$, and $\nu_2$. On the other hand, the identity \eqref{inner.product} and the assumption \eqref{proof.thm.poisson.regular.off.1} lead to  
\begin{align}
| \langle M_b G_2 (t), F_2 (t) \rangle _{L^2(\R^d)} | & \leq  \int_t^{t_0} | \langle M_{\nabla_x\cdot {\rm Im} {\bf r_2}}  \partial_s u,  e^{-(s-t)\mathcal{P}_{\mathcal{A}_1^{*}}} F_2 (t) \rangle _{L^2(\R^d)} | \dd s\nonumber \\
& \leq   \int_t^{t_0} \bigg \{  \delta \big (  \| m  \partial_s  u (s) \|_{\dot{H}^\frac{1+\alpha}{2} (\R^d)}  +  \| m^{-1} e^{-(s-t)\mathcal{P}_{\mathcal{A}_1^{*}}} F_2 (t)   \|_{\dot{H}^\frac{1+\alpha}{2} (\R^d)} \big ) \nonumber \\
& ~~~~~~~~~~~~~~ \cdot \big ( \| m  \partial_s  u (s) \|_{\dot{H}^\frac{1-\alpha}{2} (\R^d)}  +  \| m^{-1} e^{-(s-t)\mathcal{P}_{\mathcal{A}_1^{*}}} F_2 (t)   \|_{\dot{H}^\frac{1-\alpha}{2}(\R^d)} \big )   \nonumber \\
& ~~~~~~~~~~~~~~~~~~  + C_5 \big (  \| m  \partial_s  u (s) \|_{L^2 (\R^d)}  +  \| m^{-1} e^{-(s-t)\mathcal{P}_{\mathcal{A}_1^{*}}} F_2 (t)   \|_{L^2 (\R^d)}\big ) ^2 \bigg \} \dd s \label{proof.thm.poisson.regular.off.9}
\end{align} 
for any $m>0$. From Proposition \ref{prop.poisson.analytic} we have 
\begin{align}
\| \partial_s u (s) \|_{\dot{H}^\frac{1+\alpha}{2} (\R^d)}  = \| e^{-\frac{s}{4}\mathcal{P}_{\mathcal{A}}} \mathcal{P}_{\mathcal{A}} e^{-\frac{3}{4} s \mathcal{P}_{\mathcal{A}}} f \|_{\dot{H}^\frac{1+\alpha}{2} (\R^d)} & \leq  s^{-\frac{\alpha}{2}} \|\mathcal{P}_{\mathcal{A}} e^{-\frac{3}{4} s \mathcal{P}_{\mathcal{A}}} f \|_{\dot{H}^\frac{1}{2} (\R^d)} \nonumber \\
& \leq C s^{-1-\frac{\alpha}{2}} \| u(\frac{s}{2} ) \|_{\dot{H}^\frac12 (\R^d)} \leq C s^{-1-\frac{\alpha}{2}} \| u (\frac{t}{2} ) \|_{\dot{H}^\frac12 (\R^d)},\nonumber 
\end{align}
\begin{align}
\| \partial_s u (s) \|_{\dot{H}^\frac{1-\alpha}{2} (\R^d)}  \leq \| \partial_s u (s) \|_{L^2 (\R^d)} ^\alpha \| \partial_s u (s) \|_{\dot{H}^\frac{1}{2} (\R^d)} ^{1-\alpha} \leq C s^{-1 + \frac{\alpha}{2}} \| u(\frac{s}{2} ) \|_{\dot{H}^\frac12 (\R^d)} \leq C s^{-1 + \frac{\alpha}{2}} \| u (\frac{t}{2} ) \|_{\dot{H}^\frac12 (\R^d)},\nonumber 
\end{align}
while the interpolation inequality based on \eqref{proof.thm.poisson.regular.off.3} and $\| e^{-t\mathcal{P}_{\mathcal{A}_1^*}} f \|_{\dot{H}^{1/2} (\R^d)}\leq C \| f \|_{\dot{H}^{1/2} (\R^d)}$ (by Proposition \ref{prop.poisson.analytic}) yields 
\begin{align}
(s-t)^\frac{1+3 \alpha}{4} \| e^{-(s-t)\mathcal{P}_{\mathcal{A}_1^*}} F_2 (t) \|_{\dot{H}^\frac{1+\alpha}{2} (\R^d)} +  
(s-t)^\frac{1 - \alpha}{4} \| e^{-(s-t)\mathcal{P}_{\mathcal{A}_1^*}} F_2 (t) \|_{\dot{H}^\frac{1-\alpha}{2} (\R^d)} \leq C \| F_2 (t) \|_{\dot{H}^\frac{1-\alpha}{4} (\R^d)}\nonumber 
\end{align}
for $0<t<s\leq t_0/2$. Then, taking $m=t^{1/8} s^{3/8}$ in \eqref{proof.thm.poisson.regular.off.9}, we get
\begin{align}
& ~~~ | \langle M_b G_2 (t), F_2 (t) \rangle _{L^2(\R^d)} | \nonumber \\
& \leq   \int_t^{t_0}   \delta \bigg ( m^2 s^{-2} \|  u (\frac{t}{2}) \|_{\dot{H}^\frac{1}{2} (\R^d)}^2  + s^{-1+\frac{\alpha}{2}} (s-t)^{-\frac{1+3\alpha}{4}} \| u(\frac{t}{2}) \|_{\dot{H}^\frac12 (\R^d)} \| F_2 (t) \|_{\dot{H}^\frac{1-\alpha}{4} (\R^d)} \nonumber \\
& ~~~  + s^{-1-\frac{\alpha}{2}} (s-t)^{-\frac{1-\alpha}{4}} \| u(\frac{t}{2}) \|_{\dot{H}^\frac12 (\R^d)} \| F_2 (t) \|_{\dot{H}^\frac{1-\alpha}{4} (\R^d)} + m^{-2} (s-t)^{-\frac{1+\alpha}{2}} \| F_2 (t)   \|_{\dot{H}^\frac{1-\alpha}{4} (\R^d)} \bigg ) \dd s  \nonumber \\
& ~~~~~~+ C_5 \int_t^{t_0}  \big (  m^2 \| \partial_s  u (s) \|_{L^2 (\R^d)}^2  + m^{-2}  \|  e^{-(s-t)\mathcal{P}_{\mathcal{A}_1^{*}}} F_2 (t)   \|_{L^2 (\R^d)} ^2 \big ) \dd s \nonumber \\
& \leq C \bigg (  \delta \big ( \| u(\frac{t}{2} ) \|_{\dot{H}^\frac12 (\R^d)}^2 + t^{-\frac{1+\alpha}{2}} \| F_2 (t) \|_{\dot{H}^\frac{1-\alpha}{4}(\R^d)}^2\big ) +  C_5 \big ( t_0 \| u (t) \|_{\dot{H}^\frac12 (\R^d)}^2 + t^{-\frac14} t_0^\frac14 \| F_2 (t) \|_{L^2 (\R^d)}^2 \big ) \bigg ).\label{proof.thm.poisson.regular.off.10}
\end{align}
The triangle inequality and \eqref{proof.thm.poisson.regular.off.6} imply 
\begin{align}
\| F_2 (t) \|_{\dot{H}^\beta(\R^d)} & \leq \| u(t) \|_{\dot{H}^\beta(\R^d)} + \| e^{-t\mathcal{P}_{\mathcal{A}_1}} f \|_{\dot{H}^\beta(\R^d)} + \| F_1 (t) \|_{\dot{H}^\beta(\R^d)} \nonumber \\
& \leq \| u(t) \|_{\dot{H}^\beta(\R^d)} + \| e^{-t\mathcal{P}_{\mathcal{A}_1}} f \|_{\dot{H}^\beta(\R^d)} + C t^{1-\beta} t_0^{-\frac12} \| u (t) \|_{\dot{H}^\frac12 (\R^d)}\label{proof.thm.poisson.regular.off.11}
\end{align}
for $0<t\leq t_0/2$ and $\beta\in [0,1]$. Hence, if in addition $t_0$ is small so that $(C_5 + 1 ) t_0\leq \delta$, we get 
\begin{align}
| \langle M_b G_2 (t), F_2 (t) \rangle _{L^2(\R^d)} | & \leq C \bigg (  \delta \big ( \| u(\frac{t}{2} ) \|_{\dot{H}^\frac12 (\R^d)}^2 + t^{-\frac{1+\alpha}{2}} \| u(t) \|_{\dot{H}^\frac{1-\alpha}{4}(\R^d)}^2 + t^{-\frac{1+\alpha}{2}} \| F_2 (t) \|_{\dot{H}^\frac{1-\alpha}{4}(\R^d)}^2\big ) \nonumber \\
& ~~~~~ +  C_5  t^{-\frac14} t_0^\frac14  \big ( \| u(t) \|_{L^2 (\R^d)}^2 + \| f \|_{L^2 (\R^d)}^2 \big ) \bigg ).\label{proof.thm.poisson.regular.off.12}
\end{align}
Then we have from \eqref{est.prop.poisson.regular.1} with $\mathcal{A}$ replaced by $\mathcal{A}_1$, \eqref{proof.thm.poisson.regular.off.7}, \eqref{proof.thm.poisson.regular.off.8}, and \eqref{proof.thm.poisson.regular.off.12},
\begin{align}
& ~~~~ \| M_{\sqrt{b}} u (T) \|_{L^2 (\R^d)}^2 + 2\nu_1 \int_0^T \|  u (t) \|_{\dot{H}^\frac12 (\R^d)}^2 \dd t \nonumber \\
& \leq C \bigg ( \| f \|_{L^2 (\R^d)}^2 + (\delta + T t_0^{-1} ) \int_0^T \| u (t) \|_{\dot{H}^\frac12 (\R^d)}^2 \dd t + \delta \int_0^T \| F_2(t) \|_{\dot{H}^\frac{1+\alpha}{2} (\R^d)} \| F_2 (t) \|_{\dot{H}^\frac{1-\alpha}{2} (\R^d)} \dd t\nonumber \\
& ~~~~~   +  \delta \int_0^T  t^{-\frac{1+\alpha}{2}} \big ( \| u(t) \|_{\dot{H}^\frac{1-\alpha}{4}(\R^d)}^2 + \| e^{-t\mathcal{P}_{\mathcal{A}_1}} f \|_{\dot{H}^\frac{1-\alpha}{4} (\R^d)}^2 \big ) \dd t \bigg ) \nonumber \\
& ~~~~~   +  C C_5 t_0^\frac14 \int_0^T  t^{-\frac14}  \big ( \| u(t) \|_{L^2 (\R^d)}^2 + \| f \|_{L^2 (\R^d)}^2 \big ) \dd t. \label{proof.thm.poisson.regular.off.13}
\end{align}
Here $C$ depends only on $d$. $\alpha$, $\nu_1$, and $\nu_2$. Thus, by taking $T_0$ small enough depending on $\delta$ and $C_5$, and by using \eqref{est.prop.key.estimate.poisson.analytic} and \eqref{est.prop.poisson.regular.1} for $\{e^{-t\mathcal{P}_{\mathcal{A}_1}}\}_{t\geq 0}$ together with \eqref{proof.thm.poisson.regular.off.3}, we have for $0<T\leq T_0$,
\begin{align}
{\rm R.H.S.~of ~\eqref{proof.thm.poisson.regular.off.13}} & \leq C \bigg ( \| f \|_{L^2 (\R^d)}^2 + \delta \int_0^T \big ( \| u (t) \|_{\dot{H}^\frac12 (\R^d)}^2  + \| F_2 (t) \|_{\dot{H}^\frac{1+\alpha}{2} (\R^d)} \| F_2 (t) \|_{\dot{H}^\frac{1-\alpha}{2} (\R^d)} \big ) \dd t\nonumber \\
& ~~~~~   +  \delta \int_0^T  t^{-\frac{1+\alpha}{2}}  \| u(t) \|_{\dot{H}^\frac{1-\alpha}{4}(\R^d)}^2 \dd t  +  C_5 t_0^\frac14 \int_0^T  t^{-\frac14}  \| u(t) \|_{L^2 (\R^d)}^2  \dd t \bigg ). \label{proof.thm.poisson.regular.off.14}
\end{align}
Now we use 
\[
\| F_2 (t) \|_{\dot{H}^\frac{1+\alpha}{2} (\R^d)} \| F_2 (t) \|_{\dot{H}^\frac{1-\alpha}{2} (\R^d)}\leq \frac{t^\alpha}{2} \| F_2 (t) \|_{\dot{H}^\frac{1+\alpha}{2} (\R^d)} + \frac{t^{-\alpha}}{2} \| F_2 (t) \|_{\dot{H}^\frac{1-\alpha}{2} (\R^d)}
\]
and \eqref{proof.thm.poisson.regular.off.11}.  Then, again from \eqref{est.prop.key.estimate.poisson.analytic},  \eqref{est.prop.poisson.regular.1}, and \eqref{proof.thm.poisson.regular.off.3} for $\{e^{-t\mathcal{P}_{\mathcal{A}_1}}\}_{t\geq 0}$, and from \eqref{proof.thm.poisson.regular.off.6} for $F_1$, we observe that 
\begin{align}
&~~~~ {\rm R.H.S.~of ~\eqref{proof.thm.poisson.regular.off.14}} \nonumber \\
& \leq C \bigg ( \| f \|_{L^2 (\R^d)}^2 + \delta \int_0^T \big ( \| u (t) \|_{\dot{H}^\frac12 (\R^d)}^2  +  t^{-\alpha} \| u (t) \|_{\dot{H}^\frac{1-\alpha}{2} (\R^d)}^2 +  t^{-\frac{1+\alpha}{2}} \| u (t) \|_{\dot{H}^\frac{1-\alpha}{4} (\R^d)}^2  \big ) \dd t\nonumber \\
& ~~~~~   +   C_5 t_0^\frac14 \int_0^T  t^{-\frac14}  \| u(t) \|_{L^2 (\R^d)}^2  \dd t \bigg )\nonumber \\
& \leq C \bigg ( \| f \|_{L^2 (\R^d)}^2 + \delta \int_0^T \big ( \| u (t) \|_{\dot{H}^\frac12 (\R^d)}^2  + t^{-\frac{1+\alpha}{2}} \| u (t) \|_{\dot{H}^\frac{1-\alpha}{4} (\R^d)}^2  \big ) \dd t + C_5 t_0^\frac14 \int_0^T  t^{-\frac14}  \| u(t) \|_{L^2 (\R^d)}^2  \dd t \bigg ). \label{proof.thm.poisson.regular.off.15}
\end{align}
Here we have used $\| u(t) \|_{\dot{H}^{(1+\alpha)/2} (\R^d)} \leq C t^{-\alpha/2} \| u (t/2 ) \|_{\dot{H}^{1/2}(\R^d)}$ and the interpolation inequality for the term $t^{-\alpha} \| u(t) \|_{\dot{H}^{(1-\alpha)/2}(\R^d)}^2$.  By applying Proposition \ref{prop.key.estimate.poisson.analytic} for $\int_0^T t^{-(1+\alpha)/2} \| u (t) \|_{\dot{H}^{(1-\alpha)/4} (\R^d)}^2\dd t$ we get 
\begin{align}
{\rm R.H.S.~of ~\eqref{proof.thm.poisson.regular.off.15}} & \leq C \bigg ( \| f \|_{L^2 (\R^d)}^2 + \delta T^\frac{1-\alpha}{2} \| u(T) \|_{L^2 (\R^d)}^{1+\alpha} \| u(T) \|_{\dot{H}^\frac12 (\R^d)}^{1-\alpha} \nonumber \\
& ~~~~~ + \delta \int_0^T  \| u (t) \|_{\dot{H}^\frac12 (\R^d)}^2 \dd t + C_5 t_0^\frac14 \int_0^T  t^{-\frac14}  \| u(t) \|_{L^2 (\R^d)}^2  \dd t \bigg ), \label{proof.thm.poisson.regular.off.16}
\end{align}
where $C$ depends only on $d$, $\alpha$, $\nu_1$, and $\nu_2$. We take $\delta$ so small that $C\delta\leq \nu_1$. Then the Gronwall inequality and the Young inequality  give
\begin{align}
\| M_{\sqrt{b}} u(T) \|_{L^2 (\R^d)}^2 + \nu_1 \int_0^T \| u(t) \|_{\dot{H}^\frac12 (\R^d)}^2 \dd t \leq C \big ( \| f \|_{L^2 (\R^d)}^2 + \delta   \| u(T) \|_{L^2 (\R^d)}^2 + \delta T^\frac12 \| u(T) \|_{\dot{H}^\frac12 (\R^d)}^2.\label{proof.thm.poisson.regular.off.17}
\end{align}
By taking $\delta$ further small (but depending only in $d$, $\alpha$, $\nu_1$, and $\nu_2$) if necessary, the argument as in the derivation of \eqref{proof.prop.poisson.regular.5} leads to the desired estimate:
\begin{align}
\| M_{\sqrt{b}} u(T) \|_{L^2 (\R^d)}^2 + \nu_1 \int_0^T \| u(t) \|_{\dot{H}^\frac12 (\R^d)}^2 \dd t \leq C  \| f \|_{L^2 (\R^d)}^2,  ~~~~~~~ 0<T\leq T_0\label{proof.thm.poisson.regular.off.17}
\end{align}
for sufficiently small $T_0>0$. Here $C$ depends only on $d$, $\alpha$, $\nu_1$, and $\nu_2$, while  $T_0$ depends only on $d$, $\alpha$, $\delta$, $C_5$, $\nu_1$, and $\nu_2$. Clearly the same estimate is valid also for $v(t) = e^{-t\mathcal{P}_{\mathcal{A}^*}} f$. Hence Proposition \ref{thm.domain.embedding} implies 
\begin{align}
\| \nabla_x f \|_{L^2 (\R^d)} \leq C \big ( \min\big \{ \| \mathcal{P}_{\mathcal{A}} f \|_{L^2 (\R^d)}, \| \mathcal{P}_{\mathcal{A}^*}f \|_{L^2 (\R^d)}\big \}  + \| f \|_{L^2 (\R^d)} \big ), ~~~~~~~~f\in H^1 (\R^d),\label{proof.thm.poisson.regular.off.18}
\end{align}
where $C$ depends only on $d$, $\alpha$, $\delta$, $C_5$, $\nu_1$, and $\nu_2$. 

Now let us eliminate the assumption $a_{i,j}\in {\rm Lip}(\R^d)$. First we consider the mollified matrix $A^{(\epsilon)}=(a_{i,j}^{(\epsilon)})_{1\leq i,j\leq d+1}$, $a_{i,j}^{(\epsilon)}=j_\epsilon * a_{i,j}$, where $j_\epsilon$ is a standard mollifier. Note that the conditions (1) (or \eqref{proof.thm.poisson.regular.off.1}) and (2) are invariant under this mollification. We have already known that $H^1 (\R^d)= D_{L^2}(\mathcal{P}_{\mathcal{A}^{(\epsilon)}})$ by \cite[Theorem 1.2]{MM1}. Moreover, \eqref{proof.thm.poisson.regular.off.18} and the proof of Corollary \ref{cor.prop.embedding.section3.1} implies that $\| f \|_{H^1 (\R^d)} \simeq \| f \|_{D_{L^2}(\mathcal{P}_{\mathcal{A}^{(\epsilon)}})}$ with constants independent of the parameter $\epsilon>0$.  Let $f\in H^1 (\R^d)$ and set $u_\epsilon (t) = e^{-t\mathcal{P}_{\mathcal{A}^{(\epsilon)}}}f $. Then we have the uniform bounds such that
\begin{align*}
&\| \nabla u_\epsilon \|_{L^2 (\R^{d+1}_+)} \leq C \| f \|_{\dot{H}^\frac12 (\R^d)},~~~~~\sup_{0<t<T_0} \| \nabla u_\epsilon (t) \|_{L^2 (\R^d)} \leq C \| f\|_{H^1 (\R^d)},\\
&\sup_{0<t<T_0} \| u_\epsilon (t)  \|_{L^2 (\R^d)}^2  +  \int_0^{T_0} \| u_\epsilon (t)   \|_{\dot{H}^\frac12 (\R^d)}^2 \dd t \leq C \| f\|_{L^2 (\R^d)}^2.
\end{align*}
By the standard limiting procedure and the uniqueness of weak solutions we conclude that $w(t) = e^{-t\mathcal{P}_{\mathcal{A}}}f$ satisfies 
\begin{align}
\| \nabla w \|_{L^2 (\R^{d+1}_+)} & \leq C \| f \|_{\dot{H}^\frac12 (\R^d)},\\
\sup_{0<t<T_0} \| \nabla w (t) \|_{L^2 (\R^d)} & \leq C \| f\|_{H^1 (\R^d)},\label{proof.thm.poisson.regular.off.19}\\
\sup_{0<t<T_0} \| w (t)  \|_{L^2 (\R^d)}^2  +  \int_0^{T_0} \| w (t)   \|_{\dot{H}^\frac12 (\R^d)}^2 \dd t & \leq C \| f\|_{L^2 (\R^d)}^2.\label{proof.thm.poisson.regular.off.20}
\end{align}
The estimate \eqref{proof.thm.poisson.regular.off.20} implies that $\{e^{-t\mathcal{P}_{\mathcal{A}}}\}_{t\geq 0}$ is extended as a strongly continuous semigroup acting on $L^2 (\R^d)$, and by  Proposition \ref{thm.domain.embedding} the space $D_{L^2}(\mathcal{P}_{\mathcal{A}})$ is continuously embedded in $H^1 (\R^d)$. The same is true for  $\{e^{-t\mathcal{P}_{\mathcal{A}^*}}\}_{t\geq 0}$. Hence, by \eqref{proof.thm.poisson.regular.off.19} and Proposition \ref{prop.embedding.section3.2} (ii) we have $H^1 (\R^d) = D_{L^2}(\mathcal{P}_{\mathcal{A}})$ with equivalent norms. Similarly, $H^1 (\R^d) = D_{L^2}(\mathcal{P}_{\mathcal{A}^*})$ holds  with equivalent norms. The proof is complete.


%
%
%
%

\subsection{Proof of Theorem \ref{thm.factorization.strong.1}} 
\label{subsec.solvability}

We are now in position to prove 
Theorem \ref{thm.factorization.strong.1}. 
From Theorems \ref{thm.poisson.hermite} and \ref{thm.poisson.regular.off} we have $H^1 (\R^d)\hookrightarrow D_{L^2}(\Lambda_{\mathcal{A}})\cap D_{L^2}(\Lambda_{\mathcal{A}^*})$,  $\mathcal{P}_{\mathcal{A}}={\bf P}_{\mathcal{A}}$, and $\mathcal{P}_{\mathcal{A}^*} = {\bf P}_{\mathcal{A}^*}$ under the assumptions of Theorem \ref{thm.factorization.strong.1}. The factorizations \eqref{eq.factorization.strong.0} and \eqref{eq.factorization.strong.1} follow from Lemma \ref{lem.factorization.strong}.  The proof is complete.



\section{Application to boundary value problem}\label{sec.application}

In this section we consider the Dirichlet and Neumann
problems for elliptic operators in $\R^{d+1}_{+}$.

\subsection{Relation between mild solution and weak solution}
\label{subsec.mild.weak}

The aim of this section is to clarify the relation between 
the weak solution and the mild solution 
for the boundary value problems. 
\begin{thm}\label{thm.representation} Assume that the semigroups $\{e^{-t\mathcal{P}_{\mathcal{A}}}\}_{t\geq 0}$ and  $\{e^{-t\mathcal{P}_{\mathcal{A}^*}}\}_{t\geq 0}$ acting on $H^{1/2}(\R^d)$ are extended as  strongly continuous semigroups acting on $L^2 (\R^d)$. Assume further that  $D_{L^2}(\mathcal{P}_{\mathcal{A}}) = D_{L^2}(\mathcal{P}_{\mathcal{A}^*}) = H^1 (\R^d)$ with equivalent norms. Let $F\in \dot{H}^{-1}(\R^{d+1}_+), ~g\in H^{1/2}(\R^d)$, and  $u, v\in C([0,\infty); L^2 (\R^d)) \cap \dot{H}^1 (\R^{d+1}_+)$. Then the following statements hold.

\noindent {\rm (i)} If $u$ is a weak solution to \eqref{eq.dirichlet}  then  $u$ is a mild solution to  \eqref{eq.dirichlet}.

\noindent {\rm (ii)} If in addition $\partial_t F \in\dot{H}^{-1}(\R^{d+1}_+)$ and  $h =g + M_{b}  \int_0^\infty   e^{-s \mathcal{Q}_{\mathcal{A}}} M_{1/b} F (s) \dd s$ belong to the range of $\Lambda_{\mathcal{A}}$ and if $v$ is a weak solution to \eqref{eq.neumann}  then $v$ is a mild solution to \eqref{eq.neumann}.

\end{thm}

\begin{rem}{\rm In Section \ref{subsec.def.weak.sol}, we define the mild solution for $F \in L^1_{loc}(\R_+;L^2(\R^d))$. However, if the semigroup $\{e^{-t\mathcal{P_{A^*}}}\}_{t\ge0}$ is analytic in $L^2(\R^d)$ and $D(\mathcal{P_{A^*}})=H^1(\R^d)$ holds, the Duhamel terms in the integral equations still make sense even for $F\in \dot{H}^{-1}(\R^{d+1}_+)$ by the duality argument. We note that the Duhamel terms make sense also for $F$ in some Sobolev spaces of negative order in the same way.


}
\end{rem}

\noindent {\it Proof of Theorem \ref{thm.representation}.} Firstly we observe that  
\begin{align}
m_{\mathcal{A}} [F](t) = \int_t^\infty e^{-(\tau-t)\mathcal{Q}_{\mathcal{A}}} M_{1/b} F (\tau ) \dd \tau \in L^2_{loc} (\R_+; L^2 (\R^d)) ~~~~~{\rm for} ~~ F \in \dot{H}^{-1} (\R^{d+1}_+).\label{proof.thm.representation.0}
\end{align}
Indeed, we have for $\varphi\in C_0^\infty (\R^{d+1}_+)$ and $F \in L^2 (\R^{d+1}_+)\cap\dot{H}^{-1} (\R^{d+1}_+)$,
\begin{align}
 \langle m_{\mathcal{A}}[F], \varphi \rangle _{L^2 (\R^{d+1}_+)} & = \int_0^\infty \int_t^\infty \langle F (\tau), e^{-(\tau-t)\mathcal{P}_{\mathcal{A}^*}} M_{1/\bar{b}} \varphi (t) \rangle _{L^2 (\R^d)} \dd \tau \dd t\nonumber \\
& = \int_0^\infty \langle F (\tau), \int_0^\tau e^{-(\tau-t)\mathcal{P}_{\mathcal{A}^*}} M_{1/\bar{b}} \varphi (t) \dd t \rangle _{L^2 (\R^d)} \dd\tau = \langle F, m_{\mathcal{A}}^*[\varphi] \rangle _{L^2 (\R^{d+1}_+)}.\label{proof.thm.representation.1}
\end{align}
Here $m_{\mathcal{A}}^*[\varphi] (\tau ) =  \int_0^\tau e^{-(\tau-t)\mathcal{P}_{\mathcal{A}^*}} M_{1/\bar{b}} \varphi (t) \dd t$. The right-hand side of \eqref{proof.thm.representation.1} can be replaced by $\langle F, m_{\mathcal{A}}^*[\varphi] \rangle _{\dot{H}^{-1},\dot{H}^1}$, and thus, we have 
$| \langle m_{\mathcal{A}}[F], \varphi \rangle _{L^2 (\R^{d+1}_+)} |\leq \| F\|_{\dot{H}^{-1}(\R^{d+1}_+)} \| \nabla m_{\mathcal{A}}^*[\varphi] \|_{L^2 (\R^{d+1}_+)}$. From $D_{L^2}(\mathcal{P}_{\mathcal{A}^*}) = H^1(\R^d)$ and the maximal regularity,  $ \| \nabla m_{\mathcal{A}}^*[\varphi] \|_{L^2 (0,T; L^2 (\R^d))}\leq C \| \varphi \|_{L^2(0,T; L^2 (\R^d))}$ holds for all $T>0$. On the other hand, since $\varphi(t)=0$ if $t\geq T_0/2$ for some $T_0>1$, when $t\geq T_0$ we can write $m_{\mathcal{A}}^*[\varphi](t) = e^{-(t-T_0)\mathcal{P}_{\mathcal{A}^*}}h$ for some $h \in D_{L^2}(\mathcal{P}_{\mathcal{A}^*})$ satisfying $\| h \|_{H^1 (\R^d)}\leq C T_0^{1/2} \| \varphi\|_{L^2 (0,T_0; L^2 (\R^d))}$. Thus we have from \eqref{est.proof.prop.poisson.analytic.1},
\[
\| \nabla m_{\mathcal{A}}^*[\varphi]\|_{L^2(T_0,\infty; L^2 (\R^d))} \leq C\| h\|_{\dot{H}^{1/2}(\R^d)} \leq C T_0^{1/2} \| \varphi\|_{L^2 (0,T_0; L^2 (\R^d))}. 
\]
Collecting these, we conclude that $m_{\mathcal{A}}[F]\in L^2_{loc} (\R_+; L^2 (\R^d))$. Note that, if $\partial_t F\in \dot{H}^{-1}(\R^{d+1}_+)$ in addition, then $\partial_t m_{\mathcal{A}}[F]\in L^2_{loc}(\R_+; L^2 (\R^d))$ by the same calculation as above.

Now let us prove the assertion (i).  Note that $u_0 (t) = e^{-t\mathcal{P}_{\mathcal{A}}}g$, $g\in H^{1/2} (\R^d)$, is a weak solution to \eqref{eq.dirichlet} with $F=0$ belonging to $C(\overline{\R_+}; L^2 (\R^d))\cap \dot{H}^1(\R^{d+1}_+)$ due to Proposition \ref{prop.poisson.analytic.L^2}. Thus, considering $u(t)-u_0(t)$ if necessary, we may assume that $u$ is a weak solution to  \eqref{eq.dirichlet} with $g=0$ belonging to $C(\overline{\R_+}; L^2 (\R^d)) \cap \dot{H}^1 (\R^{d+1}_+)$.  From $F\in \dot{H}^{-1}(\R^{d+1}_+)$ the uniqueness of weak solutions in $C(\overline{\R_+}; L^2 (\R^d))\cap \dot{H}^1(\R^{d+1}_+)$ implies  that $u\in \dot{H}^1_0 (\R^{d+1}_+)$ and  $\langle A\nabla  u ,\nabla  \varphi \rangle _{L^2 (\R^{d+1}_+)}= \langle F , \varphi \rangle _{\dot{H}^{-1},\dot{H}^1}$ for all $\varphi \in \dot{H}^1_0 (\R^{d+1}_+)$. For a given $\psi\in C_0^\infty (\R^{d+1}_+)$ and $0<\epsilon\ll 1$ set $\varphi (t) = \int_0^t e^{-(t-s)\mathcal{P}_{\mathcal{A}^*}} e^{-\epsilon\mathcal{P}_{\mathcal{A}^*}} m_{\mathcal{A}^*} [\psi] (s) \dd s$, where $m_{\mathcal{A}^*} [\psi] (s) = \int_s^\infty e^{-(\tau-s)\mathcal{Q}_{\mathcal{A}^*}} M_{1/\bar{b}} \psi (\tau ) \dd \tau\in L^2_{loc} (\R_+; L^2 (\R^d) )$ due to \eqref{proof.thm.representation.0}. Then we have   $\partial_t \varphi,\mathcal{P}_{\mathcal{A}^*}\varphi \in L^2_{loc} (\R_+; H^{1/2} (\R^{d}))$  by the maximal regularity, for $e^{-\epsilon\mathcal{P}_{\mathcal{A}^*}} m_{\mathcal{A}^*} [\psi]$ belongs to the space $L^2_{loc} (\R_+; H^{1/2} (\R^{d}))$. Moreover, since $m_{\mathcal{A}^*} [\psi](s)=0$ for $s\geq T_0/2 \gg 1$ the function $\varphi$ is written in the form $\varphi (t) = e^{-(t- T_0+\epsilon)\mathcal{P}_{\mathcal{A}^*}} w$ with some $w\in D_{H^{1/2}}(\mathcal{P}_{\mathcal{A}^*})$ if $t\geq  T_0$. Hence, recalling \eqref{est.proof.prop.poisson.analytic.1'}, we conclude that $\varphi \in \dot{H}^1_0(\R^{d+1}_+)$ and $\partial_t\varphi, \mathcal{P}_{\mathcal{A}^*} \varphi \in L^2 (\R_+; H^{1/2} (\R^{d}))$. Then, Theorem \ref{thm.factorization} yields 
\begin{align*}
& \langle A \nabla u, \nabla \varphi\rangle _{L^2(\R^{d+1}_+)}  = \langle \nabla u, A^*\nabla \varphi\rangle _{L^2(\R^{d+1}_+)} \\
& ~~~~~~~~~~~~~ = \int_0^\infty \langle \Lambda_{\mathcal{A}} u, (\partial_t + \mathcal{P}_{\mathcal{A}^*} ) \varphi\rangle_{\dot{H}^{-\frac12},\dot{H}^\frac12} \dd t + \langle ( M_b \partial_t + M_{\bf r_2}\cdot \nabla_x ) u, (\partial_t + \mathcal{P}_{\mathcal{A}^*} ) \varphi\rangle_{L^2 (\R^{d+1}_+)}.
\end{align*}
Since $u(t) \in H^1 (\R^d)\hookrightarrow D_{L^2}(\Lambda_{\mathcal{A}})$ for a.e. $t>0$ and  $M_ b \mathcal{P}_{\mathcal{A}} u =  \Lambda_{\mathcal{A}} u + M_{{\bf r_2}}\cdot \nabla_x u$ by Corollary \ref{cor.prop.embedding.section3.1'}, we have $\langle A\nabla u , \nabla \varphi \rangle _{L^2 (\R^{d+1}_+)} =  \langle ( \partial_t  + \mathcal{P}_{\mathcal{A}} ) u, M_{\bar{b}} (\partial_t + \mathcal{P}_{\mathcal{A}^*} )\varphi \rangle _{L^2(\R^{d+1}_+)}$, that is,
\begin{align}
\int_0^\infty \langle (\partial_t  + \mathcal{P}_{\mathcal{A}} ) u, M_{\bar{b}} (\partial_t + \mathcal{P}_{\mathcal{A}^*} )\varphi \rangle _{L^2(\R^d)} \dd t= \langle F  , \varphi \rangle _{\dot{H}^{-1}, \dot{H}^1}. \label{proof.thm.representation.2}
\end{align}
By the definition of $\varphi$ and the integration by parts, the left-hand side of \eqref{proof.thm.representation.2} is calculated as 
\begin{align}
{\rm L.H.S.~of}~\eqref{proof.thm.representation.2} & =\int_0^\infty \langle (\partial_t  + \mathcal{P}_{\mathcal{A}} ) u, M_{\bar{b}} e^{-\epsilon\mathcal{P}_{\mathcal{A}^*}} m_{\mathcal{A}^*}[\psi] \rangle _{L^2(\R^d)} \dd t \nonumber \\
& \rightarrow \int_0^\infty \langle \partial_t u + \mathcal{P}_{\mathcal{A}} u, M_{\bar{b}}  m_{\mathcal{A}^*}[\psi]  \rangle _{L^2(\R^d)} \dd t = \int_0^\infty  \langle u, \psi\rangle _{L^2(\R^d)} \dd t = \langle u, \psi \rangle _{L^2 (\R^{d+1}_+)}.\label{proof.thm.representation.3}
\end{align} 
On the other hand, by using \eqref{proof.thm.representation.1} the right-hand side of  \eqref{proof.thm.representation.2} satisfies
\begin{align}
{\rm R.H.S.~of}~\eqref{proof.thm.representation.1} \rightarrow \langle F, m_{\mathcal{A}}^*\big [M_{\bar{b}} m_{\mathcal{A}^*}[\psi ] \big ]\rangle _{\dot{H}^{-1},\dot{H}^1} & =  \langle m_{\mathcal{A}}[F],   M_{\bar{b}} m_{\mathcal{A}^*}[\psi ] \rangle _{L^2 (\R^{d+1}_+)} \nonumber \\
& = \langle m_{\mathcal{A}^*}^*\big [ M_b m_{\mathcal{A}}[F] \big ],  \psi \rangle _{L^2 (\R^{d+1}_+)}.\label{proof.thm.representation.4}
\end{align}
Hence \eqref{proof.thm.representation.2}-\eqref{proof.thm.representation.4} yield $u(t) =   m_{\mathcal{A}^*}^*\big [ M_b m_{\mathcal{A}}[F] \big ] (t) = \int_0^t e^{-(t-s)\mathcal{P}_{\mathcal{A}}} m_{\mathcal{A}}[F] (s) \dd s$, which proves (i).

Next we show (ii). Let $w\in \dot{H}^1_0 (\R^{d+1}_+)$ be the weak solution to \eqref{eq.dirichlet} with $g=0$ (which uniquely exists by the assumption $F \in \dot{H}^{-1}(\R^{d+1}_+)$). Then  $w$ is a mild solution because of  (i), and since $\partial_t m_{\mathcal{A}}[F]\in L^2_{loc}(\R_+; L^2 (\R^d))$ for $\partial_t F\in \dot{H}^{-1} (\R^{d+1}_+)$, we have $\partial_t w (t) = - \mathcal{P}_{\mathcal{A}} w (t) + m_{\mathcal{A}}[F] (t)$ for all $t>0$. Applying Corollary \ref{cor.prop.embedding.section3.1'}, we see 
\begin{align*}
-M_b \partial_t w (t) - M_{{\bf r_1}}\cdot \nabla_x w (t) = \Lambda_{\mathcal{A}} w (t) - M_b m_{\mathcal{A}}[F] (t),~~~~~~~~~t>0.
\end{align*}
Hence, for all $\varphi \in H^1(\R^{d+1}_+)$ it follows that
\begin{align*} 
\langle A\nabla w, \nabla \varphi \rangle _{L^2 (\R^{d+1}_+)} & = \lim_{t\rightarrow 0} \int_t^\infty \langle A\nabla w, \nabla \varphi \rangle _{L^2 (\R^d)} \dd s \nonumber \\
& = \lim_{t\rightarrow 0} \langle -M_b \partial_t w (t) - M_{{\bf r_1}}\cdot \nabla_x w (t), \varphi (t) \rangle _{L^2 (\R^d)} + \langle F, \varphi \rangle _{\dot{H}^{-1},\dot{H}^1}\nonumber \\
& =  \lim_{t\rightarrow 0} \langle\Lambda_{\mathcal{A}} w (t) - M_b m_{\mathcal{A}}[F] (t), \varphi (t) \rangle _{L^2 (\R^d)} + \langle F, \varphi \rangle _{\dot{H}^{-1},\dot{H}^1}.
\end{align*}
Then from the estimate 
\[
|\langle\Lambda_{\mathcal{A}} w (t), \varphi (t) \rangle _{L^2 (\R^d)}|\leq C \| w(t) \|_{\dot{H}^{\frac12}(\R^d)} \|\varphi (t) \|_{\dot{H}^{\frac12}(\R^d)}\leq C t^\frac12 \| m_{\mathcal{A}}[F] \|_{L^\infty (0,t; L^2 (\R^d))} \| \varphi \|_{H^1 (\R^{d+1}_+)}
\]
for $t\in (0,1)$, we have 
\begin{align} 
\langle A\nabla w, \nabla \varphi \rangle _{L^2 (\R^{d+1}_+)} & =  - \langle   M_b m_{\mathcal{A}}[F] , \varphi  \rangle _{L^2 (\partial\R^{d+1}_+)} + \langle F, \varphi \rangle _{\dot{H}^{-1},\dot{H}^1}.\label{proof.thm.representation.5}
\end{align}
Next we set $v_0 (t) = e^{-t \mathcal{P}_{\mathcal{A}}} \Lambda_{\mathcal{A}}^{-1} h$, where $h= g + M_b m_{\mathcal{A}}[F](0)$. Then the definition of $\Lambda_{\mathcal{A}}$ implies
\begin{align}
\langle A\nabla v, \nabla \varphi \rangle _{L^2 (\R^{d+1}_+)} = \langle \Lambda_{\mathcal{A}}  \Lambda_{\mathcal{A}}^{-1} h, \varphi \rangle _{\dot{H}^{-\frac12}, \dot{H}^\frac12} = \langle h, \varphi \rangle _{L^2 (\partial\R^{d+1}_+)},~~~~~~~~\varphi \in H^1 (\R^{d+1}_+).\label{proof.thm.representation.6}
\end{align}
The equalities  \eqref{proof.thm.representation.5} and \eqref{proof.thm.representation.6} show that the function $\tilde v$ defined by $\tilde v = v_0+ w\in \dot{H}^1 (\R^{d+1}_+)$, which is a mild solution to \eqref{eq.neumann} by the construction, is a weak solution to \eqref{eq.neumann}.  Then $z=v-\tilde v$ belongs to $\dot{H}^1_0 (\R^{d+1}_+)$ and $\langle A\nabla z, \nabla \phi \rangle _{L^2 (\R^{d+1}_+)}=0$ for all $\phi\in C_0^\infty (\overline{\R^{d+1}_+})$. Hence the uniqueness of weak solutions to \eqref{eq.dirichlet} in  $\dot{H}^1_0 (\R^{d+1}_+)$ verifies the representation $z (t) = e^{-t\mathcal{P}_{\mathcal{A}}} z(0)$, i.e., $z\in \dot{H}^1(\R^{d+1}_+)$. Then we can take $\phi =z$ in the above equality, which gives $z=0$, i.e., $v=\tilde v$. The proof is complete.

\subsection{Proof of Corollary \ref{thm.solvability.intro1}}

In order to show Corollary \ref{thm.solvability.intro1}, 
we use the following elementary estimate:
\begin{lem}
\label{est.phi}
Let $u \in \dot{H}^1_0(\R^{d+1}_+)$ 
be the weak solution to \eqref{eq.dirichlet} with $F \in C_0^{\infty}(\R^{d+1}_+)$ and $g=0$. 
Then $\|\nabla u (t) \|_{L^2(\R^d)} \le C (1+t)^{3/2}$ holds,  
where $C$ depends only on $F$. 
\end{lem}
\noindent
{\it Proof.}~
By 
Theorem \ref{thm.representation} and a simple change of variables, 
we have the representation formula:
$$
u (t)=\int_0^t e^{-s\mathcal{P}_\mathcal{A}}
\int^{\infty}_0e^{-\tau \mathcal{Q_A}}M_{1/b} F(\tau+t-s)\dd \tau 
\dd s.
$$
Then since  Proposition \ref{prop.poisson.analytic.L^2} yields the estimate
$$\|e^{-t\mathcal{P_A}}\|_{L^2\rightarrow L^2}
+\|e^{-t\mathcal{Q_A}}\|_{L^2\rightarrow L^2}
\le
(1+t)^{1/2},$$
we have 
\begin{align*}
&~~~ \|\partial_t u (t)\|_{L^2(\R^d)}\\
& \le 
\|e^{-t\mathcal{P_A}}
\int^{\infty}_0
e^{-\tau \mathcal{Q_A}}F(\tau)\dd \tau
\|_{L^2(\R^d)}
+
\|\int^t_0
e^{-s\mathcal{P}_\mathcal{A}}
\int^{\infty}_0
e^{-\tau \mathcal{Q_A}}
\partial_t F(\tau+t-s)\dd \tau \dd s
\|_{L^2(\R^d)}
\\
&\le 
(1+t)^{\frac12}
\int^{\infty}_0(1+\tau)^{\frac12}
\|F(\tau)
\|_{L^2(\R^d)}\dd \tau
+
\int^t_0(1+s)^{\frac12}
\int^{\infty}_0
(1+\tau)^{\frac12}
\|
\partial_t F(\tau+t-s)\|_{L^2(\R^d)}
\dd \tau \dd s
\\
&\le C(1+t)^{\frac32},
\end{align*}
where we have used the fact that  
$F$ is smooth and has a compact support in the last inequality. 
Since $H^1(\R^d)=D_{L^2}(\mathcal{P_A})$ holds with the equivalent norm, 
the estimate for $\|\nabla_x u(t)\|_{L^2(\R^d)}$ 
can be proved in the same way. 
Thus the proof is complete.

\vspace{0.3cm}

\noindent
{\it Proof of Corollary \ref{thm.solvability.intro1}.}
Since the existence part is a direct consequence of Theorem 
\ref{thm.factorization.strong.1}, we will prove the uniqueness
of the solution in the class $C(\overline{\R_+};L^2(\R^{d}))$\ $\cap$\ $
H^1(\R^d \times (\delta,\infty))$ 
for any $\delta>0$. Let us 
consider the Dirichlet problem \eqref{eq.dirichlet}.
Let $u$, $\tilde u$ be two solutions to \eqref{eq.dirichlet} with the same boundary data 
and $F=0$, and let $w:=u-\tilde u$.
It suffices to show 
$\langle w, F \rangle _{L^2 (\R^{d+1}_+)}=0$
for all $F \in C^{\infty}_0(\R^{d+1}_+)$.
To this end, we first note that $w$ satisfies
\begin{align}
\label{eq.w}
\langle A\nabla w, \nabla (\phi \varphi) \rangle _{L^2 (\R^{d+1}_+)} 
=0~~~~~~~~{\rm for~all}~~
\varphi\in \dot{H}^1_0(\R^{d+1}_+), ~~\phi \in 
C_0^{\infty}(\R_+).  
\end{align}
For small $r \in (0,1/2)$,  
we take $\phi$ such that  $\phi(t)=1$ if 
$t \in (2r,1/r)$, 
$\phi(t)=0$ if $t \notin (r,2/r)$, 
$|\partial_t \phi(t)|\le C/r$ if $t \in (r,2r)$, 
and $|\partial_t \phi(t)| \le Cr$ if $t \in (1/r,2/r)$.
By a direct calculation we have 
\begin{align}
\label{eq.w2}
& ~~~\langle A\nabla w, \nabla (\phi \varphi) \rangle _{L^2 (\R^{d+1}_+)} \nonumber \\
&=
\langle A\nabla (\phi w), 
\nabla \varphi \rangle _{L^2 (\R^{d+1}_+)} 
-
\langle w \partial_t \phi ~{\bf a_{1, d+1}}, 
\nabla \varphi \rangle _{L^2 (\R^{d+1}_+)}  
+\langle {\bf a_{2, d+1}} \cdot \nabla w, \varphi\partial_t \phi \rangle _{L^2 (\R^{d+1}_+)} 
\notag
\\
&=
I+I\!I+I\!I\!I.
\end{align}
Here ${\bf  a_{1, d+1}} = ({\bf r_1},b)^\top$ and $ {\bf a_{2, d+1}} = ({\bf r_2},b)^\top$.  By the Lax-Milgram theorem, for all 
$F \in C_0^{\infty}(\R^{d+1}_+)$, there exists a unique 
solution $\varphi$ in the class $\dot{H}^1_0(\R^{d+1}_+)$ 
to the adjoint problem
$$
\langle A\nabla \psi, 
\nabla \varphi \rangle _{L^2 (\R^{d+1}_+)} 
=
\langle \psi, F \rangle _{L^2 (\R^{d+1}_+)} 
\quad 
{\rm for} \ \psi \in \dot{H}^1_0(\R^{d+1}_+).
$$
Now in \eqref{eq.w} taking $\varphi$  
as the solution of this problem, we have  $I=\langle \phi w, F \rangle _{L^2 (\R^{d+1}_+)}$.
On the other hand, the third term in \eqref{eq.w2} 
can be splitted as follows
$$
I\!I\!I=\int^{2r}_r \dd t +\int^{2/r}_{1/r} \dd t=I\!I\!I_a+I\!I\!I_b, 
$$
and then since the estimate 
$\|\nabla E_{\mathcal{A}} f \|_{L^2(r,\infty; L^2(\R^{d}))} 
\leq Cr^{-1}\| E_{\mathcal{A}} f \|_{L^2 (r/2,r; L^2(\R^{d}))}$
 in Remark \ref{rem.prop.poisson.analytic.1}
still holds for general weak solutions in $\dot{H}^1(\R^{d+1}_+)$, we see

\begin{align*}
I\!I\!I_a 
&\le 
 \|\nabla w(t)\|_{L^2(r,2r;L^2(\R^d))}
\sup_{r<t<2r} \| \varphi(t)\|_{L^2(\R^d)}
\|\partial_t \phi\|_{L^2(r,2r)} 
\\
&\le
Cr^{-1} \sup_{r/2<t<r}\| w(t)\|_{L^2(\R^d)}
\sup_{r<t<2r} \| \int^t_0 \partial_t \varphi(s)\dd s\|_{L^2(\R^d)}
\\
&\le 
Cr^{-1} \sup_{r/2<t<r}\| w(t)\|_{L^2(\R^d)}
\int^{2r}_0 (1+s)^{3/2}\dd s 
\\
&\le
C\sup_{r/2<t<r}\| w(t)\|_{L^2(\R^d)}
\end{align*}
for $0<r<1/2$, where we have used Lemma \ref{est.phi}
in the second to the last line. 
Since $w(t)$ converges to  $0$ in $L^2(\R^d)$ 
as $t\rightarrow 0$, $I\!I\!I_a$ converges to 0 as $r \rightarrow 0$.
Similary, we have 
\begin{align*}
I\!I\!I_b 
&\le 
\|\nabla w(t)\|_{L^2(1/r,2/r;L^2(\R^d))}
\sup_{1/r<t<2/r} \| \varphi(t)\|_{L^2(\R^d)}
(\int^{2/r}_{1/r} |\partial_t \phi(t)|^2\dd t )^{1/2}
\\
&\le 
Cr^{1/2}\|\nabla w(t)\|_{L^2(1/r,2/r;L^2(\R^d))}
 \int^{2/r}_0
\|\partial_t 
\varphi(s)\|_{L^2(\R^d)}\dd s
\\
&\le
C\|\nabla w(t)\|_{L^2(1/r,2/r;L^2(\R^d))}
 \|\partial_t \varphi\|_{L^2(\R_+;L^2(\R^d))},
\end{align*}
which converges to $0$ as $r \rightarrow 0$
by the assumption 
$w \in \dot{H}^1(\R^{d} \times (\delta,\infty))$ for any $\delta>0$.
In the same manner, we also see the term $I\!I$ converges to $0$ as 
$r \rightarrow 0$.
Thus 
taking the limit in \eqref{eq.w2} as $r\rightarrow 0$, 
we obtain $\langle w, F \rangle _{L^2 (\R^{d+1}_+)}=0$
for $F \in C^{\infty}_0(\R^{d+1}_+)$. This proves the uniqueness
of the solution to \eqref{eq.dirichlet}. The uniqueness of the solution to \eqref{eq.neumann} in $\dot{H}^1(\R^{d+1}_+)$ is straightforward since we can take $\varphi=v$ in \eqref{def.weak.neumann}. The proof is complete. 

\subsection{Solvability of 
the inhomogeneous problem}\label{subsec.solvability}

The purpose of this section 
is prove the solvability of 
the inhomogeneous boundary value problems 
under weak conditions for the inhomogeneous term.

\begin{thm}\label{thm.solvability} Assume that the semigroups $\{e^{-t\mathcal{P}_{\mathcal{A}}}\}_{t\geq 0}$ and  $\{e^{-t\mathcal{P}_{\mathcal{A}^*}}\}_{t\geq 0}$ acting on $H^{1/2}(\R^d)$ 
are extended as strongly continuous semigroups acting on 
$L^2 (\R^d)$, and that $\{e^{-t\mathcal{P}_{\mathcal{A}^*}}\}_{t\geq 0}$ is  bounded  in $L^2 (\R^d)$. Assume further that $D_{L^2}(\mathcal{P}_{\mathcal{A}}) = D_{L^2}(\mathcal{P}_{\mathcal{A}^*}) = H^1 (\R^d)$ with equivalent norms. Then for given $F \in L^1 (\R_+; L^2 (\R^d))$ and $g=0$ 
there exists a weak and mild solution $u$ to \eqref{eq.dirichlet} such that 
$$
u \in C(\overline{\R_+};L^2(\R^{d})) \quad \textit{and} \quad  
\nabla u  \in L^p_{loc}(\overline{\R_+};L^2(\R^d)) \ \ 
\textit{for any} \ \, p\in [1,\infty).
$$
If in addition $h =M_{b}  \int_0^\infty   e^{-s \mathcal{Q}_{\mathcal{A}}} M_{1/b} F (s) \dd s$ belongs to the range of $\Lambda_{\mathcal{A}}$, then there exists a weak and mild solution $v$ to \eqref{eq.neumann} such that 
$$
v \in C(\overline{\R_+};L^2(\R^{d})) \quad \textit{and} \quad 
\nabla v  \in L^p_{loc}(\overline{\R_+};L^2(\R^d))
\ \ \textit{for any} \ \, p \in [1,2).
$$

\end{thm}

\begin{rem}
{\rm Theorem \ref{thm.solvability.intro2} 
is a consequence of this result.
Indeed Proposition \ref{prop.poisson.regular.1} shows that 
under the assumptions in Theorem \ref{thm.solvability.intro2}, 
the semigroups $\{e^{-t\mathcal{P}_{\mathcal{A}}}\}_{t\geq 0}$ and 
$\{e^{-t\mathcal{P}_{\mathcal{A}^*}}\}_{t\geq 0}$ on $H^{1/2}(\R^d)$
can be extended as strongly continuous semigroups on 
$L^2 (\R^d)$,  and $\{e^{-t\mathcal{P}_{\mathcal{A}^*}}\}_{t\geq 0}$ 
becomes bounded in $L^2 (\R^d)$. 
}
\end{rem}

\smallskip

When the inhomogeneous 
term belongs to some homogeneous Sobolev space 
of negative order, 
we  have  an interpolating result between the Lax-Milgram theorem and Theorem \ref{thm.solvability}.
\begin{thm}
\label{thm.solvability.variant}
\ Suppose that the assumptions in Theorem \ref{thm.solvability} hold. Let $1<p<\infty$ and $1-1/p < r <1$.
Then for $F\in L^p (\R_+; \dot{H}^{-r}(\R^d))$ there exists a weak solution $u$ to 
\eqref{eq.dirichlet} satifying 
\begin{equation}
\label{class.solution}
u \in C(\overline{\R_+};L^2(\R^{d}))  \quad \textit{and} \quad   
\nabla u  \in L^q_{loc}(\overline{\R_+};L^2(\R^d)) ,~~~~~~\frac{1}{q}=r-1+\frac{1}{p}.
\end{equation}
If in addition 
$h =M_{b}  \int_0^\infty   e^{-s \mathcal{Q}_{\mathcal{A}}} M_{1/b} F  (s) \dd s$ belongs to the range of $\Lambda_{\mathcal{A}}$, 
then there exists a weak solution $v$ to \eqref{eq.neumann} satisfying
$$
v \in C(\overline{\R_+};L^2(\R^{d}))  \quad \textit{and} \quad   
\nabla v  \in L^{\tilde q}_{loc}(\overline{\R_+};L^2(\R^d)),~~~~~~1<\tilde q <\min\{q,2\}
$$
for the same $q$ as above.
\end{thm}

\begin{rem}{\rm It is not clear whether 
the weak solution in the class \eqref{class.solution}
is unique. It would be also an interesting question to find 
other conditions on $A$ for the boundedness of the semigroup 
$\{e^{-t\mathcal{P}_{\mathcal{A}^*}}\}_{t\geq 0}$ in $L^2(\R^d)$.
}
\end{rem}

Since Theorems \ref{thm.solvability}
and \ref{thm.solvability.variant} are proved 
in the same manner, we will only give 
a proof of the latter one.

\vspace{0.3cm}

\noindent
{\it Proof of Theorem \ref{thm.solvability.variant}.}
\ Consider the Dirichlet problem \eqref{eq.dirichlet}.
Let $\psi _N=N^{-d}\psi(\cdot /N)$ be a function such that 
$\psi \in C_0^{\infty}(\R^d)$ and let $\phi_N=\phi(t/N)$ 
be a smooth function in $\R_+$ such that $\phi(t)=1$ for $0\le t<1$ 
and $\phi(t)=0$ for $t\ge 2$. 
Then the function $F_N $ defined by 
$F_N(x,t)=\phi _N(t)(F_N(t)*\psi _N)(x)$
($N=1,2,\cdots$)  satisfies that 
$$
F_N \in \dot{H}^{-1}(\R^{d+1}_+)
\qquad {\rm and} \qquad 
F_N \rightarrow F \ \ {\rm in}\ \ L^p (\R_+; \dot{H}^r (\R^d)) 
\ \ {\rm as } \ \ N \rightarrow \infty.
$$
Therefore, for each $F_N$ ($N =1,2,\cdots$), 
the Lax-Milgram theorem gives the unique 
weak solution  $u_N$ in $\dot{H}^1_0(\R^{d+1}_+)$ such that
\begin{align}
\label{eq.u_N}
\langle A\nabla u_N, \nabla \varphi \rangle _{L^2 (\R^{d+1}_+)} = \langle F_N,\varphi \rangle _{L^2(\R^{d+1}_+)}~~~~~~~~{\rm for~all}~~\varphi\in C_0^\infty (\R^{d+1}_+)
\end{align}
and ~$u_N(t) \rightarrow  0$ as $t\rightarrow 0$ in the sense of 
distributions. By Theorem \ref{thm.representation} the weak 
solution is also the mild solution. 
Then the representation formula 
\eqref{eq.thm.representation.1} and  
the maximal regularity 
for  $\mathcal{P_A}$ yields
\begin{align}
\label{est.u_N}
& ~~~ \|\partial_t u_N \|_{L^q (0,T; L^2 (\R^d))} 
+
\| \mathcal{P_A}u_N\|_{L^q (0,T; L^2 (\R^d))}\nonumber \\
&\le C_T
\| 
\int_{\cdot}^\infty 
e^{-(\tau-\cdot)\mathcal{Q}_\mathcal{A}} M_{1/b} F_N (\tau ) 
\dd\tau 
\|_{L^q (0,T; L^2(\R^d))}.
\end{align}
Then the duality argument and  
boundedness of the semigroup $e^{-t\mathcal{P_{A^*}}}$
 in $L^2(\R^d)$ together with Remark 
\ref{rem.prop.poisson.analytic.L^2.1} yield
\begin{align*}
{\rm RHS \ of}\ \eqref{est.u_N}
&\le
C_T\|\int_0^\infty (\tau-\cdot)^{-r} \|F_N (\tau ) \|_{\dot{H}^{-r}(\R^d)}
\dd\tau \|_{L^q (\R_+)}\le C_T
\|F_N\|_{L^p (\R_+; \dot{H}^r (\R^d))}. 
\end{align*}
where we have used the Hardy-Littlewood-Sobolev inequality
in the last line.
From the assumption $D_{L^2}(\mathcal{P_A})=H^1(\R^d)$ 
with the equivalent norm, we also see 
$\nabla u_N$ converges to $\nabla u$ 
strongly in $L^q_{loc} (\overline{\R_+}; L^2(\R^d))$ as $N \rightarrow \infty$, 
where $u$ is the mild solution to  \eqref{eq.dirichlet} 
with the inhomogeneous term $F$.
Thus taking $N\rightarrow \infty$ in \eqref{est.u_N}, we see that
$u$ is the weak solution to \eqref{eq.dirichlet}.
The continuity in $L^2(\R^d)$ follows from the  interpolation between 
the estimates of $u$ and $\partial_t u$ in $L^q_{loc} (\overline{\R_+}; L^2(\R^d))$. 

For the Neumann problem \eqref{eq.neumann},  
by the assumption and 
$D_{L^2}(\Lambda_{\mathcal{A}}) \hookrightarrow H^{1/2}(\R^d)$, we have
$\Lambda_{\mathcal{A}}^{-1}h$ belongs to $H^{1/2}(\R^d)$. 
Hence from the analyticity of the semigroup and the interpolation, 
we have 
$\nabla e^{-t\mathcal{P_A}} \Lambda_{\mathcal{A}}^{-1}h \in 
L^s_{loc}(\overline{\R_+}; L^2(\R^d))$ for any $1<s<2$.  
The remaining part can be shown in the same way 
as in the Dirichlet problem. Thus the proof is complete.

\appendix

\section{Counterpart of Theorem \ref{thm.factorization.strong.1} \, for $\mathcal{A}+\lambda$}

Let $\lambda\in \C$ and we denote by $\mathcal{A}_\lambda$ the realization of $\mathcal{A} + \lambda$ in $L^2 (\R^{d+1})$.
\begin{df}\label{def.appendix} \noindent {\rm (i)} Let ${\rm Re }\,  \lambda > 0$. We denote by $E_{\mathcal{A}_\lambda}: H^{1/2} (\R^d)\rightarrow H^1 (\R^{d+1}_+)$ the $\mathcal{A}_\lambda$-extension operator, i.e., $u=E_{\mathcal{A}_\lambda} g$ is the solution to the Dirichlet problem
\begin{equation}\label{eq.dirichlet.appendix0}
\begin{cases}
& \mathcal{A}_\lambda u  = 0~~~~~~~{\rm in}~~~\R^{d+1}_+,\\
& ~ \hspace{0.3cm} u  = g~~~~~~~{\rm on}~~\partial\R^{d+1}_+
=\R^d\times \{ t=0\}.
\end{cases}
\end{equation}
The one parameter family of linear operators $\{E_{\mathcal{A}_\lambda} (t)\}_{t\geq 0}$, defined by $E_{\mathcal{A}_\lambda} (t) g = (E_{\mathcal{A}_\lambda} g )(\cdot,t)$ for $g\in H^{1/2}(\R^d)$, is called the Poisson semigroup associated with $\mathcal{A}_\lambda$.
 
\noindent {\rm (ii)} Let ${\rm Re }\, \lambda > 0$. We denote by $\Lambda_{\mathcal{A}_\lambda}: H^{1/2} (\R^d)\rightarrow H^{-1/2} (\R^d)$ the Dirichlet-Neumann map associated with $\mathcal{A}_\lambda$, which is defined through the sesquilinear form 
\begin{equation}
\langle  \Lambda_{\mathcal{A}_\lambda} g, \varphi \rangle_{H^{-\frac12},H^{\frac12}}  =  \langle  A\nabla E_{\mathcal{A}_\lambda}g, \nabla E_{\mathcal{A}_\lambda} \varphi \rangle_{L^2(\R^{d+1}_+)} +  \lambda  \langle   E_{\mathcal{A}_\lambda}g, E_{\mathcal{A}_\lambda} \varphi \rangle_{L^2(\R^{d+1}_+)},~~~~~~~~~g,\varphi \in \dot{H}^{\frac12}(\R^d).\label{def.Lambda.appendix}
\end{equation}
Here $\langle \cdot,\cdot\rangle _{H^{-1/2},H^{1/2}}$ denotes the duality coupling of  $H^{-1/2}(\R^d)$ and $H^{1/2}(\R^d)$.

\end{df}

Note that $\langle  \Lambda_{\mathcal{A}_\lambda} g, \varphi \rangle_{H^{-1/2},H^{1/2}}= \langle   g, \Lambda_{\mathcal{A}_\lambda^*}\varphi \rangle_{H^{1/2},H^{-1/2}} $ holds. The counterpart of Theorem \ref{thm.factorization.strong.1} with ${\rm Re}\,\lambda>0$ is stated as follows.
\begin{thm}\label{thm.factorization.strong.appendix}  
Let ${\rm Re}\,  \lambda > 0$. Suppose that the assumptions of Theorem \ref{thm.factorization.strong.1} hold. Then  $H^1 (\R^d)$ is continuously embedded in $D_{L^2}(\Lambda_{\mathcal{A}_\lambda})\cap D_{L^2}(\Lambda_{\mathcal{A}_\lambda^*})$, and the operators $-{\bf P}_{\mathcal{A}_\lambda}, ~-{\bf P}_{\mathcal{A}_\lambda^*}$ defined by  
\begin{align}
& D_{L^2} ({\bf P}_{\mathcal{A}_\lambda}) = H^1 (\R^d), ~~~~~~~~~ -{\bf P}_{\mathcal{A}_\lambda} f =  -M_{1/b} \Lambda_{\mathcal{A}_\lambda} f - M_{{\bf r_2}/b} \cdot \nabla _x f\,, \label{def.P.appendix.1}\\
& D_{L^2} ({\bf P}_{\mathcal{A}_\lambda^*}) = H^1 (\R^d), ~~~~~~~~~-{\bf P}_{\mathcal{A}_\lambda^*} f =  - M_{1/\bar{b}} \Lambda_{\mathcal{A}_\lambda^*} f - M_{{\bf \bar{r}_1}/\bar{b}} \cdot \nabla_x f\,,\label{def.P.appendix.2}
\end{align}
generate strongly continuous and analytic semigroups in $L^2 (\R^d)$. Moreover, the realization of $\mathcal{A}' + \lambda$ in $L^2 (\R^d)$ and the realization  $\mathcal{A}_\lambda$ in $L^2 (\R^{d+1})$  are respectively  factorized  as 
\begin{align}
\mathcal{A}' +\lambda & = M_b \mathcal{Q}_{\mathcal{A}_\lambda } {\bf P}_{\mathcal{A}_\lambda},~~~~~~~~~~~~\mathcal{Q}_{\mathcal{A}_\lambda} = M_{1/b} (M_{\bar{b}} {\bf P}_{\mathcal{A}_\lambda^*})^*, \label{eq.factorization.strong.lambda0}\\
\mathcal{A}_\lambda & = - M_b (\partial_t -  \mathcal{Q}_{\mathcal{A}_\lambda} ) (\partial_t + {\bf P}_{\mathcal{A}_\lambda} ).\label{eq.factorization.strong.lambda1}
\end{align}
Here $(M_{\bar{b}} {\bf P}_{\mathcal{A}_\lambda^*})^*$ is the adjoint of $M_{\bar{b}} {\bf P}_{\mathcal{A}_\lambda^*}$ in $L^2 (\R^d)$.

\end{thm}

\begin{rem}{\rm  Theorem \ref{thm.factorization.strong.appendix} is proved in the same manner as in the proof of Theorem \ref{thm.factorization.strong.1}. We also note that the conclusion of Theorem \ref{thm.factorization.strong.appendix} holds if $A$ is Lipschitz continuous, which can be proved by arguing as in \cite{MM1}.

}
\end{rem}

\end{document}